\documentclass[preprint,10pt]{elsarticle}

\makeatletter
\def\ps@pprintTitle{%
   \let\@oddhead\@empty
   \let\@evenhead\@empty
   \def\@oddfoot{\reset@font\hfil\thepage\hfil}
   \let\@evenfoot\@oddfoot
}
\makeatother

\usepackage{amsmath,amsfonts,mathtools}
\usepackage[margin=0.75in]{geometry}
\usepackage{booktabs}
\usepackage{subcaption}
\usepackage[section]{placeins}
\usepackage{hyperref}
\usepackage{multirow}
\hypersetup{colorlinks=true,urlcolor=black,linkcolor=black}
\usepackage[T1]{fontenc}\usepackage{lmodern}
\usepackage{array}
\usepackage{color}
\usepackage{xfrac}
\usepackage{multirow}
\usepackage{enumitem}
\usepackage{tikz}
\usepackage{pifont}
\usepackage{textcomp}
\usepackage{algorithm,algorithmic}

\newcommand{\xmark}{\ding{55}}
\usetikzlibrary{arrows}

\newcommand{\titlepaper}{Machine-learning error models for approximate solutions to\\ parameterized systems of
nonlinear equations}
\newcommand{\eqnref}[1]{Eq.~\eqref{#1}}

\definecolor{orange}{rgb}{1,0.5,0}
\definecolor{purple}{rgb}{0.5,0,0.5}
\definecolor{green}{rgb}{0,0.5,0}
\definecolor{todo}{gray}{0.75}

\newcommand{\RR}[1]{\mathbb{R}^{#1}}
\newcommand{\RRstar}[1]{\mathbb{R}_\star^{#1}}
\newcommand{\RRplus}[1]{\mathbb{R}_+^{#1}}
\newcommand{\defeq}{\vcentcolon =}
\newcommand{\zero}{\mathbf{0}}
\newcommand{\reviewerOne}[1]{#1} 
\newcommand{\reviewerTwo}[1]{#1} 
\newcommand{\reviewerThree}[1]{#1} 
\newcommand{\ourRereading}[1]{#1} 

\newcommand{\erfinv}{\text{erf}^{-1}}

\newcommand{\pti}{parameter training instances}

\newcommand{\state}{\mathbf{u}}
\newcommand{\stateOne}{\mathbf{z}_1}
\newcommand{\stateTwo}{\mathbf{z}_2}
\newcommand{\stateDummy}{\mathbf{v}}
\newcommand{\nstate}{N_{\state}}
\newcommand{\stateIt}[1]{\state^{(#1)}}
\newcommand{\stateApprox}{\tilde \state}
\newcommand{\stateApproxArg}[1]{\stateApprox^{#1}}
\newcommand{\stateRef}{\bar \state}
\newcommand{\stateLF}{\state_\text{LF}}
\newcommand{\nstateLF}{N_{\stateLF}}
\newcommand{\nfeatures}{N_{\features}}
\newcommand{\nfeaturesSplit}{N_{\text{split}}}
\newcommand{\SVRfeaturespace}{\mathcal F}
\newcommand{\innerprod}[2]{\langle #1,#2\rangle}
\newcommand{\itConverge}{K}
\newcommand{\bigO}{\mathcal{O}}
\newcommand{\nearestSetArg}[1]{\mathcal{I}(#1)}

\newcommand{\dual}{\mathbf{y}}
\newcommand{\dualEntry}[1]{y_{#1}}

\newcommand{\paramDomainPOD}{\paramDomain_\mathrm{POD}}

\newcommand{\prolongateNo}{\mathbf{p}}
\newcommand{\prolongateArg}[1]{\mathbf{p}(#1)}

\newcommand{\qoi}{s}
\newcommand{\qoiApprox}{\tilde s}

\newcommand{\qoiFunc}{g}

\newcommand{\res}{\mathbf{r}}

\newcommand{\resApprox}{{\mathbf{r}}}
\newcommand{\resApproxArg}[1]{{\mathbf{r}^{#1}}}
\newcommand{\resApproxEntry}[1]{ r_{#1}}
\newcommand{\resLF}{{\mathbf{r}}_\text{LF}}
\newcommand{\resNormNo}{\ensuremath{\|\resApprox\|}}
\newcommand{\resNorm}[1]{\ensuremath{\|\resApprox(#1)\|}}

\newcommand{\jacobian}{\mathbf{ J}}

\newcommand{\lipschitzConstant}{\beta}
\newcommand{\lipschitzConstantUB}{\lipschitzConstant_\mathrm{UB}}
\newcommand{\invlipschitzConstant}{\alpha}
\newcommand{\invlipschitzConstantLB}{\invlipschitzConstant_\mathrm{LB}}
\newcommand{\lipschitzConstantQoiFunc}{\beta_\qoiFunc}
\newcommand{\lipschitzConstantQoiFuncUB}{\beta_{\qoiFunc,\mathrm{UB}}}

\newcommand{\params}{\boldsymbol{\mu}}
\newcommand{\paramArg}[1]{{\mu}_{#1}}
\newcommand{\paramsDummy}{\boldsymbol{\nu}}
\newcommand{\nparams}{N_{\params}}

\newcommand{\napprox}{{N_\mathrm{approx}}}
\newcommand{\range}[1]{\mathrm{Ran}(#1)}
\newcommand{\card}[1]{\mathrm{card}\left(#1\right)}
\newcommand{\paramDomain}{\mathcal D}
\newcommand{\paramDomainDummy}{\bar\paramDomain}
\newcommand{\trainingData}{{\mathcal T_\mathrm{train}}}
\newcommand{\trainingDataArg}[1]{{\mathcal T_\mathrm{train}^{#1}}}
\newcommand{\hyperparams}{{\boldsymbol \theta}}
\newcommand{\hyperparamsSelect}{\hyperparams^\star}
\newcommand{\hyperparamsSet}{{\boldsymbol \Theta}}
\newcommand{\testData}{{\mathcal T_\mathrm{test}}}
\newcommand{\testDataTwo}{{\mathcal T_{\overline{\mathrm{test}}}}}
\newcommand{\errorTestNo}{{\error_{0}}}
\newcommand{\errorTest}[1]{{\error_{0,#1}}}

\newcommand{\errorTestTwo}[1]{{\error_{\bar{0},#1}}}
\newcommand{\featuresTestTwo}[1]{{\features_{\bar{0},#1}}}

\newcommand{\featuresTest}[1]{{\features_{0,#1}}}
\newcommand{\errorTrain}[1]{{\error_{#1}}}
\newcommand{\featuresTrain}[1]{{\features_{#1}}}

\newcommand{\errorDataArg}[1]{{\mathcal T_\error^{#1}}}

\newcommand{\residualDataArg}[1]{{\mathcal T_\resApprox^{#1}}}

\newcommand{\ntrainingData}{{N_\mathrm{train}}}
\newcommand{\ntestData}{N_\mathrm{test}}

\newcommand{\trainingDataRes}{\mathcal T_{\mathrm{train},\resApprox}}
\newcommand{\ntrainingDataRes}{N_{\mathrm{train},\resApprox}}

\newcommand{\paramDomainTrain}{\paramDomain_\mathrm{train}}
\newcommand{\paramDomainTest}{\paramDomain_\mathrm{test}}

\newcommand{\orthoMat}[2]{\mathbb V_{#2}(\RR{#1})}

\newcommand{\stateRed}{\hat\state}
\newcommand{\trialbasis}{\mathbf{\Phi}_\mathbf{u}}
\newcommand{\testbasis}{\mathbf{\Psi}_\mathbf{u}}
\newcommand{\basisresApprox}{\mathbf{\Phi}_{\mathbf r}}
\newcommand{\dimresApprox}{m_\resApprox}
\newcommand{\sampleMat}{\mathbf{P}}
\newcommand{\nsamples}{n_{\resApprox}}
\newcommand{\sampleMatFull}{\sampleMat_\star}
\newcommand{\resApproxRed}{\hat{\resApprox}}
\newcommand{\resApproxRedGappy}{\resApproxRed_\mathrm{g}}
\newcommand{\resApproxRef}{\bar{\resApprox}}

\newcommand{\tol}{\epsilon}

\newcommand{\error}{\delta}

\newcommand{\errorArg}[1]{\error^{#1}}
\newcommand{\errorArgs}[2]{\error_{#1,#2}}
\newcommand{\errorSubArg}[1]{\error_{#1}}
\newcommand{\errorModel}{\hat \error}
\newcommand{\errorQoi}{\error_{\qoi}}
\newcommand{\errorQoiArg}[1]{\error_{{#1}}}
\newcommand{\errorState}{\errorSubArg{\state}}

\newcommand{\errorStateVec}{\mathbf{e}}
\newcommand{\dualweightedresidual}{d}
\newcommand{\dualweightedresidualArg}[1]{\dualweightedresidual_{#1}}

\newcommand{\identity}{\mathbf{I}}
\newcommand{\identityArg}[1]{\identity_{#1}}

\newcommand{\trueRegressionModel}{f}
\newcommand{\regressionModel}{\hat f}
\newcommand{\regressionModelArg}[1]{\regressionModel_{#1}}
\newcommand{\regressionModelArgs}[2]{\regressionModel_{#1,#2}}
\newcommand{\validationLossArgs}[2]{L_{#1,#2}}
\newcommand{\ntree}{{N_\text{tree}}}

\newcommand{\features}{\mathbf{x}}
\newcommand{\featuresArg}[1]{\features^{#1}}
\newcommand{\featuresArgs}[2]{\features_{#1,#2}}
\newcommand{\featureArg}[1]{x_{#1}}
\newcommand{\noise}{\epsilon}
\newcommand{\noiseModel}{\hat\noise}
\newcommand{\expectation}[1]{\mathrm{E}[#1]}
\newcommand{\variance}[1]{\mathrm{Var}[#1]}
\newcommand{\normal}[2]{\mathcal{N}(#1,#2)}
\newcommand{\varianceEstimate}{\hat\sigma^2}
\newcommand{\expectationGiven}[2]{\expectation{#1\, |\, #2}}
\newcommand{\frequency}{\omega}

\newcommand{\weight}{\mathbf{w}}
\newcommand{\weightDummy}{\bar\weight}
\newcommand{\weightArg}[1]{w_{#1}}
\newcommand{\knnWeight}{\tau}
\newcommand{\knnWeightUniform}{\knnWeight_\text{u}}
\newcommand{\knnWeightEuclidean}{\knnWeight_2}
\newcommand{\nnFuncArg}[1]{\mathbf{g}_{#1}}
\newcommand{\nlayers}{N_\text{layers}}
\newcommand{\weightNNarg}[1]{\mathbf{W}_{#1}}

\newcommand{\nNeuron}[1]{p_{{#1}}}

\newcommand{\fcharacteristic}{Feature attribute}

\newcommand{\FeatureResNorm}{\|\mathbf{r}\|}
\newcommand{\FeatureParResSampled}{[\params;\,\sampleMat\mathbf{r}]}
\newcommand{\FeatureParResGappyPCA}{[\params;\,\hat{\mathbf{r}}_\text{g}]}
\newcommand{\FeatureParResPCA}{[\params;\,\hat{\mathbf{r}}]}
\newcommand{\FeaturePar}{\params}
\newcommand{\FeatureParResNorm}{[\params;\,\FeatureResNorm]}

\newcommand{\nsnapshots}{M_\mathbf{u}}
\newcommand{\nbasis}{m_\mathbf{u}}
\newcommand{\nbasisArg}[1]{\nbasis^{#1}}

\newcommand{\leftSing}{\bar{\mathbf{U}}}
\newcommand{\leftSingArg}[1]{\bar{\mathbf{u}}_{#1}}
\newcommand{\Sing}{\bar{\boldsymbol{\Sigma}}}
\newcommand{\SingArg}[1]{\bar{{\sigma}}_{#1}}
\newcommand{\rightSing}{\bar{\mathbf{V}}}
\newcommand{\rightSingArg}[1]{\bar{\mathbf{v}}_{#1}}

\newcommand{\energyThresholdArg}[1]{\upsilon_{#1}}
\newcommand{\diag}{\mathrm{diag}}
\newcommand{\vectorize}[1]{\mathrm{vec}(#1)}
\newcommand{\activationfunction}{h}

\newlist{Objective}{enumerate}{2}
\setlist[Objective]{label={Objective \arabic*},align=left}

\newlist{Step}{enumerate}{2}
\setlist[Step]{label={Step \arabic*},align=left}

\newlist{Attribute}{enumerate}{2}
\setlist[Attribute]{label={Attribute \arabic*},align=left}

\newlist{Data}{enumerate}{2}
\setlist[Data]{label={Data Set \arabic*},align=left}

\begin{document}

\begin{frontmatter}



\title{\titlepaper}


\author[freno]{Brian A.\ Freno}
\ead{bafreno@sandia.gov}
\author[carlberg]{Kevin T.\ Carlberg}
\ead{ktcarlb@sandia.gov}
\address[freno]{Sandia National Laboratories, PO Box 5800, MS 0828, Albuquerque, NM 87185}
\address[carlberg]{Sandia National Laboratories, 7011 East Ave, MS 9159, Livermore, CA 94550}

\begin{abstract}
	This work proposes a machine-learning framework for constructing statistical
	models of errors incurred by approximate solutions to parameterized systems
	of nonlinear equations. These approximate solutions may arise from early
	termination of an iterative method, a lower-fidelity model, or a
	projection-based reduced-order model, for example.  The proposed statistical
	model comprises the sum of a deterministic regression-function model and a
	stochastic noise model.  The method constructs the regression-function
	model by applying regression techniques from machine learning (e.g., support
	vector regression, artificial neural networks) to map features (i.e., error
	indicators such as sampled elements of the residual) to a prediction of the
	approximate-solution error. The method constructs the noise model as a
	mean-zero Gaussian random variable whose variance is computed as the sample
	variance of the approximate-solution error on a test set; this variance can
	be interpreted as the epistemic uncertainty introduced by the approximate
	solution.  This work considers a wide range of feature-engineering methods,
	data-set-construction techniques, and regression techniques that aim to
	ensure that (1) the features are cheaply computable, (2) the noise model
	exhibits low variance (i.e., low epistemic uncertainty introduced), and (3)
	the regression model generalizes to independent test data. Numerical
	experiments performed on several computational-mechanics problems and types
	of approximate solutions demonstrate the ability of the method to generate
	statistical models of the error that satisfy these criteria and
	significantly outperform more commonly adopted approaches for error
	modeling.
\end{abstract}


\end{frontmatter}


\section{Introduction} 
\label{sec:introduction}

Myriad decision-making applications in science and engineering are
\textit{many-query} in nature, as they require a parameterized computational
model to be simulated for a large number of parameter instances; examples
include \textit{design optimization}, where each parameter instance
corresponds to a different candidate system design; \textit{uncertainty
propagation}, where each parameter instance corresponds to a realization of a
random variable; and \textit{Bayesian inference}, where each parameter
instance corresponds to a sample from a posterior distribution.  Oftentimes,
the computational model in these contexts corresponds to a large-scale system
of parameterized nonlinear algebraic equations; this occurs, for example, when
the model associates with the spatial discretization of a system of
parameterized, nonlinear, stationary partial differential equations.

In such cases, it is prohibitively expensive to compute an \textit{exact
solution} to the system of nonlinear equations for each query.  Instead, 
computationally inexpensive \textit{approximate solutions} are often employed
for tractability; examples include \textit{inexact solutions} computed from
early termination of an iterative method (e.g., Newton's method); solutions
computed from \textit{lower-fidelity models} (e.g., a coarse mesh, lower
finite-element order); and solutions computed from projection-based
\textit{reduced-order models (ROMs)} (e.g., proper orthogonal decomposition
with Galerkin projection).  However, such approximations introduce an error
that should be quantified and accounted for in the resulting analysis. In the
context of uncertainty quantification, this error can be interpreted as a
source of epistemic uncertainty; as such, it is natural to quantify the error
in a statistical manner.
Researchers have developed three strategies to quantify the error in such
approximate solutions: (1) error indicators, (2) rigorous \textit{a
posteriori} error bounds, and (3)
error models. 

\textit{Error indicators} are computable quantities that are
informative of the approximate-solution error; they do not rigorously
bound the error, nor do they generally produce an unbiased prediction of the
error or an estimate of the variance in the prediction. One example is the
(dual) norm of the residual evaluated at the approximate solution. This
quantity is informative of the error, as it appears as a term in
residual-based \textit{a posteriori} error bounds. As such, it is often used
as a termination criterion for iterative linear and nonlinear solvers, as well
as for greedy methods that determine snapshot-collection parameter instances
during (offline) ROM construction
\cite{bui2008parametric,bui2008model,hinze2012residual,amsallem2015design,
wu2015adaptive,yano2018lp} and when deploying ROMs within a trust-region
setting \cite{zahr2015progressive,zahrThesis}. Unfortunately, computing the
residual norm is typically computationally costly, as its evaluation incurs an
operation count that scales with the dimension of the original problem unless
the residual is linear in the solution and affine in functions of the
parameters.  Another commonly used error indicator is the dual-weighted
residual (i.e., adjoint-based error estimator), which provides a
\reviewerOne{first-order} approximation of the error in a quantity of interest.
This approach is often employed for goal-oriented error estimation (and
adaptive refinement) for finite-element
\cite{babuvska1984post,becker1996weighted,rannacher1999dual,bangerth2003adaptive},
finite-volume \cite{venditti2000adjoint,venditti2002grid,park2004adjoint}, and
discontinuous-Galerkin discretizations
\cite{lu2005posteriori,fidkowski2007simplex}, as well as for model reduction
\cite{meyer2003efficient,carlberg2014adaptive}. Unfortunately, dual-weighted
residuals are often difficult to implement, as they require solving a dual
linear system whose matrix is the transpose of the residual Jacobian; this is
not always possible, for example, when the Jacobian is available only as a
black-box operator.  Furthermore, computing the dual-weighted residual is
computationally costly, as the dimension of the dual linear system is the same
as that of the original nonlinear system; this cost can be mitigated if the
dual linear solve is approximated using, for example, a lower-fidelity model
or reduced-order model. A recently developed approach computes randomized
approximations of normed solution errors for parameterized linear systems,
wherein the error estimate is indeed unbiased and yields probabilistic error
bounds \cite{smetana2018randomized}.

\textit{Rigorous \textit{a posteriori} error bounds} are quantities that bound
either the normed solution error or the quantity-of-interest error arising
from an approximate solution. These error bounds are typically residual based
(i.e., the (dual) norm of the residual appears as a term in the bound) and
were pioneered in the reduced-basis community \cite{GP05,rozza2007reduced}.
Although these quantities rigorously bound the error, they are typically
not sharp, often overpredicting the error by orders of magnitude. Additionally,
they can be challenging to implement, as they require Lipschitz or inf--sup
constants to be estimated or bounded, which incurs additional computational
cost and can further degrade the bounds' sharpness \cite{SCM,wirtzDeim}.
In the reduced-basis context, a recently proposed hierarchical error estimator
can yield sharper estimates without the need to compute these constants, at the
cost of solving a higher-dimensional reduced-order model
\cite{hain2018hierarchical}.  Finally, these (deterministic) bounds are of
limited utility in an uncertainty-quantification setting, where a statistical
model of the error is more readily integrable into the uncertainty analysis.

\textit{Error models} directly model the error incurred by an
approximate solution. In contrast to error indicators and rigorous error
bounds, error models often yield a more accurate approximation of the error;
furthermore, they facilitate integration of approximate solutions into uncertainty
quantification when the error model is statistical in nature. The most common
error-modeling approach is to directly model the mapping from model parameters
to the error in a quantity of interest.  This is an \textit{a priori} error
model, as it does not leverage any data generated by the approximate solution
(e.g., the residual norm)
at new parameter instances. This approach was developed
independently by the optimization community and the model-calibration
community.  In the optimization community, the technique is often referred to
as a `multifidelity correction', wherein the modeled error corresponds to the
error between low- and high-fidelity models. Here, the error model is
constructed to enforce either `global' zeroth-order consistency between the
corrected low-fidelity-model prediction and the high-fidelity-model prediction
at training points
\cite{gano2005hybrid,huang2006sequential,march2012provably,NE12},
or `local' first- or second-order consistency at trust-region centers
\cite{Alexandrov2001, Eldred2004}.  In the model-calibration community, the
technique is referred to as computing a `model inadequacy' or `model discrepancy'
function, and the modeled error corresponds to the quantity-of-interest error
between a computational model and an underlying `truth' process that can be
experimentally measured
\cite{kennedy2001bayesian,higdon2003,higdon2004combining}.  Such approaches
tend to work well when the error exhibits a lower variance than the
high-fidelity or `truth' response \cite{NE12} and when the number of model
parameters is relatively small. Variants of this approach have been pursued in
the model-reduction community: Ref.~\cite{pagani2017efficient} interpolates time-dependent error
models in the parameter space in the case of dynamical-system models, and
Ref.~\cite{moosavi2018multivariate} constructs a
mapping from (a) the parameters, (b) the ROM-training parameters, and (c) the
ROM dimension to a prediction of the ROM error. \reviewerOne{The latter error model can be
employed for generating a decomposition of the parametric domain wherein each subdomain is equipped with a tailored
local basis and reduced-order model \cite{stefanescu2018parametric}.}
We note that these also comprise \textit{a priori} error models.

More recently, researchers have leveraged methods
from machine learning to construct more accurate \textit{a posteriori} error
models that leverage data generated by the approximate solution at new
parameter instances.  First, the
reduced-order-model error surrogates (ROMES) method \cite{drohmann_2015} was
\reviewerOne{an \textit{a posteriori} error model} developed in the context of model reduction. Based on the observation that the
aforementioned error indicators and error bounds can be computed from the
approximate solution and typically are both \textit{lower
dimensional} and \textit{more informative} of the error than model parameters, they
can be viewed as better \textit{features} for performing regression. As such, the
ROMES method applies kernel-based Gaussian-process regression \cite{RW06} to
construct a mapping from these features (i.e., error indicators
or error bounds) to a (normal or log-normal) random variable for the
error. The variance of this random variable can be interpreted as
the epistemic uncertainty introduced by the approximate solution.
This work demonstrated significant improvements in accuracy with respect to
the \textit{a priori} multifidelity correction error model
in high-dimensional parameter
spaces.  Follow-on work also demonstrated the promise of this kind 
of \textit{a posteriori} error model
in an uncertainty-quantification context \cite{manzoni2014accurate}.  

Despite these promising results, the ROMES method suffers from several
drawbacks.  First, due to the poor scalability of kernel-based Gaussian-process
regression, the method requires the user to hand-select a small number of
features.  Second, the cost of computing some of the features proposed by the ROMES technique is generally non-negligible; for example, the cost of evaluating the
residual norm generally scales with the dimension of the original model. Third, the
derivation and application of the ROMES method was limited to
reduced-order-model approximate solutions; this limitation is particularly
apparent in the application of model reduction to reduce the cost of computing
dual-weighted-residual features.

To overcome some of these challenges and extend the technique to dynamical
systems, Ref.~\cite{trehan2017error} developed a machine-learning
framework for modeling surrogate-model errors in the context of dynamical
systems.  Rather than applying Gaussian-process regression, that technique
applies high-dimensional regression methods (e.g., random forests) to map a
much larger set of candidate features at a given time instance to a prediction
of the time-instantaneous surrogate-model error. That approach also provides a
mechanism for constructing local-error models that are tailored to particular
feature-space regions. Numerical results on a subsurface-flow model with a
reduced-order-model approximate solution demonstrated the ability of the
technique to significantly improve the time-instantaneous quantity-of-interest
prediction, as well as to accurately model time-averaged errors.

This work aims to build upon the progress of Ref.~\cite{trehan2017error} by
proposing a machine-learning framework for modeling approximate-solution
errors in a different context: parameterized systems of nonlinear equations.
This method is characterized by three steps: (1) \textit{feature engineering},
which aims to devise features that are cheaply computable, informative of the
error, and low-dimensional; (2) \textit{regression-function modeling}, which
applies regression methods from machine learning (e.g., support vector
regression, artificial neural networks) to construct a mapping from these
features to a (deterministic) prediction of the approximate-solution error;
and (3) \textit{noise modeling}, which models the epistemic uncertainty in
the prediction as additive mean-zero Gaussian noise whose variance is computed
as the sample variance of the approximate-solution error on a test set.  These
steps are performed to realize three objectives: (1)~the features should be
cheaply computable, (2) the noise-model variance should be low 
(i.e., low epistemic uncertainty), and (3) the error model should generalize
to independent test data.
These error models can be applied to statistically model both
quantity-of-interest errors and normed solution errors. Thus, primary new
contributions of this work include:
 \begin{enumerate} 
			\item A new machine-learning framework for statistical modeling
				of approximate-solution errors (Section
				\ref{subsec:overview}),
			\item Novel residual-based feature-engineering methods (Methods
				\ref{feat:paramres}--\ref{feat:paramsampleres} in Section
				\ref{subsec:indicators}) that trade off two key attributes: the number
				of features and the quality of features,
			\item Two methods for constructing training and test data sets when
				multiple types of approximate solutions (e.g.,
				reduced-order models of multiple dimensions) are considered (Section
				\ref{subsection:trainingdata}),
			\item Numerical experiments performed across a wide range of problems
				and types of approximate solutions, which systematically compare
				eight candidate feature-engineering methods, seven regression methods,
				and two data-set-construction methods (Section \ref{sec:experiments}).
	 \end{enumerate}
 Notably, one of the newly proposed feature-engineering methods (i.e., gappy
 principal components of the residual) significantly outperforms the standard
 features of (1) the model parameters (as employed by the multifidelity
 correction or model-discrepancy
 method), and (2) the residual norm (as employed by the ROMES method); furthermore, this
 feature-engineering method is very inexpensive to compute, as it requires
 evaluating only $\mathcal O(10)$ elements of the residual and does not
 require solving a dual linear system (as is demanded by dual-weighted
 residuals). In addition, the best-performing regression method is typically
 an artificial neural network or support vector regression with a
 radial-basis-function kernel.  In all cases, the proposed methodology
 accurately predicts the approximate-solution errors with coefficients of
 determination on an independent test set exceeding $r^2=0.996$.
 
This paper is organized as follows.  
Section~\ref{sec:pnae} discusses parameterized systems of nonlinear equations,
including approximate solutions and approaches typically employed to
quantify the approximate-solution errors.
Section~\ref{sec:approach} presents the proposed approach for constructing
statistical models of approximate-solution errors, in particular the
proposed machine-learning framework, feature-engineering techniques, methods for
constructing training and testing data sets, and
regression-function approximation.
Section~\ref{sec:experiments} performs extensive numerical experiments that
explore the tradeoffs of different feature-engineering methods, regression
techniques, and data-set methods across several computational-mechanics
problems and types of approximate solutions.
Section~\ref{sec:conclusions} provides conclusions and an outlook for future
work.

In the remainder of this paper,  matrices are denoted by capitalized bold
letters, vectors by lowercase bold letters, and scalars by unbolded letters.
Additionally, we denote 
the elements of a vector as $\mathbf{a}\equiv\left[a_{1}\ \cdots\
a_{m}\right]^T$ and
the vertical concatenation of two vectors $\mathbf{a}\in\RR{n}$
and $\mathbf{b}\in\RR{m}$ as $
[
\mathbf{a}^T\ \mathbf{b}^T
]^T
\equiv
[\mathbf{a};\ \mathbf{b}]
\in\RR{n+m}$.



\section{Parameterized systems of nonlinear equations} 
\label{sec:pnae}

This work considers parameterized systems of nonlinear equations of the form
\begin{align}
	\res(\state(\params);\params) = \zero,
\label{eq:fom_eqn}
\end{align}
where $\res:\RR{\nstate}\times\RR{\nparams}\rightarrow\RR{\nstate}$ with
$\res:(\stateDummy;\paramsDummy)\mapsto\res$ denotes the residual, which is
nonlinear in at least its first argument; $\params\in\paramDomain$ denotes the
parameters with the parameter domain $\paramDomain\subseteq\RR{\nparams}$; and
$\state:\RR{\nparams}\rightarrow\RR{\nstate}$ denotes the state (i.e.,
solution vector) implicitly defined as the solution to \eqnref{eq:fom_eqn},
given the parameters.  In many scenarios, computing a scalar-valued quantity
of interest
\begin{equation} 
	\qoi(\params) \defeq\qoiFunc(\state(\params))
\label{eq:qoi}
\end{equation} 
is the ultimate objective of the analysis. Here,
$\qoi:\RR{\nparams}\rightarrow\RR{}$ and $\qoiFunc:\RR{\nstate}\rightarrow\RR{}$ with $\qoiFunc:\stateDummy\mapsto\qoiFunc$.
We consider Eq.~\eqref{eq:fom_eqn} to be the high-fidelity model, with
$\state$ the corresponding high-fidelity solution.

Many-query problems are characterized by the need to compute $\qoi(\params)$
for a large number of parameter instances. This is typically executed by first
computing the high-fidelity solution $\state(\params)$ by solving
\eqnref{eq:fom_eqn} and subsequently computing $\qoiFunc(\state(\params))$ for
each parameter instance.  In many cases, this approach is prohibitively
expensive (e.g., when $\nstate$ is large) or unnecessary (e.g., convergence in 
PDE-constrained optimization can be guaranteed by computing inexact solutions
at intermediate iterations \cite{heinkenschloss2002analysis}). Instead, such cases
demand the computation of an inexpensive approximate solution
for computational tractability. 

\subsection{Approximate solutions} 
\label{subsec:solution_approximations}

Although an approximate solution
$\stateApprox(\params)(\approx\state(\params))$ 
with
$\stateApprox:\RR{\nparams}\rightarrow\RR{\nstate}$
is generally computationally
inexpensive to compute, it leads to an approximation of the quantity of
interest
\begin{equation*}
	\qoiApprox(\params)\defeq\qoiFunc(\stateApprox(\params)),
\end{equation*}
where $\qoiApprox:\RR{\nparams}\rightarrow\RR{}$, which incurs an error 
$	\errorQoi (\params)\defeq \qoi(\params) - \qoiApprox(\params)$ that is often
non-negligible.  
This work considers three common approaches for computing approximate
solutions: (1) employing an inexact convergence tolerance or maximum iteration
count when iteratively solving \eqnref{eq:fom_eqn}, (2) employing a
lower-fidelity model and prolongating the solution to the higher-fidelity discretization
characterizing \eqnref{eq:fom_eqn}, and (3) applying projection-based model reduction.

\subsubsection{Inexact solutions} 
\label{subsec:inexact}

Applying an iterative method (e.g., Newton's method) to solve
\eqnref{eq:fom_eqn} yields a sequence of approximations $\stateIt{k}$,
$k=0,\ldots, K$. In this context, we can obtain an approximate solution after
iteration $K$ as $\stateApprox(\params) = \stateIt{K}$, where $\stateIt{K}$ is
considered to be an `inexact solution'. The maximum iteration count $K$ can be determined
either by explicitly specifying its value (e.g., $\itConverge=2$) or by ensuring the
inexact solution satisfies an inexact tolerance $\tol>0$ (e.g., $\tol = 0.1$)
such that $\|\res(\stateIt{K};\params)\|/\|\res(\zero;\params)\| \leq \tol$.

\subsubsection{Lower-fidelity model} 
\label{subsec:lowerfidelity}
Another approach for computing an approximate solution entails employing a
computational model that exhibits lower fidelity than the original model; this
lower fidelity can be achieved by neglecting physics, by coarsening the mesh
(when the governing equations \eqref{eq:fom_eqn} correspond to the spatial discretization of
a stationary partial-differential-equations problem), or by using finite
elements with a lower polynomial order, for example. We can characterize the resulting lower-fidelity model
by a (lower-dimensional) parameterized system of nonlinear 
equations:
\begin{align}
	\resLF(\stateLF(\params);\params) = \zero,
\label{eq:LF_eqn}
\end{align}
where $\resLF:\RR{\nstateLF}\times\RR{\nparams}\rightarrow\RR{\nstateLF}$ denotes the lower-fidelity residual; $\stateLF:\RR{\nparams}\rightarrow\RR{\nstateLF}$ denotes the lower-fidelity state implicitly defined by the solution to \eqnref{eq:LF_eqn}, given the parameters; and $\nstateLF < \nstate$ is the dimension of the lower-fidelity model. 
In order to represent the lower-fidelity state in the original
$\nstate$-dimensional state space, we apply a prolongation (or interpolation)
operator $\prolongateNo:\RR{\nstateLF}\rightarrow\RR{\nstate}$, such that the
approximate solution is $\stateApprox = \prolongateArg{\stateLF}$.

\subsubsection{Model reduction} 
\label{section:ROM}

Model reduction aims to reduce the computational cost of 
solving  \eqnref{eq:fom_eqn} by applying a (Petrov--)Galerkin projection
process. First, such approaches \ourRereading{seek an} approximate solution $\stateApprox$ in an $\nbasis$-dimensional affine trial subspace $\range{\trialbasis}+ \stateRef\subseteq\RR{\nstate}$, \reviewerOne{where $\range{\mathbf{A}}$ denotes the range of matrix $\mathbf{A}$ and} $\nbasis\ll\nstate$, i.e., 
\begin{align} 
\stateApprox(\params) = \trialbasis\stateRed(\params)+ \stateRef,
\label{eq:ROMsol}
\end{align} 
where $\trialbasis\in\RRstar{\nstate\times\nbasis}$ denotes the trial-basis
matrix with $\RRstar{m\times n}$ the set of full-column-rank $m\times n$
real-valued matrices (i.e., the non-compact Stiefel manifold),
$\stateRed:\RR{\nparams}\rightarrow\RR{\nbasis}$ denotes the generalized
coordinates of the approximate solution, and  $\stateRef\in\RR{\nstate}$
denotes a prescribed reference state (e.g., the mean value of snapshots in the
case of proper orthogonal decomposition). The trial basis $\trialbasis$ can be
computed using a variety of methods, for example, proper orthogonal decomposition
(POD) \cite{sirovich_1987,holmes_1996,freno_2014} (see Algorithm~\ref{alg:pod} of \ref{app:pod}), the reduced-basis method
\cite{rozza2007reduced}, and variants that employ gradient information
\cite{noor1980rbt,ito:Hermite,carlbergCpodJour}.

Substituting the approximate solution \eqnref{eq:ROMsol} into
\eqnref{eq:fom_eqn} yields an overdetermined system of $\nstate$ equations with
$\nbasis\ll \nstate$ unknowns: $\res(\trialbasis\stateRed(\params)+
\stateRef;\params) = \zero$, which may not have a solution. Thus, the second
step of model reduction enforces the orthogonality of the residual to an
$\nbasis$-dimensional test subspace $\range{\testbasis}\subseteq\RR{\nstate}$,
and the generalized coordinates $\stateRed(\params)$ are implicitly defined as the solution to 
\begin{align} 
\testbasis^T\res(\trialbasis\stateRed(\params)+ \stateRef;\params) = \zero,
\label{eq:rom_eqn}
\end{align} 
where $\testbasis\in\RRstar{\nstate\times\nbasis}$ denotes the test-basis matrix. Common choices for the test basis include Galerkin projection (i.e., $\testbasis = \trialbasis$) and least-squares Petrov--Galerkin (LSPG) projection (i.e., $\displaystyle\testbasis = \frac{\partial\res}{\partial\stateDummy}(\trialbasis\stateRed(\params)+ \stateRef;\params)\trialbasis$) \cite{bui2008model,LeGresleyThesis,CarlbergGappy,carlbergJCP,carlbergGalDiscOpt}.
The nonlinear terms in \eqnref{eq:rom_eqn} can be further approximated using `hyper-reduction' techniques to ensure that solving the ROM
equations incurs an $\nstate$-independent computational cost; these techniques include
collocation \cite{astrid2007mpe}, gappy POD
\cite{sirovichOrigGappy,carlbergJCP},
and the empirical interpolation method (EIM)
\cite{barrault2004eim,chaturantabut2010journal}. 
Note that model reduction constitutes a particular type of lower-fidelity model
characterized by $\resLF:(\stateDummy;\paramsDummy)\mapsto
\testbasis^T\res(\trialbasis\stateDummy+ \stateRef;\paramsDummy)$, $\stateLF =
\stateRed$, $\nstateLF = \nbasis$, and
$\prolongateNo:\stateDummy\mapsto\trialbasis\stateDummy+\stateRef$.

\subsection{Typical approaches for error quantification} 
\label{subsec:error_analysis}

Regardless of the particular approach used to generate the approximate
solution $\stateApprox$, it is essential to quantify the error incurred by
employing the approximate solution $\stateApprox$ in lieu of the high-fidelity solution $\state$. This work focuses on quantifying two such errors: 
\begin{enumerate} 
\item 
the error in the quantity of interest
$	\errorQoi (\params)\defeq \qoi(\params) - \qoiApprox(\params)$,
and 
\item the normed solution error
	$
	\errorState(\params) \defeq \|\errorStateVec(\params)\|
	$ with $\errorStateVec(\params)  \defeq \state(\params)
	-\stateApprox(\params)$.
 \end{enumerate}
We now describe two approaches that are typically employed to quantify these
errors: (1)
the dual-weighted residual, which is an \textit{error indicator} that
comprises a \reviewerOne{first-order} approximation of the quantity-of-interest
error $\errorQoi$, and (2) \textit{\textit{a posteriori} error bounds}, which
can be derived both for the solution error $\errorState$ and for the absolute
value of the quantity-of-interest error $|\errorQoi|$.

\subsubsection{Error in the quantity of interest: the dual-weighted residual error indicator} 
\label{subsubsec:local_errors}
The quantity-of-interest error $\errorQoi (\params)$ can be approximated using Taylor-series expansions of the residual and the quantity of interest via the dual-weighted residual. In particular, if the residual is twice continuously differentiable, it can be approximated \reviewerOne{to first order} about the approximate solution $\stateApprox$ by
\begin{align}
	\res (\state(\params) ;\params ) = \zero  =  \resApprox(\params)  +
	\jacobian(\params) \errorStateVec(\params)  +
	\bigO (\|\errorStateVec(\params)  \|^2),
	\quad\text{as $\|\errorStateVec(\params)  \|\rightarrow 0$},
\label{eq:resFirstOrder}   
\end{align}
where $\resApprox(\params)  \defeq \res (\stateApprox(\params);\params
)\in\RR{\nstate}$ denotes the residual evaluated at the approximate solution
and $\displaystyle\jacobian(\params)\defeq\frac{\partial \res
}{\partial\stateDummy}(\stateApprox(\params);\params)\in\RR{\nstate\times\nstate}$
denotes the Jacobian of the residual evaluated at the approximate solution.
\eqnref{eq:resFirstOrder} can be solved for an approximation of the
solution error:
\begin{align}
	\errorStateVec(\params)  =  -\jacobian
	(\params)^{-1}\resApprox(\params)+\bigO (\|\errorStateVec(\params)
	\|^2),
	\quad\text{as $\|\errorStateVec(\params)  \|\rightarrow 0$}.
\label{eq:e_oe}
\end{align}
%

If the quantity of interest is twice continuously differentiable, it also can
be approximated \reviewerOne{to first order} by
\begin{align}
	\qoi(\params) = \qoiApprox(\params) + 
	\frac{\partial \qoiFunc }{\partial \stateDummy
	}(\stateApprox(\params))\errorStateVec(\params)  +
	\bigO (\|\errorStateVec (\params) \|^2),
	\quad\text{as $\|\errorStateVec(\params)  \|\rightarrow 0$}.
\label{eq:qoiFunc}
\end{align}
Substituting the approximation of the solution error $\errorStateVec(\params)$ from \eqnref{eq:e_oe} into \eqnref{eq:qoiFunc} yields
\begin{align}
	\errorQoi(\params)=
-\frac{\partial \qoiFunc }{\partial \stateDummy
	}(\stateApprox(\params))\jacobian (\params)^{-1}\resApprox(\params)
	+ \bigO (\|\errorStateVec(\params)  \|^2).
	\quad\text{as $\|\errorStateVec(\params)  \|\rightarrow 0$}.
\label{eq:g2}
\end{align}
Defining the dual or adjoint 
$\dual:\RR{\nparams}\rightarrow\RR{\nstate}$
as the solution to the dual linear system
\begin{equation}
	\label{eq:dual_solve}
	\jacobian (\params)^{T}\dual(\params) =
-\frac{\partial \qoiFunc }{\partial \stateDummy
	}(\stateApprox(\params))^T,
\end{equation}
and substituting the dual $\dual$ into \eqnref{eq:g2} yields
\begin{align}
	\errorQoi(\params) = \dualweightedresidual(\params)  +
	\bigO (\|\errorStateVec(\params)  \|^2),
	\quad\text{as $\|\errorStateVec(\params)  \|\rightarrow 0$},
\label{eq:g3}
\end{align}
where the  dual-weighted residual 
$\dualweightedresidual:\RR{\nparams}\rightarrow\RR{}$ is defined as a weighted
sum of residual elements, that is,
\begin{align*}
	\dualweightedresidual(\params) \defeq \dual(\params)^T\resApprox(\params)=
	\sum_{i=1}^{\nstate} \dualEntry{i}(\params)\resApproxEntry{i}(\params) .
\end{align*}

Dual-weighted residuals are commonly employed for goal-oriented error
estimation and adaptive refinement, because they provide a \reviewerOne{first-order}
approximation of the error in the quantity of interest. Noting that 
\begin{align*}
	|\dualweightedresidual(\params)| \leq 	\sum_{i=1}^{\nstate}
	|\dualEntry{i}(\params)||\resApproxEntry{i}(\params) |,
\end{align*}
the absolute values of the dual elements $|\dualEntry{i}(\params)|$ can be interpreted as \textit{indicators} that inform the extent to which the absolute value of the associated residual element $|\resApproxEntry{i}(\params)|$ contributes to the quantity-of-interest error; this provides guidance for adaptive refinement.

In the present context, employing dual-weighted residuals $\dualweightedresidual(\params)$ for error quantification poses several challenges. First,
dual-weighted residuals are computationally costly to compute, as solving the dual
linear system
\eqref{eq:dual_solve} requires solving an $\nstate$-dimensional system of
linear equations; this cost is typically reduced by employing a lower-fidelity
model  (e.g., coarse mesh) for the dual solve and subsequently prolongating the
dual to the $\nstate$-dimensional state space (e.g., fine mesh)
\cite{venditti2000adjoint}. Second, dual-weighted residuals can pose an
implementation challenge, as computing the dual in \eqnref{eq:dual_solve}
requires the transpose of the Jacobian, which is not always
available, such as when the Jacobian is available only as a black-box operator in
a Jacobian-free Newton--Krylov setting \cite{knoll2004jacobian}. Third, for the
purpose of uncertainty quantification, there is no assurance that the
dual-weighted residual $\dualweightedresidual(\params)$ will be a low-bias
estimate of the quantity-of-interest error $\errorQoi$. Fourth, it provides a deterministic approximation of the
error; a statistical model is needed to quantify the epistemic uncertainty
introduced by the approximate solution.  Despite these challenges,
dual-weighted residuals remain \textit{informative} of the quantity-of-interest
error. Section \ref{subsec:indicators} describes how they can be used
as features in a machine-learning regression setting to construct accurate
statistical models of this error.

\subsubsection{Normed solution error: \textit{a posteriori} error bound} 
\label{subsubsec:normed_state_space_errors}
Residual-based bounds for the normed solution error $\errorState(\params)$
constitute a common approach for \textit{a posteriori} error quantification.
These approaches first assume that the residual $\res(\cdot;\params)$ is both inverse Lipschitz
continuous (i.e., inf--sup stable) and Lipschitz continuous, i.e., 
\begin{equation}
\invlipschitzConstant(\params)\|\stateOne-\stateTwo
\|\leq
\|\res(\stateOne;\params)-\res(\stateTwo;\params)\|\leq
\lipschitzConstant(\params)\|\stateOne-\stateTwo\|,\quad
\forall\stateOne,\stateTwo\in\RR{\nstate},
\label{eq:lgen1}
\end{equation}
where $\invlipschitzConstant:\RR{\nparams}\rightarrow\RRplus{}$ and
$\lipschitzConstant:\RR{\nparams}\rightarrow\RRplus{}$ denote the inverse
Lipschitz (i.e., inf--sup) and Lipschitz constants, respectively.
Note that if the residual is linear in its first argument and the norm is
taken to be the Euclidean norm, then the Lipschitz constants
$\invlipschitzConstant(\params)$ and $\lipschitzConstant(\params)$ correspond
to the minimum and maximum singular values of the Jacobian $\partial
\res/\partial\stateDummy(\params)$, respectively.
Substituting $\stateOne\leftarrow\state(\params)$ and
$\stateTwo\leftarrow\stateApprox(\params)$ in Inequalities \eqref{eq:lgen1}
with \eqnref{eq:fom_eqn} yields
\begin{align}
	\frac{\resNorm{\params}}{\lipschitzConstant(\params )}\le
	\errorState(\params) \le
	\frac{\resNorm{\params}}{\invlipschitzConstant(\params)}.
\label{eq:lgen2}
\end{align} 

Similarly, assuming the quantity-of-interest functional $\qoiFunc$ is
Lipschitz continuous, i.e.,
\begin{equation}
|\qoiFunc(\stateOne)-\qoiFunc(\stateTwo)
|\leq
\lipschitzConstantQoiFunc(\params)\|\stateOne-\stateTwo\|,\quad
\forall\stateOne,\stateTwo\in\RR{\nstate},
\label{eq:lgenQoiFunc1}
\end{equation}
where $\lipschitzConstantQoiFunc:\RR{\nparams}\rightarrow\RRplus{}$ denotes the
Lipschitz constant, then substituting
$\stateOne\leftarrow\state(\params)$ and
$\stateTwo\leftarrow\stateApprox(\params)$ in Inequality
\eqref{eq:lgenQoiFunc1} and making use of Inequality \eqref{eq:lgen2}
yields
\begin{equation}
|\errorQoi (\params)
|\leq
\frac{\lipschitzConstantQoiFunc(\params)}{\invlipschitzConstant(\params)}\resNorm{\params}.
\label{eq:lgenQoiFunc2}
\end{equation}

Although rigorous, the error bounds \eqref{eq:lgen2} and
\eqref{eq:lgenQoiFunc2} pose several challenges in the present context. First,
they are not typically sharp, as the upper (resp.~lower) bounds can
significantly overpredict (resp.~underpredict) the actual error. Second, they are often challenging to implement, as it is typically
difficult to compute the true Lipschitz/inf--sup constants for a given problem; rather,
bounds
$\invlipschitzConstantLB(\params)\leq \invlipschitzConstant(\params)$,
$\lipschitzConstantUB(\params)\geq\lipschitzConstant(\params)$, and
$\lipschitzConstantQoiFuncUB(\params)\geq\lipschitzConstantQoiFunc(\params)$
 must be computed
and employed within the bounds. These Lipschitz/inf--sup constant bounds often lack
sharpness, and improving their sharpness can be computationally costly \cite{SCM,wirtzDeim}. Third, these bounds do not
produce a statistical distribution over the normed solution error, which
precludes quantifying the epistemic uncertainty introduced by the
approximation. While the lower and upper bounds can in principle define the
interval for a uniform distribution, this distribution is typically not
representative of the actual behavior of the error, as demonstrated in Ref.~\cite{drohmann_2015}.
Again, despite these challenges, \textit{a posteriori} error bounds remain
\textit{informative} of the normed solution error, and Section
\ref{subsec:indicators} describes how they can be used in a machine-learning
regression setting to construct accurate statistical models of this error.

\section{Machine-learning error models} 
\label{sec:approach}
In the spirit of the ROMES method \cite{drohmann_2015}, this work aims to
construct statistical models of the quantity-of-interest error $\errorQoi$ and
normed solution error $\errorState$. However, in contrast to the original
ROMES method that relies on Gaussian-process regression and the hand-selection
of a small number of error indicators that can be relatively costly to
compute, the proposed method applies high-dimensional regression methods from
machine learning. This enables a larger number of inexpensive error indicators
to be considered, resulting in less costly, more accurate error models. In
this way, the proposed method exhibits similarities to the approach proposed in
Ref.~\cite{trehan2017error} for modeling the error in dynamical-system
surrogates.

\subsection{Machine-learning framework} 
\label{subsec:overview}

We begin by assuming that $\nfeatures$ \textit{error indicators} or \textit{features}
$\features(\params)\in\RR{\nfeatures}$ with $\features\equiv[\featureArg{1}\
\cdots\ \featureArg{\nfeatures}]^T$---which are 
informative of the
error of interest $\errorQoi(\params)$ or $\errorState(\params)$---can be
cheaply computed from the solution approximation $\stateApprox(\params)$.
 Section \ref{subsec:indicators} discusses how
these features can be engineered based on the analysis performed in Section
\ref{subsec:error_analysis}.

Given the ability to compute these error indicators $\features(\params)$, we
model the nondeterministic mapping $\features(\params)\mapsto
\error(\params)$, where we have denoted 
the
error of interest $\errorQoi(\params)$ or $\errorState(\params)$
generically as
$\error(\params)\in\RR{}$, using an additive error model comprising the sum of a deterministic
regression function $\trueRegressionModel$ and stochastic
noise\footnote{This is often referred to as the `error' in the machine-learning literature. We refer to it as
	`noise' to distinguish it from the approximate-solution error $\error$ we
aim to model using regression.} $\noise$, as 
\begin{equation} \label{eq:trueErrorModel}
	 \error(\params) = \trueRegressionModel(\features(\params)) +
	 \noise(\features(\params)).
 \end{equation} 
 Here, the noise $\noise$ is a mean-zero random variable accounting for
 irreducible error in Formulation \eqref{eq:trueErrorModel} due 
 to omitted explanatory variables; thus, we consider it to represent epistemic
 uncertainty, as including additional features can in principle enable zero
 noise (see the discussion of Feature-Engineering Method
 \ref{feat:params} in Section \ref{subsec:indicators}). Its
 dependence on
 the features enables feature-dependent distribution parameters (e.g.,
 variance in the case of Gaussian noise) to be considered. Thus, the regression function defines
 the conditional expectation of the error given the features, i.e.,
 \begin{equation}
	 \expectationGiven{\error(\params)}{\features(\params)} = 
	 \trueRegressionModel(\features(\params)).
 \end{equation}

The proposed methodology constructs models of both the
deterministic regression function
$\regressionModel(\approx\kern-0.25em\trueRegressionModel)$ and the stochastic noise
$\noiseModel(\approx\noise)$, which in turn yield a statistical model for the
approximate-solution error:
 \begin{equation} 
\errorModel(\params) = \regressionModel(\features(\params)) + \noiseModel(\features(\params)).
 \end{equation} 
The methodology aims to
construct this regression model $\errorModel$ such that it satisfies three
objectives:
\begin{Objective} 
	\item \label{objective:cheap}The regression model should employ \textit{cheaply computable} features $\features$,
	\item \label{objective:lowvar} The regression model should exhibit
		\textit{low noise variance}; that is,
	$\variance{\noiseModel}$ should be small, as the noise variance
	quantifies
	the epistemic uncertainty introduced by the approximate solution, and 
\item \label{objective:validated}The regression model should \textit{generalize}; it should be
	numerically validated such that the empirical distributions of
	$\errorModel$ and $\error$ are `close' on an independent test set
	$\testData\defeq\{(\errorTest{i},\featuresTest{i})\}_{i=1}^{\ntestData}$
	that is not
	used to train the model. A generalizable model is one that has not been
	overfitted on training data.
\end{Objective}
We note that these objectives were originally described by the ROMES
method~\cite{drohmann_2015}.

We propose to construct the regression model in three steps, which address
these objectives:
\begin{Step} 
\item \textit{Feature engineering}. Devise features $\features$ that are 
	cheaply computable	(\ref{objective:cheap}),  informative of the error such that a
	low-noise-variance model can be constructed (\ref{objective:lowvar}), and low dimensional (i.e., $\nfeatures$ small) such that less training data is
	needed to obtain a
	generalizable model (\ref{objective:validated}). We note that in
	the case of a very large training set (i.e., $\ntrainingData
	$ large), 
	regression methods from representation learning (e.g., deep neural networks)
	could be applied with an extremely large set of candidate features,
	rendering feature engineering less critical. However,
this scenario is not expected to occur in the present context, as
	each training data point requires solving the high-fidelity-model equations
	(unless multiple types of approximate solutions are simultaneously considered; see
	\ref{data:pooled} in Section \ref{subsection:trainingdata}). Section
\ref{subsec:indicators}
describes feature engineering.
\item \textit{Regression-function modeling}. Construct the
(deterministic) regression-function model $\regressionModel$ by
applying regression methods from machine learning to approximate the mapping
from features $\features$ to	approximate-solution error $\error$ using a
training set
$\trainingData\defeq\{(\errorTrain{i},\featuresTrain{i})\}_{i=1}^{\ntrainingData}$.
Sections \ref{subsection:trainingdata} and
\ref{subsection:regressionfunctionapprox} describe this step.
\item\label{step:noiseApproximation} \textit{Noise modeling}. Model the
	stochastic noise
model $\noiseModel$ as a mean-zero Gaussian random variable with constant variance,
i.e.,
$\noiseModel\sim\normal{0}{\varianceEstimate}$, where $\varianceEstimate$ is
the sample variance of the error on the test set
$\testData$, i.e., $\varianceEstimate =
\frac{1}{\ntestData}\sum_{i=1}^{\ntestData}
(\errorTest{i}-\regressionModel(\featuresTest{i}))^2$.
	\footnote{The
denominator corresponds to $\ntestData$ and not $\ntestData-1$, as the mean is
known to be zero.} 
	Note that more complex noise models could be considered.  For example, 
	the ROMES method \cite{drohmann_2015} uses
	Gaussian-process
regression, which enables modeling heteroscedastic noise.
\end{Step}

\subsection{Feature engineering} 
\label{subsec:indicators}

The goal of feature engineering is to devise a set of cheaply computable
regression-model inputs $\features$ that satisfy 
\ref{objective:cheap} to \ref{objective:validated}.  We propose a range of
possible features that accomplish this by balancing two attributes:
\begin{Attribute}
\item\label{char:nfeatures} \textit{Number of features $\nfeatures$}.
Using a large number of
	features can generally lead to a low-noise-variance regression
	model, as including more explanatory variables enables reduction of the
	epistemic uncertainty and thus the noise in the regression-model formulation
	\eqref{eq:trueErrorModel} (i.e., \ref{objective:lowvar} is bolstered). However, employing
	a large number of features incurs two drawbacks: 
	(1) the features may no longer be cheaply computable (i.e., 
	\ref{objective:cheap} suffers), and 
	(2) using more features usually implies a higher capacity (i.e., 
	lower bias and
	higher variance) regression model, which in turn requires more training data
	(i.e., $\ntrainingData$ large) to generalize (i.e., 
	\ref{objective:validated} suffers without sufficient training data). Employing regularization while
	training the regression model mitigates the latter effect.
\item\label{char:featurequality}\textit{Quality of features}.
Employing high-quality features can lead
	to a low-noise-variance regression model, as such features provide a strong
	explanation of the response quantity and thus reduce the epistemic
	uncertainty and therefore the noise in the regression-model Formulation 
	\eqref{eq:trueErrorModel} (i.e., \ref{objective:lowvar} is bolstered).
	However, computing high-quality features is generally computationally
	expensive (i.e., \ref{objective:cheap} suffers).  Thus, it is often
	advantageous to reduce the computational cost at the expense of a
	higher-noise-variance model by employing lower-quality features that are
	cheaply computable, yet remain informative of the response quantity.
\end{Attribute}
We now describe feature engineering for modeling the error. This
task employs the  error analysis presented in
Section~\ref{subsec:error_analysis}, as it provides insight into 
quantities that are informative for quantifying both quantity-of-interest
errors and normed solution errors.
Table~\ref{tab:featuresQoI} summarizes the candidate feature-engineering
methods. 

\begin{table}[htbp] 
\small
\centering 
\begin{tabular}{c l c c c c c c}
\toprule
\begin{tabular}{@{}c@{}} Method \\ Index \end{tabular} & 
\multicolumn{1}{c}{Method Name} & 
\begin{tabular}{@{}c@{}}Features \\ $\features(\params)$ \end{tabular}&
\begin{tabular}{@{}c@{}}   Number of\\ features \\$\nfeatures$ \end{tabular} & 
\begin{tabular}{@{}c@{}}   Number of\\ residual \\ elements\\required
\end{tabular}&
\begin{tabular}{@{}c@{}}Dual \\ solve \\ required? \end{tabular}&
\begin{tabular}{@{}c@{}}   Applicable\\ error $\error$ \end{tabular}&
\begin{tabular}{@{}c@{}} $\nstate$-indep. \\ cost? \end{tabular}
\\
\midrule
\ref{feat:params} &
Parameters & 
$\params$ & 
$\nparams$ & 
$0$ & \xmark&$\errorQoi$, $\errorState$ & \checkmark
\\[.75em] 
\ref{feat:DWR} &
Dual-weighted residual & 
$\dualweightedresidual(\params)$ & 
1 & 
$\nstate$ 
&\checkmark& $\errorQoi$ & \xmark
\\[.75em] 
\ref{feat:parameters_resnorm}& Parameters and residual norm & 
$[\params;\, \resNorm{\params}]$ & 
$\nparams+1$ &
$\nstate$ & \xmark&$\errorQoi$, $\errorState$
& \xmark
\\[.75em] 
\ref{feat:resnorm}& Residual norm & 
$\resNorm{\params}$ & 
$1$ &
$\nstate$ & \xmark&$\errorQoi$, $\errorState$
& \xmark
\\[.75em] 
\ref{feat:paramres} &
Parameters and residual & 
$[\params;\,\resApprox(\params)]$ & 
$\nparams + \nstate$ & 
$\nstate$ & 
\xmark&
$\errorQoi$, $\errorState$& \xmark
\\[.75em] 
\ref{feat:paramresPCA} &
\begin{tabular}{@{}l@{}}Parameters
	and residual \\principal components \end{tabular}&
$[\params;\,\resApproxRed(\params)]$ & 
$\nparams + \dimresApprox$ & 
	$\nstate$ &\xmark& $\errorQoi$, $\errorState$
& \xmark
\\[1.5em] 
\ref{feat:paramgappyPODres} &
\begin{tabular}{@{}l@{}}Parameters and residual \\ gappy principal components
\end{tabular}&
$[\params;\,\resApproxRedGappy(\params)]$ & 
$\nparams + \dimresApprox$ & 
	$\nsamples$ & \xmark&$\errorQoi$, $\errorState$
& \checkmark
\\[1.5em] 
\ref{feat:paramsampleres} &
\begin{tabular}{@{}l@{}}Parameters \\ and residual samples \end{tabular}& 
$[\params;\,\sampleMat\resApprox(\params)]$ & 
$\nparams + \nsamples$ & 
$\nsamples$ &\xmark& $\errorQoi$, $\errorState$
& \checkmark
\\ 
\bottomrule 
\end{tabular} 
\caption{Proposed feature-engineering methods.  Note that
$\dimresApprox\leq\nsamples$ and typically $\nsamples,\nparams\ll\nstate$.} 
\label{tab:featuresQoI} 
\end{table} 
We consider the following feature-engineering methods:
\begin{enumerate} [leftmargin=1.4em]
	\item \label{feat:params} \textbf{Parameters}. 
Because the mapping $\params\mapsto \error(\params)$ is deterministic for both
the quantity-of-interest error $\error=\errorQoi$ and normed solution error
$\error=\errorState$, we
could employ the parameters as features for modeling these errors:
\begin{equation}
\features(\params)=\params.
\end{equation}
 Then, with a sufficiently large amount of training data and a
 regression-model form with sufficient capacity,
it is possible to construct a regression
model with zero noise variance; this is the optimal outcome for  \ref{objective:lowvar}.
However, in practice, this is unlikely to be effective due to lack of feature
quality, as the mapping
$\params\mapsto \error(\params)$ is often complex and difficult to model. 
For example, in the case of reduced-order models,
the error is usually zero for points in the parameter space where data
were collected \cite{drohmann_2015,NE12}, leading to a highly oscillatory 
$\params\mapsto \error(\params)$ mapping.
Thus, this approach typically corresponds to a 
small number ($\nfeatures=\nparams$) of low-quality features.

As previously mentioned, this choice for features corresponds to the approach
		taken by `multifidelity correction' or
		`model discrepancy'	methods, wherein kriging or radial basis functions are
		typically employed to construct the regression-function model. This choice
		also yields an \textit{a priori} error model, as the parameters $\params$
		are known before computing the approximate solution, and the approach does
		not make use of any data generated by the approximate solution.
\item \label{feat:DWR} \textbf{Dual-weighted residual}.  \eqnref{eq:g3} in
	Section \ref{subsubsec:local_errors} demonstrates that the dual-weighted
	residual $\dualweightedresidual(\params)$ is a \reviewerOne{first-order}
	approximation of the quantity-of-interest error $\errorQoi(\params)$.
	Consequently, the dual-weighted residual can be chosen as a 
	feature for modeling the quantity-of-interest error $\errorQoi$:
\begin{align*}	
	\features(\params)=\dualweightedresidual(\params)\defeq \dual(\params)^T\resApprox(\params).
\end{align*}
This approach corresponds to a small number ($\nfeatures=1$) of high-quality
features.  Thus, this approach is expected to lead to a
low-noise-variance regression model and requires a relatively small amount of
training data to generalize.
However, as discussed in Section \ref{subsubsec:local_errors}, using this high-quality feature is
computationally costly, as computing the required dual vector $\dual(\params)$
entails solving the $\nstate$-dimensional dual linear system~\eqref{eq:dual_solve}; furthermore, it can pose an implementation
challenge, as the transpose of the Jacobian---which is needed to compute the
dual $\dual(\params)$---is not always available.  
This approach was pursued by the ROMES method
\cite{drohmann_2015}, which reduced the cost of computing the dual
by applying model reduction to the dual linear system \eqref{eq:dual_solve}.
\item \label{feat:parameters_resnorm}
\textbf{Parameters and residual norm}.
Inequalities \eqref{eq:lgen2} demonstrate that the normed solution error $\errorState(\params)$ can be bounded
from above (resp.\ below) by a parameter-dependent constant
${1}/{\invlipschitzConstant(\params)}$
(resp.\ ${1}/{\lipschitzConstant(\params )}$) multiplied by the 
residual norm $\resNorm{\params}$. In addition, Inequality
\eqref{eq:lgenQoiFunc2} demonstrates that the absolute value of the
quantity-of-interest error
$|\errorQoi (\params)
|$ can be bounded from above by a parameter-dependent constant
${\lipschitzConstantQoiFunc(\params)}/{\invlipschitzConstant(\params)}$
multiplied by the residual norm $\resNorm{\params}$.
Thus, the residual norm and parameters are
sufficient quantities for bounding these errors and
can thus be employed as features in the proposed approach: 
\begin{align}
\features(\params)=[\params;\,\resNorm{\params}].
\label{eq:parameters_resnorm}
\end{align}
\noindent In contrast to Feature-Engineering Method \ref{feat:DWR}, this approach avoids
the need for any dual solves;
furthermore, it requires a relatively small number of features, as $\nfeatures
= \nparams + 1$. However, these may be low quality when
modeling the quantity-of-interest error $\errorQoi$, as the
residual norm is informative of only the \textit{absolute value} 
$|\errorQoi|$,
according to Inequality~\eqref{eq:lgenQoiFunc2}; thus, it may
exhibit difficulties in discerning the sign of the quantity-of-interest error
$\errorQoi$.
%
\item \label{feat:resnorm} \textbf{Residual norm}.  
Rather than employing both the parameters $\params$ and residual norm
$\resNormNo$ as is done by Feature-Engineering Method
\ref{feat:parameters_resnorm}, we can instead simply employ the residual norm as a feature:
\begin{align}
	\features(\params)= \resNorm{\params} .
\label{eq:resnorm}
\end{align}
This method shares many of the same attributes as
Feature-Engineering Method \ref{feat:parameters_resnorm}: it avoids dual
solves, it employs a small number of features, and it may be low quality for
modeling the quantity-of-interest error $\errorQoi$. However, this feature may be more
appropriate if the parameters are low quality (as mentioned
in the discussion of Feature-Engineering Method \ref{feat:params}), or if
the problem is characterized by a high-dimensional parameter space (i.e.,
$\nparams$ large) and there is an insufficient amount of training data.
This feature choice was also pursued by the ROMES method
\cite{drohmann_2015} for modeling normed solution errors
$\errorState$.
\item \label{feat:paramres}
\textbf{Parameters and residual}. 
Noting that the (high-quality) dual-weighted residual
$\dualweightedresidual(\params) \defeq \dual(\params)^T\resApprox(\params)$ is
simply a weighted sum of the elements of residual
vector with parameter-dependent weights, the parameters $\params$ and residual $\resApprox(\params)$ can
be employed directly as features for modeling the quantity-of-interest
error $\errorQoi$: 
\begin{align}
\features(\params)=[\params;\,\resApprox(\params)].
\label{eq:res_parameters}
\end{align}
This feature choice can also be justified for modeling the normed solution
error $\errorState$, as the features of Feature-Engineering Method
\ref{feat:parameters_resnorm}---which are informed by error-bound
analysis---can be directly recovered from these quantities.
This approach amounts to employing a large number (i.e.,
$\nfeatures=\nparams+\nstate$) of low-quality features, as each element of the
parameters or residual may itself be a poor predictor of the
error. Thus, this approach may require a significant
amount of training data (i.e., $\ntrainingData$ large) to avoid overfitting.
Furthermore, computing the features is computationally costly, as it requires
computing all $\nstate$ elements of the residual $\resApprox(\params)$.
However, it does not require any dual solves, which constitutes a practical and
computational-cost improvement over Feature-Engineering Method \ref{feat:DWR}.
\par
We now propose several techniques for reducing the number of features
$\nfeatures$ in order to reduce both the amount of required training data
and the cost of applying the regression model online. In principle, we could
apply a number of standard techniques, including subset selection (e.g., forward-
and backward-stagewise regression), shrinkage (e.g., ridge, lasso), or derived
inputs (e.g., principal-components regression, partial least squares), for this
purpose. We focus on employing principal components, as this technique mirrors
the gappy POD approach \cite{sirovichOrigGappy} often employed in model reduction.
\item \label{feat:paramresPCA}
\textbf{Parameters and residual principal components}.  Because the elements
of the residual vector $\resApprox(\params)$ tend to be highly correlated,
rather than employing the entire $\nstate$-dimensional residual vector
$\resApprox(\params)$ as features, we can instead represent the residual in
terms of its principal components:
\begin{align*}
\features(\params)=[\params;\,\resApproxRed(\params)],
\end{align*}
with
\begin{align}
\resApproxRed(\params) \defeq \basisresApprox^T(\resApprox(\params) - \resApproxRef). 
\label{eq:redRes}
\end{align} 
Here, $\basisresApprox\in\orthoMat{\nstate}{\dimresApprox}$ denotes the matrix
whose columns comprise the first $\dimresApprox\ll \nstate$ principal
components of the training data
$\trainingDataRes\defeq\{\resApprox_i\}_{i=1}^{\ntrainingDataRes}$;
$\orthoMat{n}{m}\subset\RR{n\times m}$ denotes the Stiefel manifold, which is
the set of all real-valued $n\times m$ matrices with orthonormal columns;  and
$\resApproxRef\defeq\frac{1}{\ntrainingDataRes}\sum_{i=1}^{\ntrainingDataRes}\resApprox_i\in\RR{\nstate}$. This approach reduces the number of features to
$\nfeatures = \nparams + \dimresApprox$; thus it will likely require less
training data to generalize (i.e., $\ntrainingData$ smaller) than
Feature-Engineering Method~\ref{feat:paramres}. However, its online cost remains large, as
computing $\resApproxRed(\params)$ via \eqnref{eq:redRes} requires
first evaluating the entire $\nstate$-dimensional residual vector
$\resApprox(\params)$. Fortunately, this cost can be reduced using the gappy
POD method. 
\item \label{feat:paramgappyPODres}
\textbf{Parameters and residual gappy principal components}. 
The gappy POD method \cite{sirovichOrigGappy} reconstructs vector-valued data
that have `gaps', that is, entries with unknown or uncomputed values. It is
equivalent to least-squares regression in one discrete-valued variable using
an empirically computed basis; it was devised by Everson and Sirovich
\cite{sirovichOrigGappy} for the purpose of image reconstruction. It has also
been employed for flow-field reconstruction
\cite{bui2003proper,venturi2004gappy,willcox2006ufs}, inverse design
\cite{bui2004aerodynamic}, \reviewerOne{compressed sensing \cite{peherstorfer_2016},} and for decreasing the spatial
\cite{astrid2007mpe,bos2004als,CarlbergGappy,carlbergJCP,carlbergGalDiscOpt}
and temporal \cite{carlberg2015decreasing, carlbergParareal,choiCarlberg}
complexity in model reduction.  In the present context, this technique can be
applied to approximate the generalized coordinates $\resApproxRed(\params)$ in
\eqnref{eq:redRes} from a \textit{sampled} subset of elements of
$[\resApprox(\params) - \resApproxRef]$. In particular, the method
approximates these generalized coordinates as 
\begin{align*} 
\resApproxRedGappy(\params) \defeq (\sampleMat\basisresApprox)^+\sampleMat[\resApprox(\params) - \resApproxRef]\approx\resApproxRed(\params),
\end{align*} 
where the superscript $+$ denotes the Moore--Penrose pseudoinverse and
$\sampleMat\in\{0,1\}^{\nsamples\times\nstate}$ denotes a sampling matrix
comprising $\nsamples$ rows of $\identityArg{\nstate}$,
where $\identityArg{n}$ denotes the $n\times n$
identity matrix and 
$\dimresApprox\leq\nsamples\ll\nstate$.\footnote{Note that gappy
POD is equivalent to the (discrete) empirical interpolation method
\cite{barrault2004eim,chaturantabut2010journal} when $\nsamples =
\dimresApprox$, as the pseudo-inverse is equal to the inverse in this case.}
Thus, this approach employs features
\begin{align*}
\features(\params)=[\params;\,\resApproxRedGappy(\params)],
\end{align*}
which are less expensive to compute than the features of
Feature-Engineering Method~\ref{feat:paramresPCA}, as computing $\resApproxRedGappy(\params)$
requires computing only $\nsamples\ll \nstate$ elements of the vector
$\resApprox(\params) - \resApproxRef$; however, these features are expected to
be of (slightly) lower quality, as the generalized coordinates have been approximated.

The sampling matrix $\sampleMat$ can be computed by a variety of
techniques. In this work, we consider two approaches. First, we consider
\textit{q-sampling} \cite{drmacGugercinReview}, wherein the sampling matrix
$\sampleMat$ consists of the transpose of the first $\nsamples$ columns of the
permutation matrix
$\sampleMatFull\in\{0,1\}^{\nstate\times\nstate}$ that arises from the QR
factorization  $\basisresApprox^T\sampleMatFull=\mathbf{Q}\mathbf{R}$ (see
\cite[Algorithm 1]{drmacGugercinReview}). Here,
the pivoting provided by $\sampleMatFull$ ensures the diagonal elements of
$\mathbf{R}$ are non-increasing.
This approach was devised in the model-reduction community
as a less computationally expensive alternative (with sharper \textit{a
priori} error bounds) to more traditional greedy methods
\cite{barrault2004eim,chaturantabut2010journal}.

Second, we consider \textit{k-sampling}, a linear univariate feature selection
approach developed in the machine-learning community.  
Here, the sampling matrix $\sampleMat$ comprises the $\nsamples$ rows of
$\identityArg{\nstate}$ corresponding to the $\nsamples$
features with the highest scores on an
$F$-test~\cite{scikit-learn,sklearn_api} with respect to the response $\error$.  
The scores are computed by
\begin{align*}
F_i = \frac{\rho_i^2}{1-\rho_i^2}\left(\ntrainingData-2\right),
\end{align*}
where
\begin{align*}
\rho_i=\frac{\mathbf{e}_i^T\left(\resApprox-\bar{\resApprox}\right)\left(\error-\bar{\error}\right)}{\sqrt{\variance{\mathbf{e}_i^T\resApprox}\variance{\error}}},
\end{align*}
and $\mathbf{e}_i$ is the $i$th column of $\identityArg{\nstate}$ identity matrix. 
\item \label{feat:paramsampleres} 
\textbf{Parameters and residual samples}.
Noting that the features $\resApproxRedGappy(\params)$ employed by
Feature-Engineering Method~\ref{feat:paramgappyPODres} correspond to the sampled centered residual
$\sampleMat[\resApprox(\params) - \resApproxRef]$ premultiplied by a constant
matrix $(\sampleMat\basisresApprox)^+\in\RR{\dimresApprox\times\nsamples}$,
the sampled residual can instead be used as features:
\begin{align*}
\features(\params)=[\params;\,\sampleMat\resApprox(\params)].
\end{align*}
While the number of features $\nfeatures = \nparams + \nsamples$ is slightly
larger than in the case of Feature-Engineering Methods~\ref{feat:paramresPCA} and
\ref{feat:paramgappyPODres}, this approach does not in principle require the
computation of the basis $\basisresApprox$ (although the samples $\sampleMat$
are computed from the basis $\basisresApprox$ in the case of
\textit{q}-sampling as previously discussed).
Thus, if the computation of $\basisresApprox$ can be avoided with this
approach (as in the case of \textit{k}-sampling), it can lead to a lower-cost
training stage relative to Feature-Engineering Methods \ref{feat:paramresPCA} and \ref{feat:paramgappyPODres}.
\end{enumerate}

\subsection{Training and test data} 
\label{subsection:trainingdata}
This section describes the construction of the training set $
\trainingData
\defeq\{(\errorTrain{i},\featuresTrain{i})\}_{i=1}^{\ntrainingData}
$
and test set $\testData \defeq\{(\errorTest{i},\featuresTest{i})\}_{i=1}^{\ntestData}$, which are required for training the
regression-function approximation and noise approximation, respectively (see
Section \ref{subsec:overview}). \reviewerThree{Algorithm \ref{alg:dataGen} of
\ref{app:methodalgorithms} provides the resulting data-generation algorithm.}

We propose two approaches for constructing these data sets:
\begin{Data} 
\item\label{data:pooled} \textit{Pooled data set}. This approach constructs these data sets as
		the union of data sets generated by $\napprox\geq 1$ different types of
		approximate solutions $\stateApproxArg{i}(\params)$,
		$i=1,\ldots,\napprox$; these different solutions can correspond to
		reduced-order models of different dimensions $\nbasisArg{i}$,
		$i=1,\ldots,\napprox$, for example.  The resulting error model can then be
		deployed across all $\napprox$ different approximate solutions. The
		benefit of this approach is access to more training data, as $\napprox$
		training data points can be generated from a single (costly-to-compute)
		high-fidelity solution $\state(\params)$. However, the resulting error
		model may not generalize well, as feature--error relationships for
		different kinds of approximate solutions may exhibit different
		characteristics.
\item\label{data:unique} \textit{Unique data set.} This approach constructs a unique error model 
for each type of approximate solution considered, and thus employs a separate data 
	set for
	each type of approximate solution. This approach suffers
	from limited training data, as each error model has access to only one
	training point per high-fidelity solution $\state(\params)$. However,
	because the training data are specialized to a single type of approximate
	solution, the resulting error model is more likely to generalize well.
 \end{Data}

%
\begin{itemize}[leftmargin=1.4em]
\item \textbf{Training and test data $\trainingData$ and $\testData$}.
We define a set of parameter training instances
$\paramDomainTrain\subset\paramDomain$
and parameter test instances $\paramDomainTest\subset\paramDomain$ such that $
\paramDomainTrain\cap\paramDomainTest=\emptyset$, where---for each parameter
instance---both the
error $\error$ (i.e., the response to predict) and the features
$\features$ are computed for each of the $\napprox$ types of approximate
solutions.  \reviewerTwo{These instances are randomly sampled from the
		expected probability distribution (e.g., normal distribution, uniform
		distribution) defined on the parameter domain.}

Denoting by a
superscript $i$ the quantity computed by approximate
solution $\stateApproxArg{i}(\params)$, we define 
\begin{align*}
	\errorDataArg{i}(\paramDomainDummy)\defeq\{(\errorArg{i}(\params),\featuresArg{i}(\params))\,
|\, \params\in\paramDomainDummy\}
\end{align*}
as the set of feature--error pairs generated by the $i$th approximate solution on parameter
set $\paramDomainDummy\subseteq\paramDomain$. Then, 
the training and
test sets 
for \ref{data:pooled} correspond to
\begin{align}
\trainingData
=\bigcup\limits_{i=1}^\napprox\errorDataArg{i}(\paramDomainTrain)\quad\text{and}\quad
\testData
=\bigcup\limits_{i=1}^\napprox\errorDataArg{i}(\paramDomainTest),
\label{eq:trainingTestDataset1}
\end{align}
respectively, and
the training and test sets of \ref{data:unique} for the $i$th approximate
solution correspond to
\begin{align}
\trainingData
=\errorDataArg{i}(\paramDomainTrain)\quad\text{and}\quad
\testData
=\errorDataArg{i}(\paramDomainTest),
\label{eq:trainingTestDataset2}
\end{align}
respectively.
Note that for \ref{data:pooled}, 
$\ntrainingData = \napprox\times\card{\paramDomainTrain}$
and $\ntestData =\napprox\times\card{\paramDomainTest}$,
while for \ref{data:unique}, 
$\ntrainingData=\card{\paramDomainTrain}$
and $\ntestData =\card{\paramDomainTest}$.
\item \textbf{Training data for residual PCA $\trainingDataRes$}. Feature-Engineering Methods~\ref{feat:paramresPCA} and \ref{feat:paramgappyPODres} require the
	principal components $\basisresApprox$; Feature-Engineering Method~\ref{feat:paramsampleres}
	also requires these principal components if they are employed to define the
	sampling matrix $\sampleMat$ (e.g., via \textit{q}-sampling).  
Define 
\begin{align*} 
\residualDataArg{i}(\paramDomainDummy)\defeq\{\resApproxArg{i}(\params)\,|\,
\params\in\paramDomainDummy\}
\end{align*} 
as the set of residuals 
	generated by the $i$th approximate solution over parameter set
	$\paramDomainDummy\subseteq\paramDomain$,
	where $\resApproxArg{i}(\params)\defeq\res (\stateApproxArg{i}(\params);\params
)$. Then, because
principal component regression (PCR) typically employs the same data
used to train the regression model as those employed to compute the principal
components, we define the training data employed to compute these principal
components for \ref{data:pooled}
as
\begin{align} 
\trainingDataRes = 
\bigcup\limits_{i=1}^\napprox\residualDataArg{i}(\paramDomainTrain).
\label{eq:trainingResDataset1}
\end{align} 
The corresponding training data for \ref{data:unique} for the $i$th
approximate solution are
\begin{align} 
\trainingDataRes = 
\residualDataArg{i}(\paramDomainTrain).
\label{eq:trainingResDataset2}
\end{align} 
%
Because principal component
analysis is an unsupervised learning method, a test set is not required.
\end{itemize}
%


\subsection{Regression-function approximation} 
\label{subsection:regressionfunctionapprox}

We consider several different techniques---each of which exhibits a different
level of capacity---to construct the regression-function approximation
$\regressionModel$ from training data $\trainingData$.  In machine learning,
high-capacity models tend to be low bias and high variance and can thus lead
to very accurate models with \reviewerOne{low mean-squared error on an
independent test set and thus} low noise variance (i.e., \ref{objective:lowvar}
is bolstered); however, they generally require a large number of training
examples $\ntrainingData$ to generalize (i.e., \ref{objective:validated}
suffers without sufficient training data) and are thus prone to overfitting.
For many regression models, increasing the number of considered features
$\nfeatures$ leads to a higher-capacity model as mentioned in Section
\ref{subsec:indicators}. 

On the other hand, lower-capacity methods make stronger structural
assumptions and typically lead to higher bias and lower variance.  These
models typically require less training data to generalize (i.e.,
\ref{objective:validated} is bolstered) and are therefore less prone to
overfitting.  However, the high bias of these models can result in significant
prediction errors, ultimately resulting in a large noise variance (i.e.,
\ref{objective:lowvar} suffers), even when the amount of training data is
large.

For all regression-function approximations, we employ cross-validation within
the training set in order to tune hyperparameters characterizing each
regression model. \reviewerThree{Algorithm \ref{alg:regressionNoiseModel} of
\ref{app:methodalgorithms} provides the resulting algorithm for training a
regression model.}

\subsubsection{Ordinary least squares} 
\label{section:OLS}
Ordinary least squares (OLS) corresponds to a regression model that is linear
in the features $\features$; it takes the form
\begin{align}
	\regressionModel(\features;\weight) = \weightArg{0} + \sum_{p=1}^{\nfeatures}
	\weightArg{p} \featureArg{p},
\label{eq:ols_linear}
\end{align}
where $\weight\equiv\left[\weightArg{0}\
\cdots\ \weightArg{\nfeatures}\right]^T\in\mathbb{R}^{\nfeatures+1}$
denotes the weighting coefficient vector.
We denote this
method as OLS: Linear.  We also consider a regression model that is quadratic
in the features:
\begin{align}
\regressionModel(\features;\weight) = 
\weightArg{0} + \sum_{p=1}^{\nfeatures}
	\weightArg{p} \featureArg{p}+
\sum_{p=1}^{\nfeatures}
\sum_{q=p}^{\nfeatures}\weightArg{q+\frac{p}{2}(2\nfeatures+1-p)}
\featureArg{p} \featureArg{q}.
\end{align}
Here, $\weight\in\mathbb{R}^{(\nfeatures+1)(\nfeatures+2)/2}$.  We denote this
method as OLS: Quadratic.

In either case, we can compute the weights $\weight$ as the solution to the
mean-squared-error (MSE) minimization problem
\begin{align}
	\underset{\weightDummy}{\mathrm{minimize}}\
	\sum_{i=1}^\ntrainingData\left(\regressionModel(\featuresTrain{i};\weightDummy) -
	\errorTrain{i}\right)^2.
\label{eq:ols_constraint}
\end{align}
%
This is equivalent to maximum likelihood estimation (MLE) if the noise
$\noise$ in Eq.~\eqref{eq:trueErrorModel} is assumed to have a mean-zero Gaussian distribution with constant
variance, which we assume in our noise-approximation approach (see \ref{step:noiseApproximation} in Section
\ref{subsec:overview}).

When \eqnref{eq:ols_constraint} is underdetermined (i.e.,
$\nfeatures+1<\ntrainingData$ for OLS: Linear; 
$(\nfeatures+1)(\nfeatures+2)/2<\ntrainingData$ for OLS: Quadratic), a unique solution can be
obtained by a penalty formulation. We consider the ridge penalty, which
reformulates Problem \eqref{eq:ols_constraint} as
\begin{align}
\begin{split}
	\underset{\weightDummy}{\mathrm{minimize}}&\
\left\|\weightDummy\right\|_2,
	\\
	\text{subject to}& \
	\regressionModel(\featuresTrain{i};\weightDummy) -
	\errorTrain{i}=0,\quad i=1,\ldots,\ntrainingData.
\end{split}
\label{eq:ols_constraint2}
\end{align}

OLS: Linear is a relatively low-capacity model, as it imposes a linear
relationship between the features and the response;
OLS: Quadratic is higher capacity.
Because we only consider a penalty formulation in the event of an
underdetermined least-squares problem \eqref{eq:ols_constraint}, these
approaches do not have any hyperparameters to tune during cross-validation.


\subsubsection{Support vector regression} 
\label{section:SVR}
Support vector regression (SVR) seeks to develop a model~\cite{smola_2004}
\begin{align}\label{eq:SVRone}
	\regressionModel(\features;\weight) =
	\innerprod{\weight}{\boldsymbol{\phi}(\features)}+b,
\end{align}
where $\boldsymbol{\phi}:\RR{\nfeatures}\rightarrow\SVRfeaturespace$,
$\weight\in\SVRfeaturespace$, and $\SVRfeaturespace$ is a (potentially
unknown) feature space
equipped with inner product $\innerprod{\cdot}{\cdot}$.  SVR aims to compute a
`flat' function (i.e., $\innerprod{\weight}{\weight}$
small) that only penalizes prediction errors that exceed a threshold
$\epsilon$ (i.e., soft margin loss).  SVR employs slack variables $\boldsymbol{\xi}$ and
$\boldsymbol{\xi}^\star$  to address deviations exceeding
$\epsilon$ and employs a parameter $C$ to penalize these deviations,
leading to the primal problem~\cite{chang_2011}
\begin{align*}
	\begin{split}
		\underset{\weightDummy,b,\boldsymbol{\xi},\boldsymbol{\xi}^\star}{\mathrm{minimize}}\
		&\frac{1}{2}\innerprod{\weightDummy}{\weightDummy}+C\sum_{i=1}^{\ntrainingData}(\xi_i+\xi_i^\star),\\
		\text{subject to}\
		&\errorTrain{i}-\innerprod{\weightDummy}{\boldsymbol{\phi}(\featuresTrain{i})}-b \le \epsilon +
		\xi_i,\quad i=1,\ldots,\ntrainingData,
\\
&\innerprod{\weightDummy}{\boldsymbol{\phi}(\featuresTrain{i})}+b-\errorTrain{i} {}\le \epsilon +
\xi_i^\star,\quad i=1,\ldots,\ntrainingData,
\\
&\boldsymbol{\xi},\,\boldsymbol{\xi}^\star {}\ge \zero.
	\end{split}
\end{align*}
The corresponding dual problem is
\begin{align*}
\begin{split}
	\underset{\boldsymbol{\alpha},\boldsymbol{\alpha}^\star}{\text{minimize}} \
	&\frac{1}{2}(\boldsymbol{\alpha}-\boldsymbol{\alpha}^\star)^T\mathbf{Q}(\boldsymbol{\alpha}-\boldsymbol{\alpha}^\star)
{}+\epsilon \mathbf{1}^T(\boldsymbol{\alpha}+\boldsymbol{\alpha}^\star)
{}-\boldsymbol{\error}^T(\boldsymbol{\alpha}-\boldsymbol{\alpha}^\star),\\
\text{subject to}\ &\mathbf{1}^T(\boldsymbol{\alpha}-\boldsymbol{\alpha}^\star){}=0,
\\
&0 \le \alpha_i,\,\alpha_i^\star{}\le C,\quad i=1,\ldots,\ntrainingData.
\end{split}
\end{align*}
Here, $Q_{ij}\defeq K(\featuresTrain{i},\,\featuresTrain{j})\defeq
\innerprod{\boldsymbol{\phi}(\featuresTrain{i})}{\boldsymbol{\phi}(\featuresTrain{j})}$ and
$\weight=\sum_{i=1}^\ntrainingData(\alpha_i-\alpha_i^\star)\boldsymbol{\phi}(\featuresTrain{i})$.
The resulting model \eqref{eq:SVRone} can be equivalently expressed as
\begin{align}
\regressionModel(\features;\weight) = \sum_{i=1}^{\ntrainingData}
\left(\alpha_i-\alpha_i^\star\right) K(\featuresTrain{i},\features) + b.
\label{eq:SVRtwo}
\end{align}
There exist many kernel functions $K(\featuresTrain{i},\,\featuresTrain{j})$ that
correspond to an inner product in a feature space $\SVRfeaturespace$.
We consider two in this work: the linear kernel
$K(\featuresTrain{i},\,\featuresTrain{j})=\featuresTrain{i}^T\featuresTrain{j}$,
in which case we denote the method as SVR: Linear, and the
Gaussian radial basis function kernel
$K(\featuresTrain{i},\,\featuresTrain{j})=\exp(-\gamma\|\featuresTrain{i}-\featuresTrain{j}\|^2)$,
in which case we denote the method as
SVR: RBF.

SVR: Linear exhibits a similar (low) capacity as the OLS: Linear approach,
while the SVR: RBF method is higher capacity.  Hyperparameters for this
approach are the penalty parameter $C$ and the margin $\epsilon$. The
method SVR: RBF also considers  $\gamma$ to be a hyperparameter.

\subsubsection{Random forest} 
\label{section:RF}
Random forest (RF) regression \cite{breiman_2001} 
uses decision trees constructed by decomposing the feature space along
canonical directions in a manner that greedily minimizes the mean-squared prediction errors over the training data. The prediction generated by a decision tree corresponds to
the average value of the response over the training data that reside within the same
feature-space region as the prediction point. 

Because decision trees are high
capacity (i.e., low bias and high variance),
random forests employ bootstrap aggregating (i.e., bagging)---which is a
variance-reduction mechanism---to reduce
the prediction variance. Here,
$\ntree$ different data sets are generated by sampling the original training
set $\trainingData$ with replacement; a decision tree is subsequently
constructed from each of these training sets, yielding $\ntree$ different
regression functions
$\regressionModelArg{i}$, $i=1,\ldots,\ntree$. The final regression
function corresponds to the average prediction across the ensemble such that 
\begin{align*}
	\regressionModel(\features) = \frac{1}{\ntree}\sum_{1=1}^\ntree \regressionModelArg{i}(\features).
\end{align*}
To decorrelate each of the decision trees and further reduce prediction
variance, random forests introduce another source of randomness: when training
each tree, a random subset of $\nfeaturesSplit\ll \nfeatures$ features is considered for
performing the feature-space split at each node in the tree.

Random forest regression yields a high-capacity model with relatively low
variance due to the variance-reduction mechanisms it employs.  The
hyperparameters for this approach correspond to the number of trees in the
ensemble $\ntree$ and the size of the feature subset $\nfeaturesSplit$
considered for splitting during training.

\subsubsection{$k$-nearest neighbors} 
\label{section:knn}

$k$-nearest neighbors ($k$-NN) produces predictions arising from a
weighted average of the responses corresponding to the $k$-nearest training 
points in feature space:
\begin{align*}
	\regressionModel(\features) = \sum_{i\in\nearestSetArg{k}}
	\knnWeight(\featuresTrain{i},\features)
	\errorTrain{i}.
\end{align*}
Here, $\nearestSetArg{k}\subseteq\{1,\ldots,\ntrainingData\}$ with
$\card{\nearestSetArg{k}}=k(\leq\ntrainingData)$ satisfies
$\|\features-\featuresTrain{i}\|_2\leq 
\|\features-\featuresTrain{j}\|_2$ for all $i\in\nearestSetArg{k}$,
$j\in\{1,\ldots,\ntrainingData\}\setminus \nearestSetArg{k}$.
Choices for the weights $\knnWeight$ include uniform weights:
$\knnWeight=\knnWeightUniform$ with
\begin{align*}
\knnWeightUniform(\featuresTrain{j},\features)\defeq\frac{1}{k}, 
\end{align*}
and weights determined from the Euclidean distance:
$\knnWeight=\knnWeightEuclidean$ with 
\begin{align*}
\knnWeightEuclidean(\featuresTrain{j},\features)\defeq\frac{\left(\sum_{i=1}^{k}
\|\featuresTrain{i}-\features\|_2^{-1}\right)^{-1}}{\|\featuresTrain{j}-\features\|_2}.
\end{align*}

The capacity of $k$-NN increases as the number of nearest neighbors $k$
decreases. Hyperparameters for this approach are the number of nearest
neighbors $k$ and the choice of weights $\knnWeight$.

\subsubsection{Artificial neural network} 
\label{sec:ANN}

A feed-forward artificial neural network (ANN) with multiple layers---also known as a multilayer
perceptron (MLP)~\cite{rumelhart_1985}---generates a regression function by composition such that
\begin{equation}
	\regressionModel(\features;\weight) =
	\nnFuncArg{\nlayers}(\cdot;\weightNNarg{\nlayers})\circ
	\nnFuncArg{\nlayers-1}(\cdot;\weightNNarg{\nlayers-1})\circ
	\cdots\circ\nnFuncArg{1}(\features;\weightNNarg{1}),
\end{equation}
where $\nnFuncArg{i}(\cdot\ourRereading{;\weightNNarg{i}}):\RR{\nNeuron{i-1}}\rightarrow\RR{\nNeuron{i}}$, $i=1,\ldots, \nlayers$
denotes the function applied at layer $i$ of the neural network; 
$\weightNNarg{i}
\in\RR{\nNeuron{i}\times\ourRereading{(\nNeuron{i-1}+1)}}$, $i=1,\ldots,\nlayers$ denotes the weights
employed at
layer $i$; $\nNeuron{i}$ denotes the number of neurons at layer $i$;
and $\weight\equiv[\vectorize{\weightNNarg{1}};\ \ldots\ ;
\vectorize{\weightNNarg{\nlayers}}]$ in this case.  
The input layer ($i=0$) consists of \ourRereading{the features $\features$ as the $\nNeuron{0} = \nfeatures$ neurons, as well as unity}, 
and the final (output) layer ($i=\nlayers$)
produces the prediction
$\regressionModel(\features)$ such that $\nNeuron{\nlayers}=1$.  The
intermediate layers 1 to $\nlayers-1$ are considered hidden and include unity as an input.
Each neuron $j$ in layers 1 to $\nlayers$ applies an activation function 
$\activationfunction$
to
a linear combination of the outputs from the previous layer such that
\begin{equation}
	\nnFuncArg{i}(\mathbf{y}_{i-1};\weightNNarg{i}) =
		h(\weightNNarg{i}\ourRereading{[1;\mathbf{y}_{i-1}]}),
\end{equation}
where $\mathbf{y}_{i-1}\in\RR{\nNeuron{i-1}}$ is the output from the
previous layer, and $\mathbf{y}_{0}\defeq\ourRereading{\features}$.  The activation function is applied element-wise to the vector argument.

\begin{figure}
\centering
\begin{tikzpicture}

\node () at (0,-5) {Input Layer};

\makeatletter
\DeclareRobustCommand{\rvdots}{%
  \vbox{
    \baselineskip4\p@
    \hbox{.}\hbox{.}\hbox{.}
  }}
\makeatother


 \definecolor{color1}{HTML}{F3C300}
 \definecolor{color2}{HTML}{875692}
 \definecolor{color3}{HTML}{F38400}
 \definecolor{color4}{HTML}{A1CAF1}
 \definecolor{color5}{HTML}{BE0032}
 \definecolor{color6}{HTML}{C2B280}
 \definecolor{color7}{HTML}{848482}
 \definecolor{color8}{HTML}{008856}
 \definecolor{color9}{HTML}{E68FAC}
\definecolor{color10}{HTML}{0067A5}
\definecolor{color11}{HTML}{F99379}
\definecolor{color12}{HTML}{604E97}
\definecolor{color13}{HTML}{F6A600}
\definecolor{color14}{HTML}{B3446C}
\definecolor{color15}{HTML}{DCD300}



\draw[draw=color1] (0,-0) ellipse (0.7 and 0.35);
\draw[draw=color2] (0,-1) ellipse (0.7 and 0.35);
\draw[draw=color3] (0,-2) ellipse (0.7 and 0.35);
\draw[draw=color5] (0,-4) ellipse (0.7 and 0.35);

\node (x) at (0,-0) {\color{color1} \strut 1};
\node (x) at (0,-1) {\color{color2} \strut$\phantom{_1}x_1$};
\node (x) at (0,-2) {\color{color3} \strut$\phantom{_2}x_2$};
\node (x) at (0,-3) {\color{color4} $\rvdots$};
\node (x) at (0,-4) {\color{color5} \strut$\phantom{_{N_x}}x_{N_x}^{}$};

\draw[->,>=stealth, line cap=round, draw=color1] (0.7,-0) -- (2.3,-1);
\draw[->,>=stealth, line cap=round, draw=color1] (0.7,-0) -- (2.3,-2);
\draw[->,>=stealth, line cap=round, draw=color1] (0.7,-0) -- (2.3,-3);
\draw[->,>=stealth, line cap=round, draw=color1] (0.7,-0) -- (2.3,-4);
           
\draw[->,>=stealth, line cap=round, draw=color2] (0.7,-1) -- (2.3,-1);
\draw[->,>=stealth, line cap=round, draw=color2] (0.7,-1) -- (2.3,-2);
\draw[->,>=stealth, line cap=round, draw=color2] (0.7,-1) -- (2.3,-3);
\draw[->,>=stealth, line cap=round, draw=color2] (0.7,-1) -- (2.3,-4);
            
\draw[->,>=stealth, line cap=round, draw=color3] (0.7,-2) -- (2.3,-1);
\draw[->,>=stealth, line cap=round, draw=color3] (0.7,-2) -- (2.3,-2);
\draw[->,>=stealth, line cap=round, draw=color3] (0.7,-2) -- (2.3,-3);
\draw[->,>=stealth, line cap=round, draw=color3] (0.7,-2) -- (2.3,-4);
              
\draw[->,>=stealth, line cap=round, draw=color4] (0.7,-3) -- (2.3,-1);
\draw[->,>=stealth, line cap=round, draw=color4] (0.7,-3) -- (2.3,-2);
\draw[->,>=stealth, line cap=round, draw=color4] (0.7,-3) -- (2.3,-3);
\draw[->,>=stealth, line cap=round, draw=color4] (0.7,-3) -- (2.3,-4);
            
\draw[->,>=stealth, line cap=round, draw=color5] (0.7,-4) -- (2.3,-1);
\draw[->,>=stealth, line cap=round, draw=color5] (0.7,-4) -- (2.3,-2);
\draw[->,>=stealth, line cap=round, draw=color5] (0.7,-4) -- (2.3,-3);
\draw[->,>=stealth, line cap=round, draw=color5] (0.7,-4) -- (2.3,-4);

\node () at (3,-5) {Hidden Layer 1};
\draw[draw=color6] (3,-0) ellipse (0.7 and 0.35);
\draw[draw=color7] (3,-1) ellipse (0.7 and 0.35);
\draw[draw=color8] (3,-2) ellipse (0.7 and 0.35);
\draw[draw=color10] (3,-4) ellipse (0.7 and 0.35);
\node (1) at (3,-0) {\color{color6} \strut 1};
\node (x) at (3,-1) {\color{color7} \strut $\phantom{_1^{(1)}}y_1^{(1)}$};
\node (x) at (3,-2) {\color{color8} \strut $\phantom{_2^{(1)}}y_2^{(1)}$};
\node (x) at (3,-3) {\color{color9} $\rvdots$};
\node (x) at (3,-4) {\color{color10}\strut $\phantom{_{p_1}^{(1)}}y_{p_1}^{(1)}$};

\draw[->,>=stealth, line cap=round, draw=color6] (3.7,-0) -- (5.3,-1);
\draw[->,>=stealth, line cap=round, draw=color6] (3.7,-0) -- (5.3,-2);
\draw[->,>=stealth, line cap=round, draw=color6] (3.7,-0) -- (5.3,-3);
\draw[->,>=stealth, line cap=round, draw=color6] (3.7,-0) -- (5.3,-4);
                   
\draw[->,>=stealth, line cap=round, draw=color7] (3.7,-1) -- (5.3,-1);
\draw[->,>=stealth, line cap=round, draw=color7] (3.7,-1) -- (5.3,-2);
\draw[->,>=stealth, line cap=round, draw=color7] (3.7,-1) -- (5.3,-3);
\draw[->,>=stealth, line cap=round, draw=color7] (3.7,-1) -- (5.3,-4);
                
\draw[->,>=stealth, line cap=round, draw=color8] (3.7,-2) -- (5.3,-1);
\draw[->,>=stealth, line cap=round, draw=color8] (3.7,-2) -- (5.3,-2);
\draw[->,>=stealth, line cap=round, draw=color8] (3.7,-2) -- (5.3,-3);
\draw[->,>=stealth, line cap=round, draw=color8] (3.7,-2) -- (5.3,-4);
               
\draw[->,>=stealth, line cap=round, draw=color9] (3.7,-3) -- (5.3,-1);
\draw[->,>=stealth, line cap=round, draw=color9] (3.7,-3) -- (5.3,-2);
\draw[->,>=stealth, line cap=round, draw=color9] (3.7,-3) -- (5.3,-3);
\draw[->,>=stealth, line cap=round, draw=color9] (3.7,-3) -- (5.3,-4);
                 
\draw[->,>=stealth, line cap=round, draw=color10] (3.7,-4) -- (5.3,-1);
\draw[->,>=stealth, line cap=round, draw=color10] (3.7,-4) -- (5.3,-2);
\draw[->,>=stealth, line cap=round, draw=color10] (3.7,-4) -- (5.3,-3);
\draw[->,>=stealth, line cap=round, draw=color10] (3.7,-4) -- (5.3,-4);

\node () at (6,-5) {Hidden Layer 2};
\draw[draw=color11] (6,-0) ellipse (0.7 and 0.35);
\draw[draw=color12] (6,-1) ellipse (0.7 and 0.35);
\draw[draw=color13] (6,-2) ellipse (0.7 and 0.35);
\draw[draw=color15] (6,-4) ellipse (0.7 and 0.35);

\node (1) at (6,-0) {\color{color11} \strut 1};
\node (x) at (6,-1) {\color{color12} \strut $\phantom{_1^{(2)}}y_1^{(2)}$};
\node (x) at (6,-2) {\color{color13} \strut $\phantom{_2^{(2)}}y_2^{(2)}$};
\node (x) at (6,-3) {\color{color14} $\rvdots$};
\node (x) at (6,-4) {\color{color15} \strut $\phantom{_{p_2}^{(2)}}y_{p_2}^{(2)}$};

\draw[->,>=stealth, line cap=round, draw=color11] (6.7,-0) -- (8.3,-2);
\draw[->,>=stealth, line cap=round, draw=color12] (6.7,-1) -- (8.3,-2);
\draw[->,>=stealth, line cap=round, draw=color13] (6.7,-2) -- (8.3,-2);
\draw[->,>=stealth, line cap=round, draw=color14] (6.7,-3) -- (8.3,-2);
\draw[->,>=stealth, line cap=round, draw=color15] (6.7,-4) -- (8.3,-2);

\node () at (9,-3) {Output Layer};
\draw[rotate=0] (9,-2) ellipse (0.7 and 0.35);
\node (x) at (9,-2) {\strut $\regressionModel(\features;\weight)$};

\end{tikzpicture}
\caption{Feedforward artificial neural network with two hidden layers.}
\label{fig:mlp}
\end{figure}
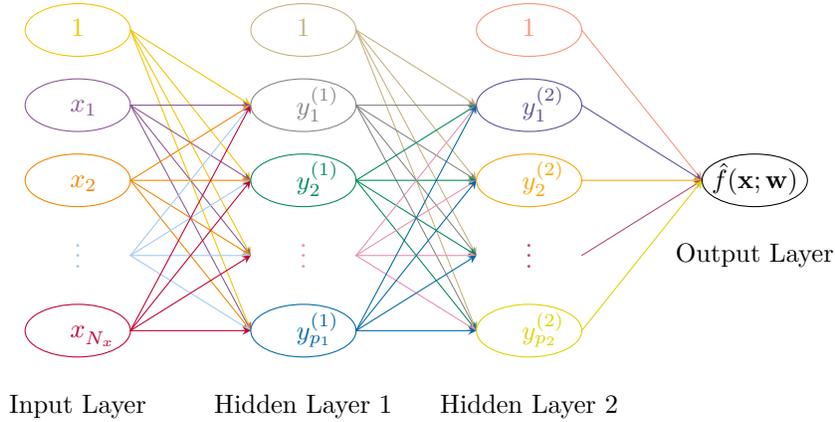

Common activation functions include
\begin{alignat*}{4}
\text{Identity:} \quad & \activationfunction(y) &&{}= y,
\\
\text{Logistic Sigmoid:} \quad & \activationfunction(y) &&{}= \frac{1}{1+e^{-y}},
\\
\text{Hyperbolic Tangent:} \quad & \activationfunction(y) &&{}= \tanh y,
\\
\text{Rectified Linear Unit:} \quad & \activationfunction(y) &&{}= \max\{0,\,y\}.
\end{alignat*}
For regression problems, the output layer typically employs an identity
activation function.
Figure~\ref{fig:mlp} provides a diagram of a feedforward artificial neural network with two hidden layers.  
%

To train the neural network, the weights $\weight$ are often computed by minimizing
the mean-squared error over the training data according to
Eq.~\eqref{eq:ols_constraint}. Because neural networks tend to be high
capacity, regularization is often employed; this work considers a ridge
formulation wherein the weights are the solution to the minimization problem
\begin{align}
	\underset{\weightDummy}{\mathrm{minimize}}\
	\sum_{i=1}^\ntrainingData\left(\regressionModel(\featuresTrain{i};\weightDummy) -
	\errorTrain{i}\right)^2 + \alpha\|\weightDummy\|_2^2,
\label{eq:ols_penalty}
\end{align}
where $\alpha\in\RRplus{}$ is a regularization penalty term.

The capacity of the ANN increases as the number of layers and number of
neurons per layer increase, and as the regularization parameter $\alpha$ decreases.
Hyperparameters for this approach include the neural network architecture
(i.e., number of layers $\nlayers$ and number of neurons at each hidden layer
$\nNeuron{i}$, $i=1,\ldots, \nlayers-1$), the choice of
activation function $h$, and the penalty parameter $\alpha$.

\section{Numerical experiments} 
\label{sec:experiments}

This section demonstrates the methodology proposed in Section~\ref{sec:approach} 
on several computational-mechanics examples.

\subsection{Machine-learning error models} 
The numerical experiments assess all feature-engineering methods described in
Section \ref{subsec:indicators}, as well as all regression techniques
introduced in Section \ref{subsection:regressionfunctionapprox}. We now
describe details related to the construction of the error models.

\subsubsection{Preprocessing} 
Prior to training and applying the regression techniques, the training and
test data are shifted and scaled as follows.  
First, for Feature-Engineering Methods
\ref{feat:paramres}--\ref{feat:paramsampleres}, the elements of $\resApprox$
for which the variance of the training data is zero are removed.  
For Feature-Engineering Methods
\ref{feat:paramresPCA} and \ref{feat:paramgappyPODres}, we perform
PCA after subtracting the training mean from the training data.
For all feature-engineering methods, once the features are computed, each
feature is standardized by subtracting the training mean subsequently dividing
by the training variance.

\subsubsection{Settings, hyperparameters, and cross-validation} 
\label{sec:s_hp_cv}
We employ \texttt{scikit-learn}~\cite{scikit-learn,sklearn_api} (with default
options unless otherwise specified) to construct the regression-function
approximations described in Section \ref{subsection:regressionfunctionapprox}. 

For regression methods with hyperparameters, we perform a grid search on the
\reviewerOne{cross-validation} grids specified below using five-fold cross-validation with data shuffling on
the training data.  The coefficient of determination $r^2$ computed on the
validation sub-sample is used to score the hyperparameter combinations.  The
hyperparameter values yielding the highest average score across all folds are
selected, and the final model is trained using these parameter values on the
entire training set.

\paragraph{OLS: Quadratic} 
The OLS methods described in Section \ref{section:OLS} have no
hyperparameters. However, the number of effective features for OLS: Quadratic
$(\nfeatures+1)(\nfeatures+2)/2$ is quite large for many feature-engineering
methods.  To keep this number tractable, if $\nfeatures>100$, we perform a
univariate $F$-test (as described for Feature-Engineering Method
\ref{feat:paramgappyPODres} in Section \ref{subsec:indicators}) on the
original features $\features$ to determine the 100 most important features to
consider, resulting in 5151 effective features. 

\paragraph{SVR: Linear} 
The SVR techniques described in Section \ref{section:SVR} have two
hyperparameters:
the penalty parameter $C$ and margin $\epsilon$. We
employ cross-validation grids of $C =
10^n$, $n\in\{-2,\,-1,\,\hdots,\,4\}$, and $\epsilon = 10^n$, 
$n\in\{-3,\,-2,\,\hdots,\,0\}$ for these. SVR: RBF has an additional
hyperparameter, which is the parameter $\gamma$; we employ a
cross-validation grid of 
$\gamma = 10^n$, $n\in\{-5,\,-4,\,\hdots,\,1\}$ for this parameter.

\paragraph{RF}
As described in Section \ref{section:RF}, the hyperparameters for RF include
the number of trees in the
ensemble $\ntree$ and the size of the feature subset $\nfeaturesSplit$
considered for splitting during training. For these parameters, we employ
cross-validation grids of 
$\ntree \in\{25,\,50,\,\hdots,\,150\}$ and $\nfeaturesSplit\in\{n,\,\sqrt{n},\,\log_2 n\}$.

\paragraph{$k$-NN}
As discussed in Section \ref{section:knn}, the hyperparameters for $k$-NN
include
the number of nearest
neighbors $k$ and the choice of weights $\knnWeight$; we employ
cross-validation grids of 
$k\in
\{1,\,2,\,\hdots,\,\min(10,\,\frac{4}{5}\ntrainingData)\}$
and $\knnWeight\in\{\knnWeightUniform,\knnWeightEuclidean\}$.
In addition, to make distance computation less computationally expensive, we
perform a univariate $F$-test (as described for Feature-Engineering Method
\ref{feat:paramgappyPODres} in Section \ref{subsec:indicators}) to determine
the most important features. We consider the number of retained features to
be a hyperparameter with \reviewerOne{cross-validation} grid $\{1,\,2,\,\hdots,\,\min(10,\,\nfeatures)\}$.
We do this for all feature-engineering methods except for Feature-Engineering Methods
\ref{feat:paramresPCA} and \ref{feat:paramgappyPODres}, in which case the
number of residual principal components is varied instead (see ``Residual
principal components'' below).

\paragraph{ANN} For ANN, the activation function is used as a hyperparameter,
where considered activation functions are identity, logistic sigmoid,
hyperbolic tangent, and rectified linear unit. In addition, the $L^2$ penalty
regularization term $\alpha = 10^n$ is a hyperparameter; we employ a
cross-validation grid of $n\in\{-8,\,-6,\,\hdots,\,0\}$.  One hidden layer is
considered with 100 neurons.  The solution is obtained using a limited-memory
BFGS
optimization algorithm with a tolerance of $10^{-5}$ and a maximum allowable iteration count of 1000.

\paragraph{Residual principal components}
Feature-Engineering Methods \ref{feat:paramresPCA} and
\ref{feat:paramgappyPODres} employ features
$\FeatureParResPCA$ and $\FeatureParResGappyPCA$, which are equipped with
their own hyperparameter: the number of residual principal components
$\dimresApprox$. When these features are employed with any regression
technique, we employ a cross-validation grid  of
$\dimresApprox\in\{1,\,2,\,3,\,4,\,5,\,10,\,15,\,20,\,25,\,30\}$. For
Feature-Engineering Method \ref{feat:paramgappyPODres}, we consider only
values of $\dimresApprox$
that satisfy $\dimresApprox\leq\nsamples$ to ensure that 
$(\sampleMat\basisresApprox)^+\sampleMat\basisresApprox=\identityArg{\dimresApprox}$, which is required by gappy POD.

\subsubsection{Performance metrics} 
\label{sec:performanceMetrics}

To quantify the performance of different feature-engineering methods and
machine-learning regression techniques, we compute two quantities on the test data: (1) the mean squared
error (MSE) and (2) the fraction of variance unexplained (FVU).  The test MSE is
defined by
\begin{align*}
\text{MSE}\defeq 
\frac{1}{\ntestData}\sum_{i=1}^{\ntestData}
(\errorTest{i}- \regressionModel(\featuresTest{i}))^2.
\end{align*}
When \ref{data:unique} is employed, the reported test MSE corresponds to the
average test MSE value computed over all $\napprox$ approximate solutions.
Note that we employ $\varianceEstimate =\text{MSE}$ as the variance
in the noise approximation 
$\noiseModel\sim\normal{0}{\varianceEstimate}$;
see \ref{step:noiseApproximation} in Section \ref{subsec:overview}. Thus, a
small MSE implies a smaller variance in the noise approximation and lower
epistemic uncertainty.
The test FVU is defined as 
\begin{align*}
\text{FVU}\defeq \frac{\text{MSE}}{\text{Var}(\errorTestNo)} = 1-r^2,
\end{align*}
with $r^2$ the coefficient of determination (i.e., the fraction of variance
explained).

\subsection{Cube: reduced-order modeling} 
\label{subsec:cube}
The first experiment considers the mechanically induced deformation of a cube
with an approximate solution provided by a Galerkin reduced-order model as
described in Section \ref{section:ROM}.
We conduct this experiment using Albany, an implicit, unstructured grid,
finite-element code for the solution and analysis of multiphysics
problems~\cite{albany_2016}\reviewerTwo{, with eight-node hexahedral elements}.
\reviewerOne{
We solve all systems of nonlinear algebraic equations using a damped Newton's
method, i.e., Newton's method with a fixed step size less than or equal to one.
}
\FloatBarrier

\subsubsection{Overview} 
\label{subsec:cube_overview}
\begin{figure}
\centering
\begin{subfigure}[b]{.49\textwidth}
\includegraphics[scale=.3]{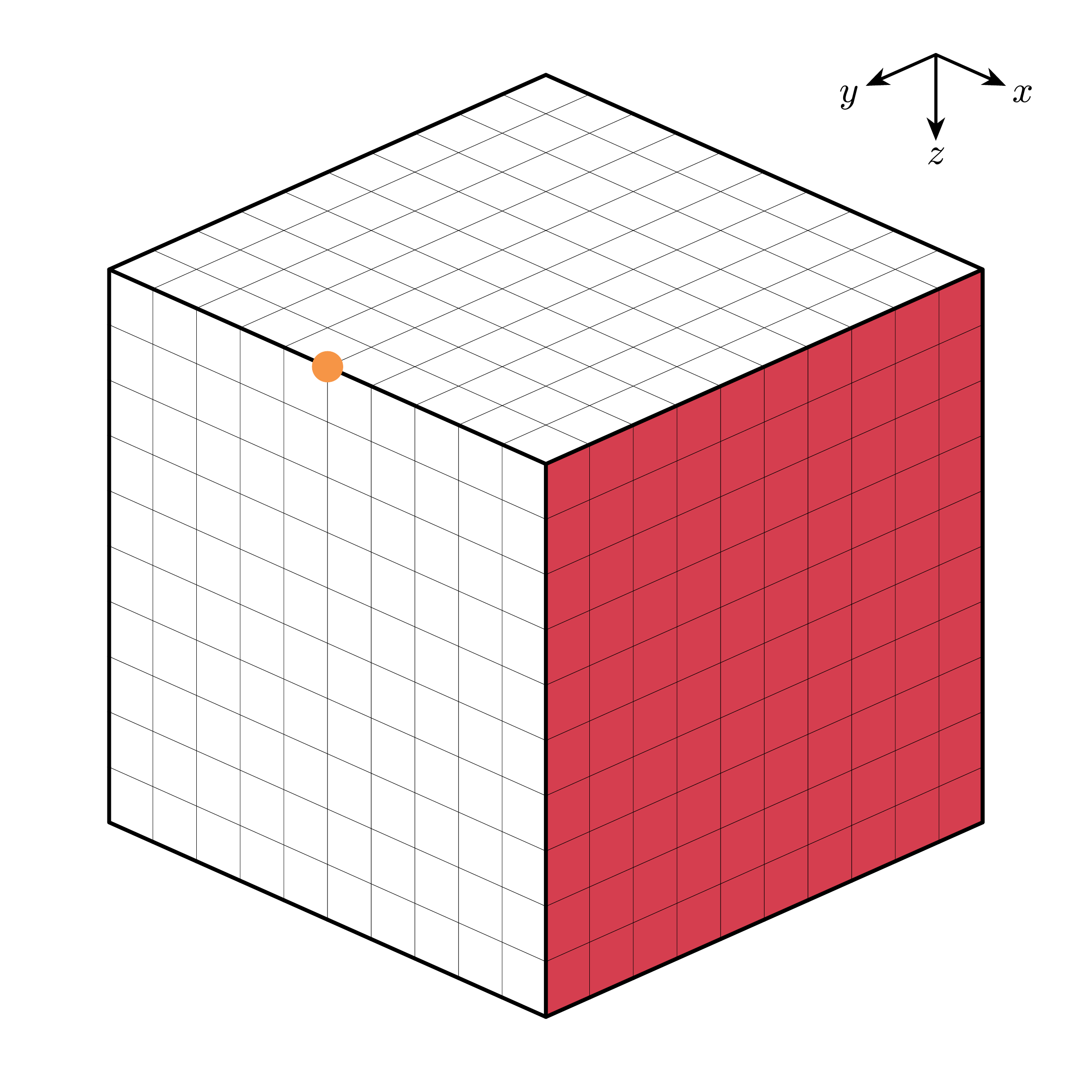}
\caption{View 1}
\label{fig:complex_cube/mesh_plots/mesh_bc_a}
\vspace{-0.4em}
\end{subfigure}
\begin{subfigure}[b]{.49\textwidth}
\includegraphics[scale=.3]{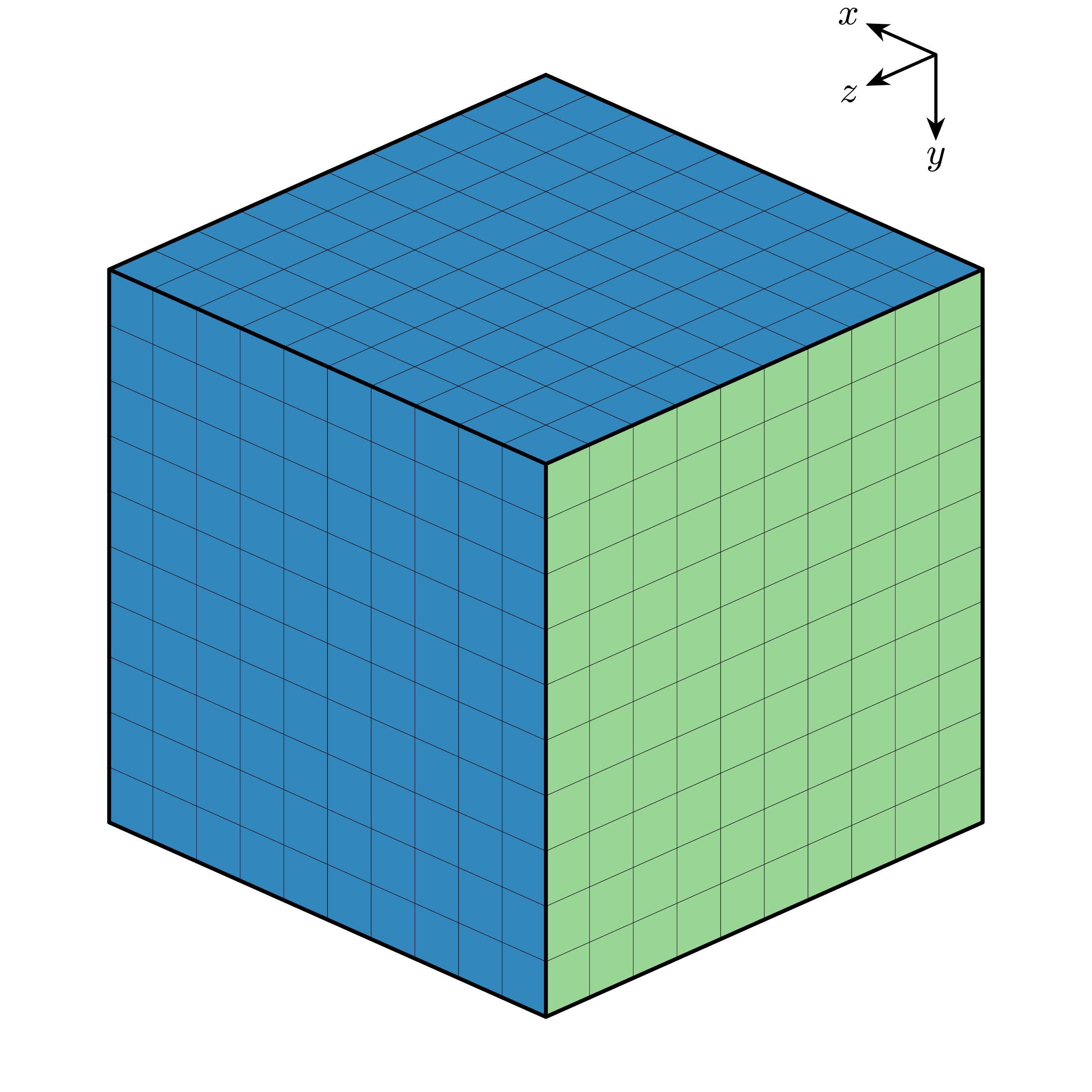}
\caption{View 2}
\label{fig:complex_cube/mesh_plots/mesh_bc_b}
\vspace{-0.4em}
\end{subfigure}
\caption{Cube: Mesh with boundary conditions and node of interest: red denotes
applied traction (Neumann boundary condition), blue denotes
planar-displacement constraint
(Dirichlet boundary condition), green denotes zero-displacement constraint (Dirichlet
boundary condition), and orange denotes the node of interest.}
\label{fig:complex_cube/mesh}
\end{figure}

The undeformed cube is one cubic meter in size with domain $[0\,\text{m}, 1\,\text{m}]^3$.  As
shown in Figure~\ref{fig:complex_cube/mesh}, the cube is discretized using 10
elements along each edge. This discretization is deliberately coarse to enable
computational tractability of the
dual-weighted-residual computations discussed in Section~\ref{subsubsec:local_errors}.

A traction of magnitude $t$ is applied as a Neumann boundary condition to the cube face with an
outward normal in the positive $x$-direction.  The nodes on the opposite face,
with an outward normal in the negative $x$-direction, are constrained to zero 
$x$-, $y$-, and $z$-displacements by homogeneous Dirichlet boundary
conditions.  The nodes on the faces with outward normals in the negative $y$-
and positive $z$-directions are constrained to planar motion along their
respective planes by homogeneous Dirichlet boundary conditions.  The resulting
dimension of the model is $\nstate=3410$.

The node of interest is
located
midway along the edge shared by the faces with outward normals in the positive
$y$- and negative $z$-directions at coordinate ($\frac{1}{2}$ m, 1 m, 0 m) when the
cube is undeformed.  Since the displacement in the $z$-direction of the node
of interest is equal to the displacement in the negative $y$-direction, the
displacements in the $x$- and $y$-directions, denoted by $u_x$ and $u_y$, are the
quantities of interest; that is, $\qoi = u_x$ and $\qoi=u_y$.  We denote their
errors by $\errorQoiArg{u_x}$ and $\errorQoiArg{u_y}$, respectively.
Figure~\ref{fig:complex_cube/mesh} identifies
the boundary conditions and node of interest.

This experiment considers $\nparams=3$ parameters.
The elastic modulus $\paramArg{1}=E$ is varied
between 75.0 and 125.0 GPa, the Poisson ratio $\paramArg{2}=\nu$ is varied between 0.200
and 0.350, and the traction $\paramArg{3}=t$ is varied between 40.0 and 60.0
GPa, resulting in a parameter domain $\paramDomain = \left[
75.0\,\text{GPa}, 125.0\, \text{GPa}\right]\times
\left[
0.200, 0.350
\right]\times \left[40.0\,\text{GPa},60.0\,\text{GPa}\right]
	$.  Across
the parameter domain $\paramDomain$, quantities of interest $u_x$ and $u_y$
vary between limits of 0.242 m and 0.906 m, and
$-0.0814$ m and $-0.237$ m, respectively.

\begin{figure}
\centering
\begin{subfigure}[b]{.49\textwidth}
\includegraphics[scale=.3]{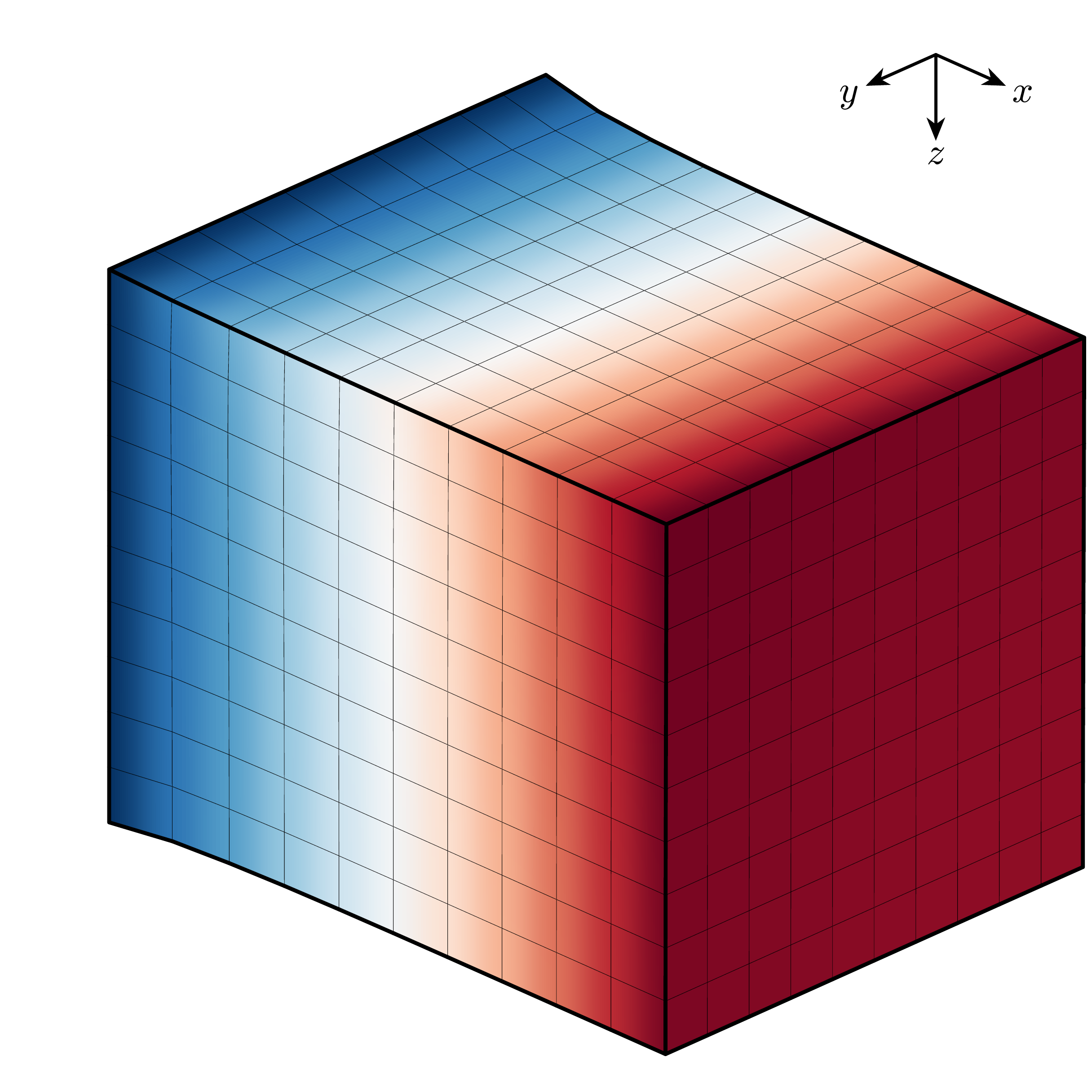}
\caption{First \reviewerTwo{POD basis vector}}
\label{fig:complex_cube/mesh_plots/basis1}
\vspace{-0.4em}
\end{subfigure}
\begin{subfigure}[b]{.49\textwidth}
\includegraphics[scale=.3]{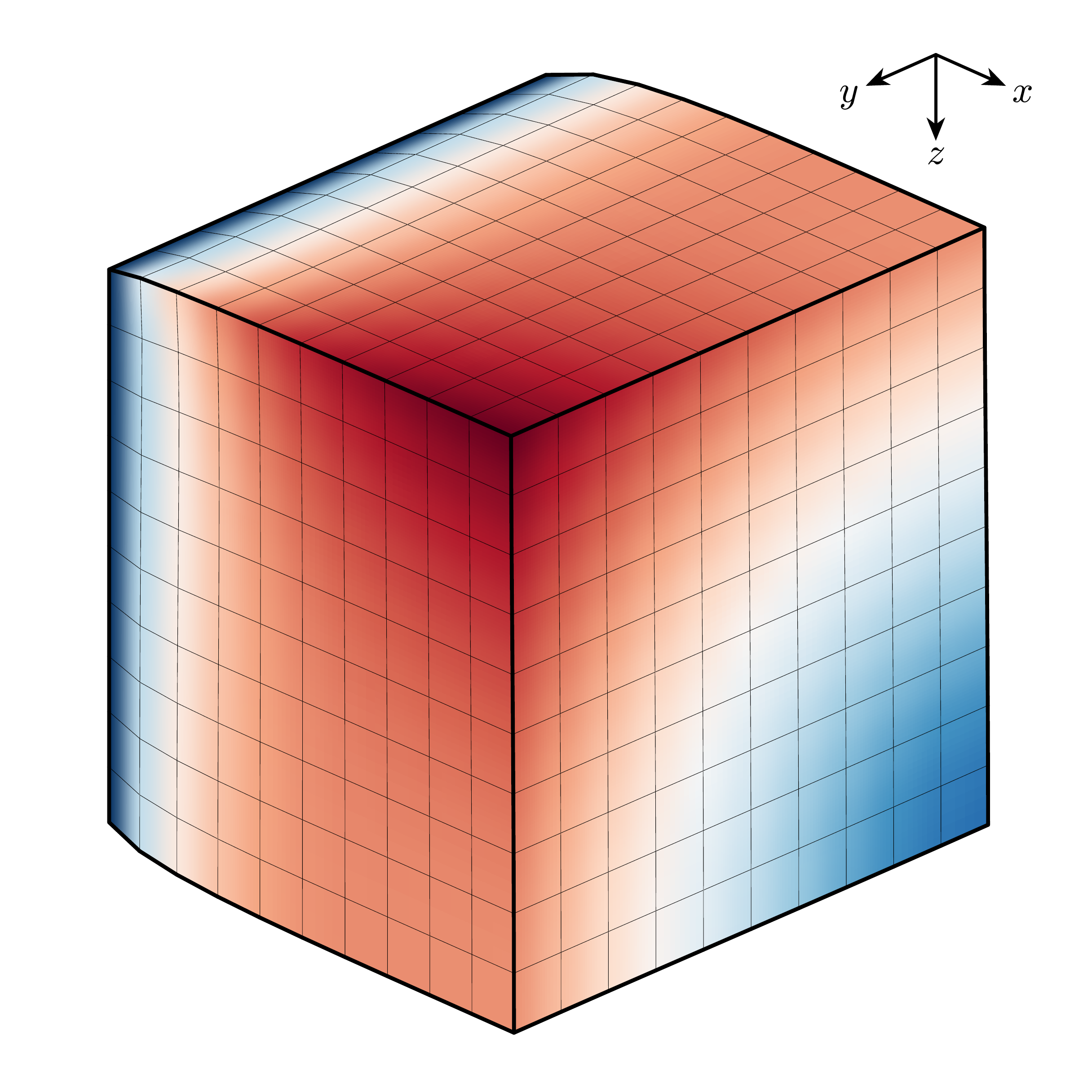}
\caption{Second \reviewerTwo{POD basis vector}}
\label{fig:complex_cube/mesh_plots/basis2}
\vspace{-0.4em}
\end{subfigure}
\caption{Cube: POD basis vectors.  Red indicates larger deformations; blue
indicates smaller deformations.}
\label{fig:complex_cube/basis_functions}
\end{figure}

To construct the requisite snapshots, we first compute the high-fidelity-model
solution $\state(\params)$ for $\params\in\paramDomainPOD\subset \paramDomain$,
where $\paramDomainPOD$ comprises eight Latin hypercube samples 
such that $\nsnapshots=8$. We subsequently compute 
the trial-basis matrix $\trialbasis$  using proper orthogonal
decomposition (POD) by applying Algorithm~\ref{alg:pod} of \ref{app:pod} with inputs
$\{\state(\params)\}_{\params\in\paramDomainPOD}$.
Due to the simplicity of this problem, the first POD vector captures 95.87\% of
the statistical energy such that $\energyThresholdArg{1}=0.9587$ (see output of Algorithm \ref{alg:pod}), and the
second captures 3.62\% such that $\energyThresholdArg{2} = 0.9949\%$. 
Figure~\ref{fig:complex_cube/basis_functions} depicts
the first two POD \reviewerTwo{basis} vectors.

One set of 100 parameter instances is randomly sampled from the parameter domain
$\paramDomain$ to serve as the set of \pti{} $\paramDomainTrain$.  To assess
method performance for smaller amounts of training data,
we randomly create nested training sets from this set. Another set of 100
parameter instances is
randomly sampled to serve as the set of
parameter testing instances $\paramDomainTest$.
We consider approximate solutions generated by a Galerkin reduced-order model
(i.e., $\testbasis = \trialbasis$)
with $\stateRef=\zero$ using 
$\napprox=2$ different
basis dimensions: $\nbasisArg{1}=1$ and $\nbasisArg{2}=2$.
We consider both \ref{data:pooled} (pooled) and \ref{data:unique} (unique) as described in Section \ref{subsection:trainingdata}.
The former case yields $\ntestData=200$ testing data points and up to
$\ntrainingData=200$ training points, while the latter yields
$\ntestData=100$ testing points and up to
$\ntrainingData=100$ testing points.

When using Feature-Engineering Methods~\ref{feat:paramgappyPODres} and \ref{feat:paramsampleres},
three numbers of sample points are used: $\nsamples\in\{10,\, 100,\,1000\}$.

\subsubsection{Results} 
\label{cube_results}
\newif\ifshowlocal
\showlocaltrue
\newif\ifhasgauss
\hasgausstrue
\newif\ifispcap
\ispcapfalse
\def \casetxt {complex_cube}
\def \casetex {Cube}

\def \vartxta {ex}
\def \vartexa {$\errorQoiArg{u_x}$}
\def \vartxtb {ey}
\def \vartexb {$\errorQoiArg{u_y}$}
\def \vartxtc {epsilon}
\def \vartexc {$\errorState$}
\def \nvar    {3}
\def\ntrain/{100}
\def\kpercentbetter/{36.51\%}
\def\qpercentbetter/{61.90\%}
\def\extratext/{Instances for which $\text{MSE}_q=\text{MSE}_k$ are not
plotted and comprise 1.59\% of cases.}
\def\localtex/{\ref{data:unique}}
\def\globaltex/{\ref{data:pooled}}

We first assess the difference in performance between the $k$- and
$q$-sampling approaches for computing the sampling matrix $\sampleMat$
employed by Feature-Engineering Methods \ref{feat:paramgappyPODres} and \ref{feat:paramsampleres}.
Figure~\ref{fig:\casetxt/kq_histogram} provides a comparison of the test MSEs
that arise from using $k$- and $q$-sampling with $\card{\paramDomainTrain}=\ntrain/$ parameter training instances.  
This figure is generated from 
252 total data points, which aggregate test MSE values over three errors: 
$\errorQoiArg{u_x}$, $\errorQoiArg{u_y}$, and $\errorState$;
Feature-Engineering Methods 
\ref{feat:paramgappyPODres} and \ref{feat:paramsampleres};
 the seven
regression techniques discussed
in~\ref{subsection:regressionfunctionapprox}; three numbers of sample points
$\nsamples\in\{10, 100, 1000\}$; and the two data-set approaches discussed in Section \ref{subsection:trainingdata}.
%
Overall, $q$-sampling outperforms $k$-sampling in this case;
therefore, the remaining results for Feature-Engineering Methods
\ref{feat:paramgappyPODres} and \ref{feat:paramsampleres} for this experiment
consider only $q$-sampling.

For each of the ML regression techniques, Figure~\ref{fig:\casetxt/best_ind} reports how the
test FVU varies with respect to the number of \pti{}
$\card{\paramDomainTrain}$
when using the best-performing feature-engineering method for each regression technique.  
Generally, as the number of \pti{} 
$\card{\paramDomainTrain}$
increases, the test FVU decreases,
and the best regression techniques are those that enable higher capacity:
ANN, SVR: RBF, and OLS: Quadratic. On the other hand, performance of the low-capacity regression
methods (OLS: Linear and SVR: Linear) saturates as the amount of training data increases.  This saturation can be
attributed to the greater structure they enforce, thereby leading to a high
bias that produces errors that cannot be reduced with additional training data.  Nonetheless, for
a small amount of training data, these low-capacity methods perform more competitively, as
the higher capacity regression methods overfit the data in this case. RF and
$k$-NN exhibit similar performance; they perform poorly with a small amount of training due to their high capacity.
\ref{data:pooled} and \ref{data:unique} perform similarly in this case.
Overall, the test FVU values are fairly small, and all techniques work
relatively well with a modest number of parameter training instances.  Note that predicting
$\errorState$ is more challenging than predicting $\errorQoiArg{u_x}$ or
$\errorQoiArg{u_y}$; the low-capacity regression methods 
(OLS: Linear and SVR: Linear)
fail to reduce the test FVU below 0.9 when using \ref{data:pooled} as reported in Figure~\ref{fig:\casetxt/best_ind_epsilon_global}.


Whereas Figure~\ref{fig:\casetxt/best_ind} compares each regression technique when
using the best feature-engineering method for that technique,
Figure~\ref{fig:\casetxt/best_ml} compares each feature-engineering method when using the
best regression technique for that method.
It is immediately clear that using the feature $\FeatureResNorm$ alone yields
the highest test FVU, which does not improve as the amount of training data increases.
This is expected since $\FeatureResNorm$ is a single feature of relatively low quality.
For responses $\errorQoiArg{u_x}$ and $\errorQoiArg{u_y}$,
features 
$\dualweightedresidualArg{u_x}$ and $\dualweightedresidualArg{u_y}$ tend to perform the best in the presence of limited training data.
This is expected since these correspond to single high-quality features.
As with Figure~\ref{fig:\casetxt/best_ind}, predicting $\errorState$ is more
challenging than predicting $\errorQoiArg{u_x}$ or $\errorQoiArg{u_y}$; in this case, 
the feature-engineering methods perform similarly with the exception of
$\FeatureResNorm$.
In Figures~\ref{fig:\casetxt/best_ml_ex_local}
and~\ref{fig:\casetxt/best_ml_ey_local}, features $\FeatureParResSampled$ and
$\FeatureParResGappyPCA$ typically outperform the commonly used features $\FeaturePar$ by roughly
an order of magnitude.


Figure~\ref{fig:\casetxt/matrix_t\ntrain/} reports the test FVU for each
combination of feature-engineering method, regression technique, and data-set
approach for $\card{\paramDomainTrain}=\ntrain/$ \pti{}.
Features $\FeaturePar$ require SVR: RBF or ANN to perform well.
Regression methods SVR: RBF and ANN consistently yield the best performance, whereas OLS: Quadratic yields very inconsistent performance.
Once again, there is not a clear performance discrepancy between the two data-set approaches.
Generally, features $\FeatureParResGappyPCA$ perform better than features $\FeatureParResSampled$; this is particularly clear in Figures~\ref{fig:\casetxt/matrix_ex_global_t\ntrain/}, \ref{fig:\casetxt/matrix_ey_global_t\ntrain/}, and~\ref{fig:\casetxt/matrix_ey_local_t\ntrain/}. 
Because features $\FeatureResNorm$, $\dualweightedresidualArg{u_x}$, and
$\dualweightedresidualArg{u_y}$ correspond to only a single feature, their
performances are nearly insensitive to the chosen regression technique, as all
regression techniques yield similarly low capacities for a small number of features.  As
demonstrated previously, feature $\FeatureResNorm$ yields large values of the test
FVU.
For this experiment, full sampling (i.e., $\nsamples=\nstate$) and subsampling
with $\nsamples=1000$ yield little benefit to significant subsampling with $\nsamples=10$; this implies that excellent performance can be obtained
using only a small number of cheaply computable features with the proposed
methodology.


For $\card{\paramDomainTrain}=\ntrain/$ \pti{}, Figure~\ref{fig:\casetxt/exact_v_pred} compares the
values of the responses $\errorQoiArg{u_x}$, $\errorQoiArg{u_y}$, and
$\errorState$ predicted by various regression methods with their exact values over the test set. The figure reports results for 
conventional feature choices (i.e., $\FeatureResNorm$; $\FeaturePar$; 
dual-weighted residual, where applicable), as well as $\FeatureParResGappyPCA$
with only $\nsamples=10$ sampled residual elements, where performance of the best
regression technique for each of these feature-engineering methods is reported.
%
%
In each case, 
$\FeatureParResGappyPCA$ performs better
than all conventional approaches, with $r^2>0.996$ in every case.  We note that
although commonly used features $\FeaturePar$ with ANN and SVR: RBF regression
perform relatively well in this case, this does not occur in
subsequent experiments.

\reviewerThree{
Figure~\ref{fig:\casetxt/hist_gauss} assesses the accuracy of the noise model
$\noiseModel\sim\normal{0}{\varianceEstimate}$ computed in
\ref{step:noiseApproximation} in Section \ref{subsec:overview} for $\FeatureParResGappyPCA$
with $\nsamples=10$ using the best-performing regression technique. In particular,
the figure reports the standard normal distribution compared to a histogram of
the prediction errors, which have been standardized according to the
hypothesized distribution. These prediction errors are computed 
on a second set of test data
$\testDataTwo\defeq\{(\errorTestTwo{i},\featuresTestTwo{i})\}_{i=1}^{\card{\testData}}$,
which is constructed in an identical way to the the first set of test data
$\testData$ used to train the noise model, but is independent of $\testData$
and the training data $\trainingData$ (used to train the regression model).
Table~\ref{tab:\casetxt/validation_frequencies} lists the corresponding
validation frequencies $\frequency_\error(\frequency)$ for multiple
$\frequency$-prediction intervals, where
\begin{align}
	\frequency_\error(\frequency) \defeq
	\card{\left\{(\error,\features)\in\testDataTwo\ |\ {(\error-
	\regressionModel(\features))}/{\varianceEstimate}\in
	C(\frequency)\right\}}/\card{\testDataTwo}
\label{eq:validation_frequency}
\end{align}
and the standard $\frequency$-prediction interval is
\begin{align}
C(\frequency) \defeq \left[-\sqrt{2}\erfinv(\frequency),\,\sqrt{2}\erfinv(\frequency)\right].
\label{eq:prediction_interval}
\end{align}

\begin{table}[ht!]
\centering
\reviewerThree{
\begin{tabular}{c c c c}
\toprule
& \multicolumn{3}{c}{$\frequency_\error(\frequency)$}\\
\cmidrule(l){2-4}
$\frequency$ & $\errorQoiArg{u_x}$ & $\errorQoiArg{u_y}$ & $\errorState$ \\
\midrule
0.80 & 0.96 & 0.86 & 0.96 \\
0.90 & 0.97 & 0.92 & 0.97 \\
0.95 & 0.97 & 0.93 & 0.98 \\
0.99 & 0.98 & 0.96 & 0.98 \\
\bottomrule
\end{tabular}
}
\vspace{-.5em}
\caption{\casetex: Validation frequencies for models using $\FeatureParResGappyPCA$
($\nsamples=10$) with SVR: RBF, SVR: RBF, and ANN, respectively for $\errorQoiArg{u_x}$, $\errorQoiArg{u_y}$, and $\errorState$.}
\label{tab:\casetxt/validation_frequencies}
\end{table}

These results show that the data are not quite Gaussian and are characterized
by heavier tails than would be predicted by a Gaussian distribution. This motivates the
need for a more sophisticated error model to accurately model non-Gaussian
noise, which is the subject of future work. However, the validation
frequencies show that---while all the prediction intervals of the noise model
are not accurate---a subset are accurate. In particular, the 0.99-prediction
interval, the 0.95-prediction interval, and the 0.99-prediction interval are
reasonably accurate for the models associated with
responses $\errorQoiArg{u_x}$, $\errorQoiArg{u_y}$, and $\errorState$,
respectively. Thus, despite the non-Gaussian behavior of the error, these
prediction intervals can still be used to make accurate statistical
predictions.
}

\def\scalefactor/{0.5}

\def\lefttrim/{0in}
\def\righttrim/{0in}
\begin{figure}
   \centering
   \includegraphics[scale=\scalefactor/,clip=true,trim=\lefttrim/ 0in \righttrim/ 0in]{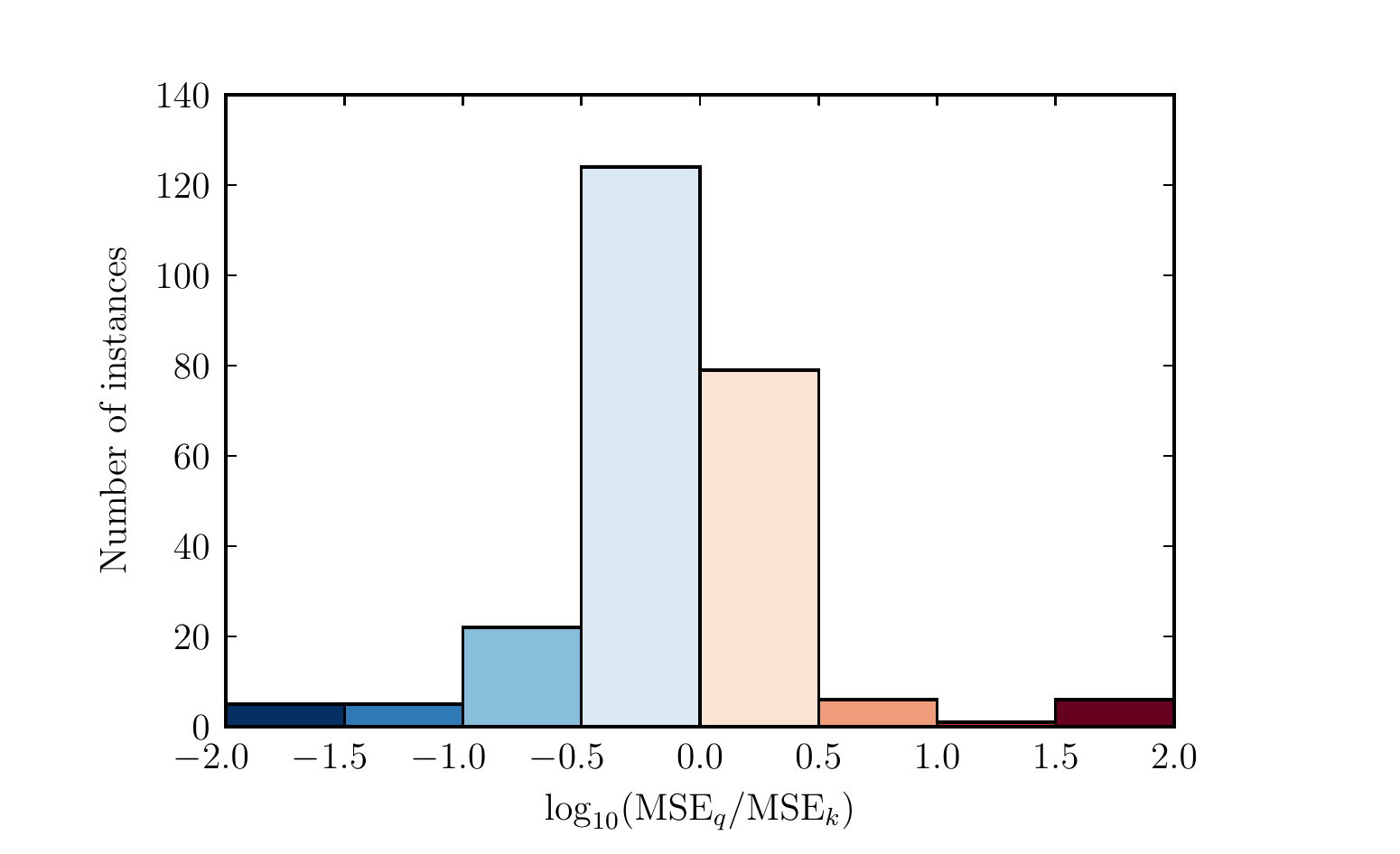}
   \vspace{-0.4em}
   \caption{\casetex: Comparison of MSE using $k$- and $q$-sampling approaches
	 with $\card{\paramDomainTrain}=\ntrain/$ training parameter instances.
 When $\text{MSE}_q/\text{MSE}_k<1$, $q$-sampling outperforms $k$-sampling;
 these cases are depicted in blue and comprise \qpercentbetter/ of cases. 
The converse holds when $\text{MSE}_q/\text{MSE}_k>1$; these cases are
depicted in red and comprise \kpercentbetter/ of cases. \extratext/
 }
   \label{fig:\casetxt/kq_histogram}
\end{figure}

\ifispcap
\def\lefttrim/{0in}
\def\righttrim/{0in}
\def\traininfo/{_t\ntrain/}
\def \descriptiontxt{time_matrix}
\def\descriptiontex/{\reviewerThree{Comparison of wall time, in seconds, required to train regression models offline, including cross-validation (Steps \ref{cv_loop}--\ref{compute_regression} of Algorithm~\ref{alg:regressionNoiseModel} of \ref{app:methodalgorithms}), with $\card{\paramDomainTrain}=\ntrain/$ training parameter instances.  After training, all combinations took less than 0.2 seconds to predict during the online stage.}}
\begin{figure}
   \centering
   \begin{subfigure}[b]{\widthfraction/\textwidth}
      \centering
      \includegraphics[scale=\scalefactor/,clip=true,trim=\lefttrim/ 0in \righttrim/ 0in]{fig/\casetxt/\descriptiontxt_\vartxta_global\traininfo/.pdf}
      \caption{Response = \vartexa, \globaltex/}
      \label{fig:\casetxt/\descriptiontxt_\vartxta_global\traininfo/}
      \vspace{-0.4em}
   \end{subfigure}%
   \ifshowlocal%
   \hspace*{-1.5em}
   \begin{subfigure}[b]{\widthfraction/\textwidth}
      \centering
      \includegraphics[scale=\scalefactor/,clip=true,trim=\lefttrim/ 0in \righttrim/ 0in]{fig/\casetxt/\descriptiontxt_\vartxta_local\traininfo/.pdf}
      \caption{Response = \vartexa, \localtex/}
      \label{fig:\casetxt/\descriptiontxt_\vartxta_local\traininfo/}
      \vspace{-0.4em}
   \end{subfigure}%
   \fi
   \ifnum \nvar>1
   \\[1em]
   \begin{subfigure}[b]{\widthfraction/\textwidth}
      \centering
      \includegraphics[scale=\scalefactor/,clip=true,trim=\lefttrim/ 0in \righttrim/ 0in]{fig/\casetxt/\descriptiontxt_\vartxtb_global\traininfo/.pdf}
      \caption{Response = \vartexb, \globaltex/}
      \label{fig:\casetxt/\descriptiontxt_\vartxtb_global\traininfo/}
      \vspace{-0.4em}
   \end{subfigure}%
   \ifshowlocal
   \hspace*{-1.5em}
   \begin{subfigure}[b]{\widthfraction/\textwidth}
      \centering
      \includegraphics[scale=\scalefactor/,clip=true,trim=\lefttrim/ 0in \righttrim/ 0in]{fig/\casetxt/\descriptiontxt_\vartxtb_local\traininfo/.pdf}
      \caption{Response = \vartexb, \localtex/}
      \label{fig:\casetxt/\descriptiontxt_\vartxtb_local\traininfo/}
      \vspace{-0.4em}
   \end{subfigure}%
   \fi
   \fi
   \ifnum \nvar>2
   \\[1em]
   \begin{subfigure}[b]{\widthfraction/\textwidth}
      \centering
      \includegraphics[scale=\scalefactor/,clip=true,trim=\lefttrim/ 0in \righttrim/ 0in]{fig/\casetxt/\descriptiontxt_\vartxtc_global\traininfo/.pdf}
      \caption{Response = \vartexc, \globaltex/}
      \label{fig:\casetxt/\descriptiontxt_\vartxtc_global\traininfo/}
      \vspace{-0.4em}
   \end{subfigure}%
   \ifshowlocal
   \hspace*{-1.5em}
   \begin{subfigure}[b]{\widthfraction/\textwidth}
      \centering
      \includegraphics[scale=\scalefactor/,clip=true,trim=\lefttrim/ 0in \righttrim/ 0in]{fig/\casetxt/\descriptiontxt_\vartxtc_local\traininfo/.pdf}
      \caption{Response = \vartexc, \localtex/}
      \label{fig:\casetxt/\descriptiontxt_\vartxtc_local\traininfo/}
      \vspace{-0.4em}
   \end{subfigure}%
   \fi
   \fi
   \caption{\casetex: \descriptiontex/}
   \label{fig:\casetxt/\descriptiontxt\traininfo/}
\end{figure}

\fi

\def\traininfo/{}
\ifshowlocal
   \def\lefttrim/{2.4in}   
   \def\righttrim/{1.2in} 
   \def\widthfraction/{0.5}
\else
   \def\lefttrim/{1.3in}  
   \def\righttrim/{1.3in} 
   \def\widthfraction/{1}
\fi
\def \best {ind}
\def \descriptiontxt{best_\best}
\def\bestapproach/{feature choice}
\def\fromapproach/{ML technique}
\def\descriptiontex/{Comparison of FVU using best \bestapproach/ for each \fromapproach/.}

\def \best {ml}
\def \descriptiontxt{best_\best}
\def\bestapproach/{ML technique}
\def\fromapproach/{feature choice}
\def\descriptiontex/{Comparison of FVU using best \bestapproach/ for each \fromapproach/.}

\def\lefttrim/{0in}
\def\righttrim/{0in}
\def\traininfo/{_t\ntrain/}
\def \descriptiontxt{matrix}
\def\descriptiontex/{Comparison of FVU \ourRereading{between} each \bestapproach/ and \fromapproach/ with $\card{\paramDomainTrain}=\ntrain/$ training parameter instances.}

\def\lefttrim/{0in} 
\def\righttrim/{0in} 
\def \descriptiontxt{exact_v_pred}
\def\descriptiontex/{Comparison of predicted responses, relative to exact.}

\begin{figure}
   \centering
   \begin{subfigure}[b]{.99\textwidth}
      \centering
      \includegraphics[scale=\scalefactor/,clip=true,trim=\lefttrim/ 0in \righttrim/ 0in]{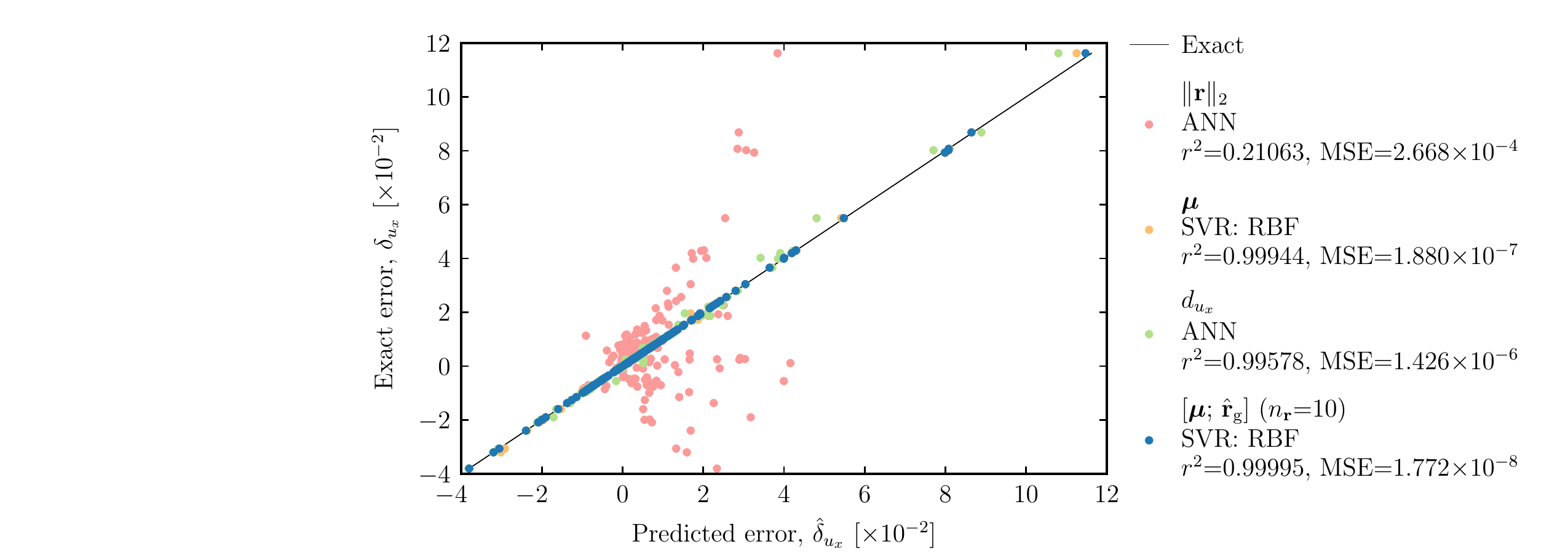}
      \caption{Response = \vartexa}
      \label{fig:\casetxt/\descriptiontxt_\vartxta}
      \vspace{-0.4em}
   \end{subfigure}%
   \ifnum \nvar>1
   \\[1em]
   \begin{subfigure}[b]{.99\textwidth}
      \centering
      \includegraphics[scale=\scalefactor/,clip=true,trim=\lefttrim/ 0in \righttrim/ 0in]{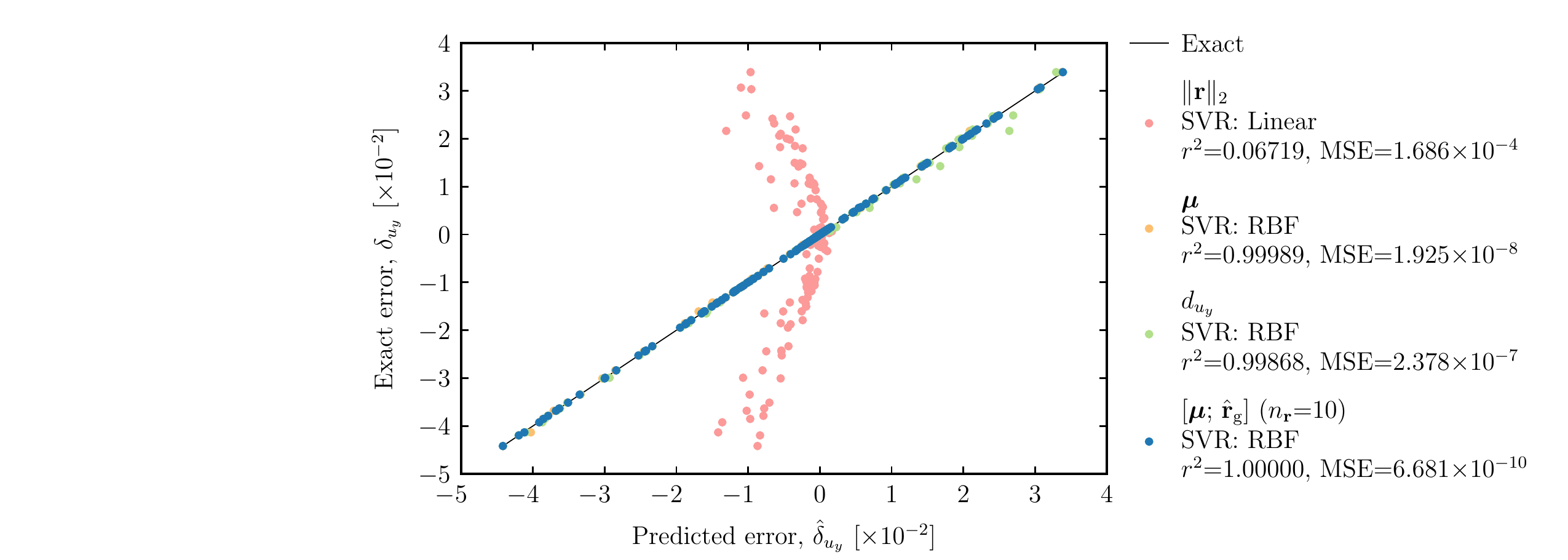}
      \caption{Response = \vartexb}
      \label{fig:\casetxt/\descriptiontxt_\vartxtb}
      \vspace{-0.4em}
   \end{subfigure}%
   \fi
   \ifnum \nvar>2
   \\[1em]
   \begin{subfigure}[b]{.99\textwidth}
      \centering
      \includegraphics[scale=\scalefactor/,clip=true,trim=\lefttrim/ 0in \righttrim/ 0in]{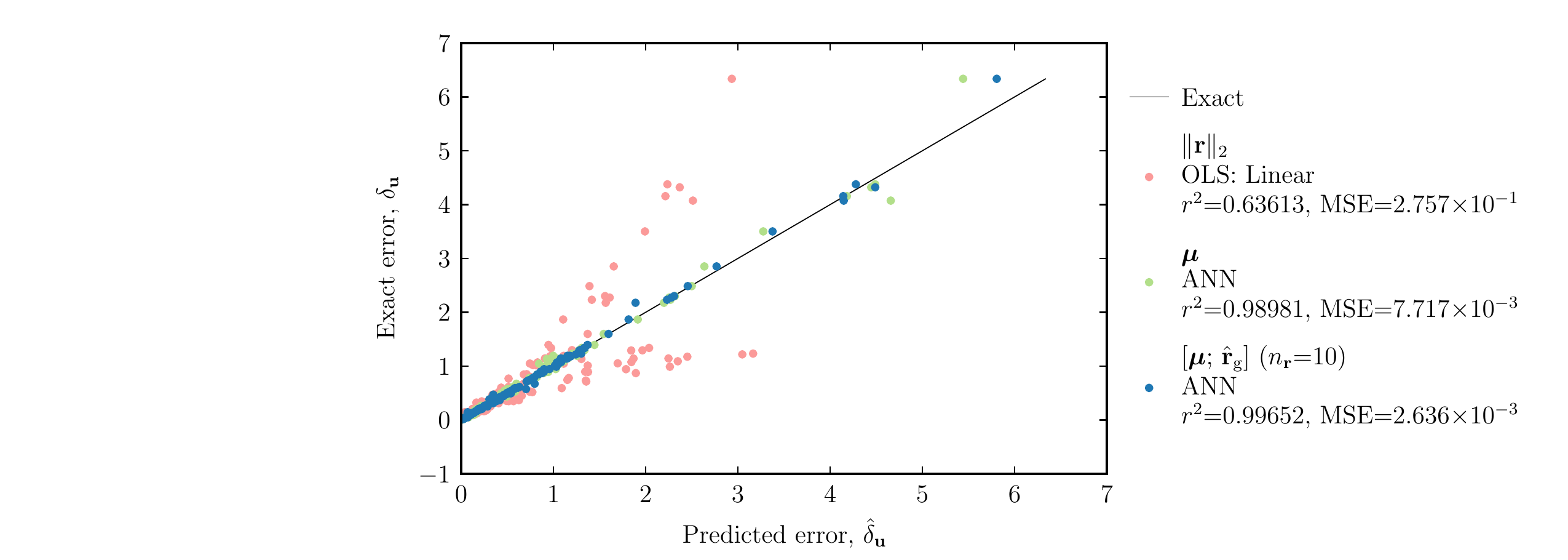}
      \caption{Response = \vartexc}
      \label{fig:\casetxt/\descriptiontxt_\vartxtc}
      \vspace{-0.4em}
   \end{subfigure}%
   \fi
   \caption{\casetex: \descriptiontex/}
   \label{fig:\casetxt/\descriptiontxt}
\end{figure}

\ifhasgauss
\def\lefttrim/{0in} 
\def\righttrim/{0in} 

\def \descriptiontxt{hist_gauss}
\def\descriptiontex/{Error in predicted responses for second set of test data, scaled by $\varianceEstimate =\text{MSE}$ from the first set of test data.}

\begin{figure}
   \centering
   \begin{subfigure}[b]{.99\textwidth}
      \centering
      \includegraphics[scale=\scalefactor/,clip=true,trim=\lefttrim/ 0in \righttrim/ 0in]{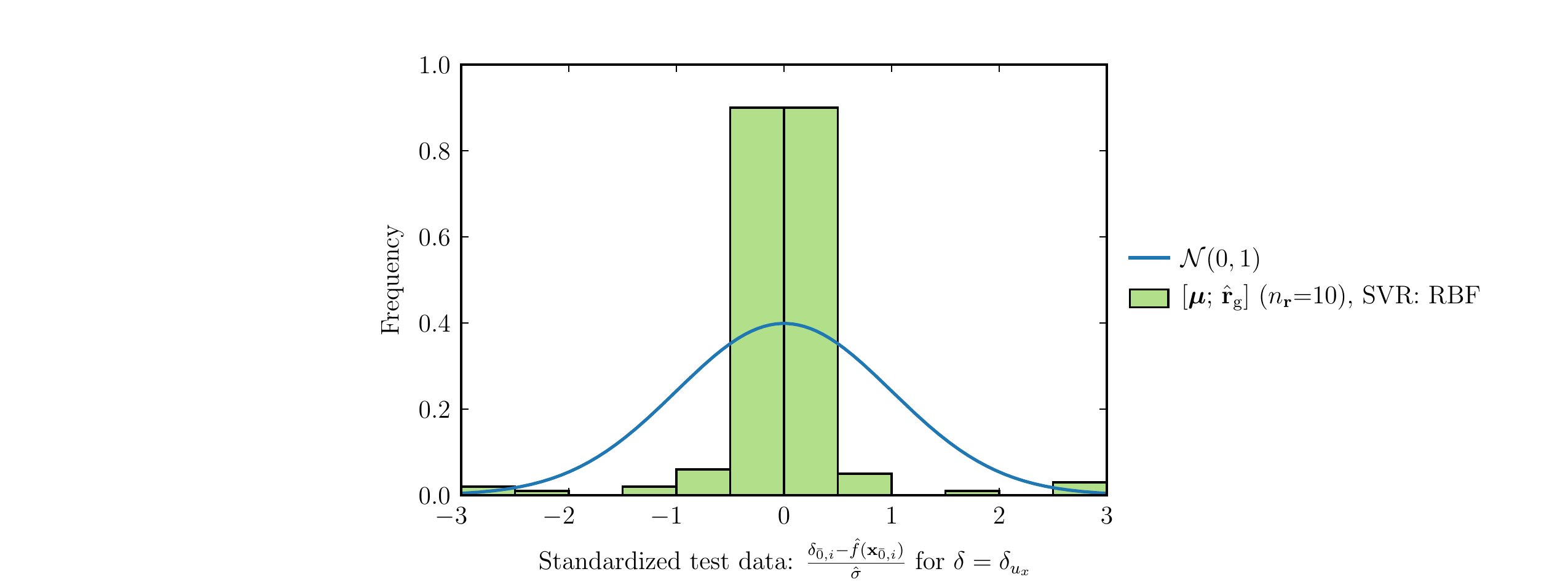}
      \caption{Response = \vartexa}
      \label{fig:\casetxt/\descriptiontxt_\vartxta}
      \vspace{-0.4em}
   \end{subfigure}%
   \ifnum \nvar>1
   \\[1em]
   \begin{subfigure}[b]{.99\textwidth}
      \centering
      \includegraphics[scale=\scalefactor/,clip=true,trim=\lefttrim/ 0in \righttrim/ 0in]{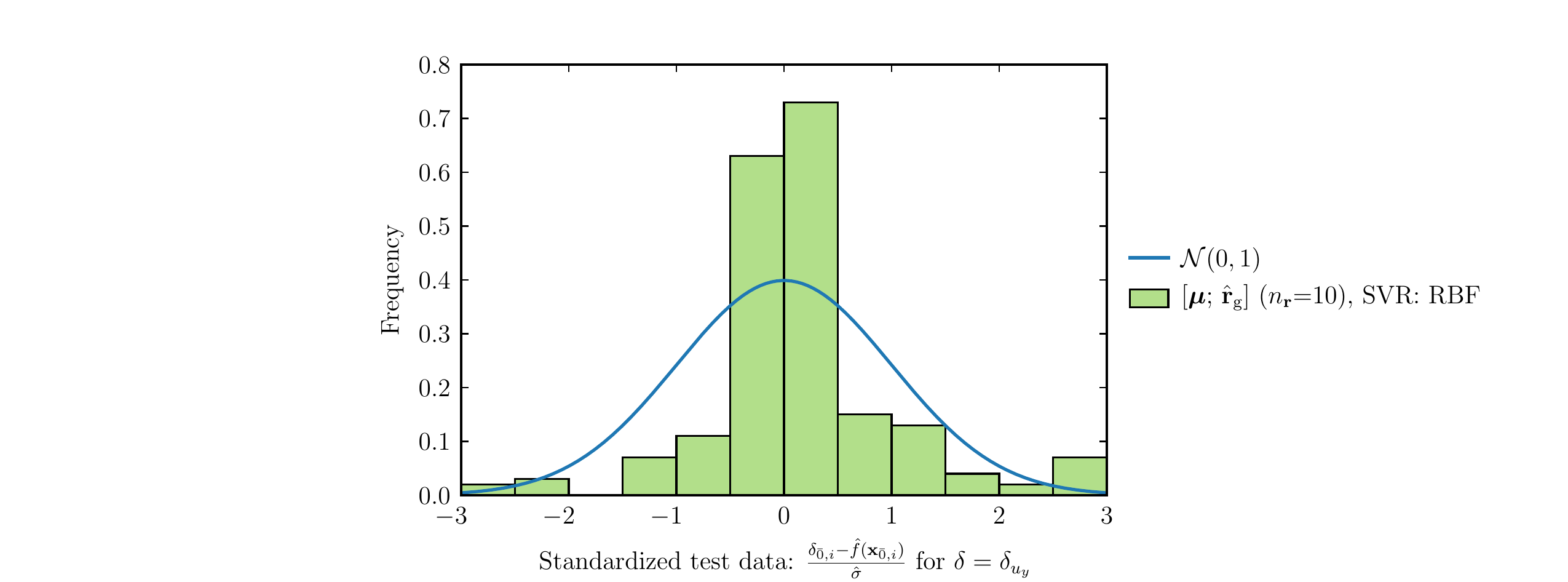}
      \caption{Response = \vartexb}
      \label{fig:\casetxt/\descriptiontxt_\vartxtb}
      \vspace{-0.4em}
   \end{subfigure}%
   \fi
   \ifnum \nvar>2
   \\[1em]
   \begin{subfigure}[b]{.99\textwidth}
      \centering
      \includegraphics[scale=\scalefactor/,clip=true,trim=\lefttrim/ 0in \righttrim/ 0in]{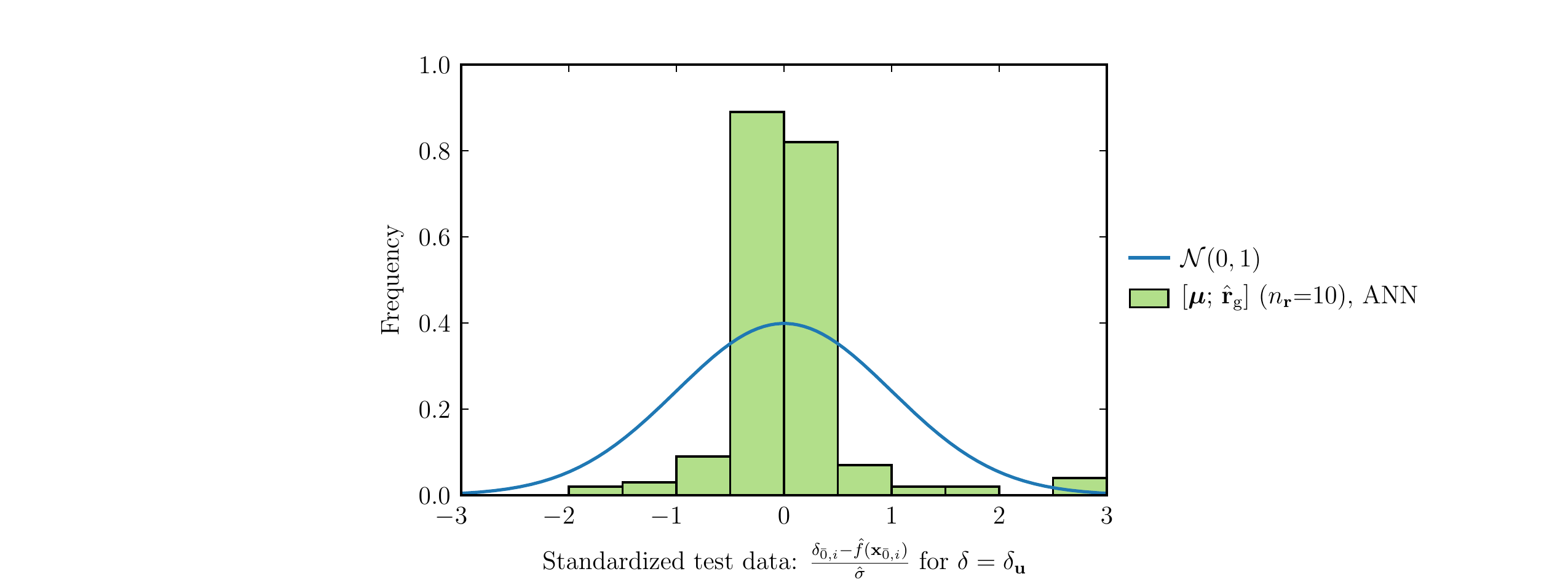}
      \caption{Response = \vartexc}
      \label{fig:\casetxt/\descriptiontxt_\vartxtc}
      \vspace{-0.4em}
   \end{subfigure}%
   \fi
   \caption{\casetex: \descriptiontex/}
   \label{fig:\casetxt/\descriptiontxt}
\end{figure}
\fi

\clearpage
\subsection{PCAP: reduced-order modeling} 
\label{subsec:pcap}
The second experiment considers the mechanically induced deformation of the
Predictive Capability Assessment Project (PCAP) test case with an approximate
solution provided by a Galerkin reduced-order model as described in Section
\ref{section:ROM}.
We also conduct this experiment using Albany~\cite{albany_2016} \reviewerTwo{with eight-node hexahedral elements}.\reviewerOne{
We compute solutions with a single nonlinear solve using a damped Newton's method.
}

\FloatBarrier 

\subsubsection{Overview} 
\label{subsec:pcap_overview}
\begin{figure}
\centering
\includegraphics[scale=.3]{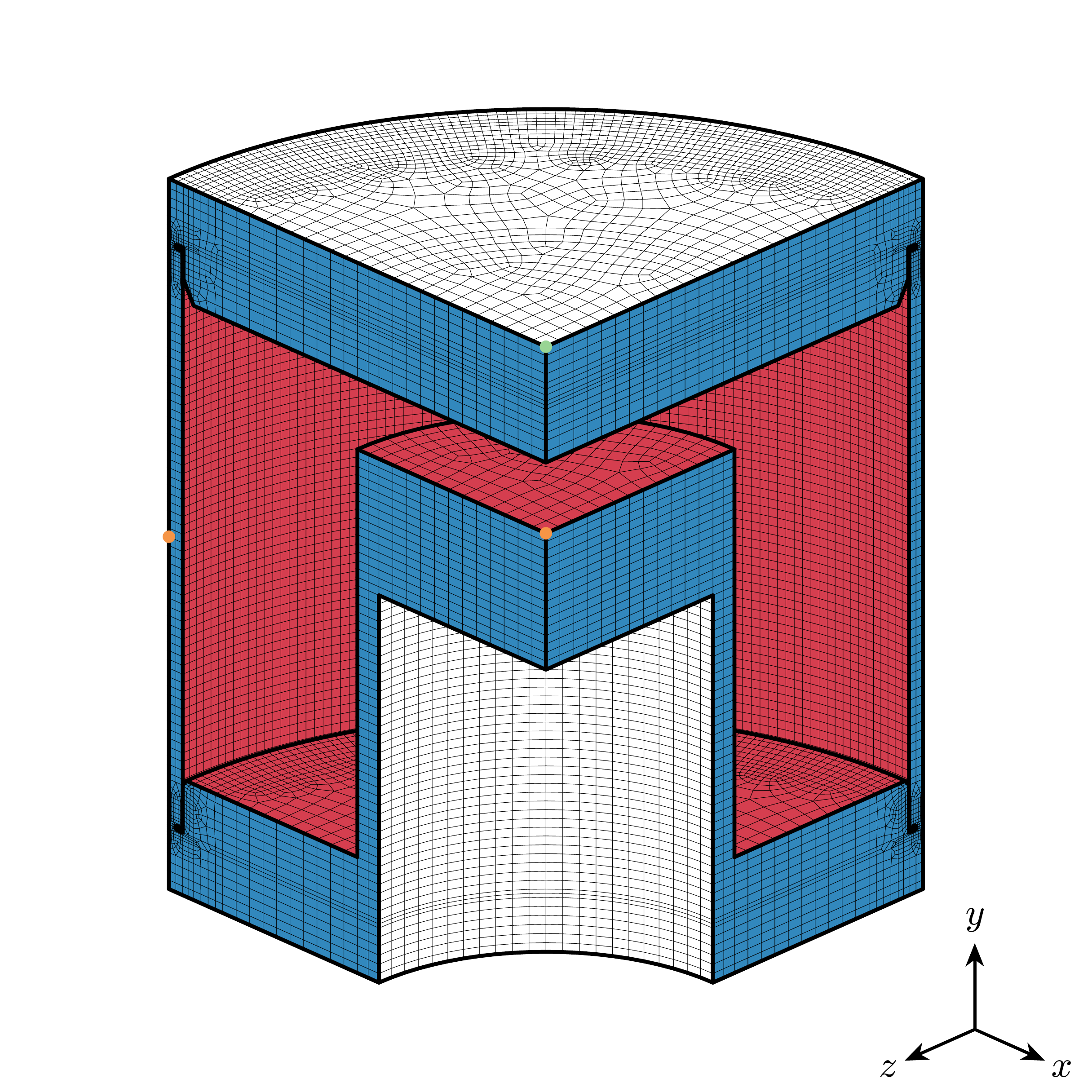}
\vspace{-0.4em}
\caption{PCAP: Mesh with boundary conditions and nodes of interest: red
denotes applied pressure (Neumann boundary condition), blue denotes
planar-displacement constraint
(Dirichlet boundary condition), green denotes zero-displacement constraint (Dirichlet
boundary condition), and orange denotes nodes of interest.}
\label{fig:pcap/mesh}
\end{figure}

The PCAP test case is an axisymmetric structure consisting of a tube with
different lids attached at each end.  The radius of the test case is 4.44 cm,
and the height is 6.61 cm.  As illustrated in Figure~\ref{fig:pcap/mesh}, a
quarter of the test case is discretized with 77,768 elements and 92,767 nodes.

Pressure is applied as a Neumann boundary condition on the interior of the
structure, and the center of the top of the upper lid is further constrained
in the $y$-direction.  Due to the axisymmetry of the problem, we introduce
Dirichlet boundary conditions on the sides to enforce planar displacement.  Therefore, the size of the problem is $\nstate=274,954$.

We consider two quantities of interest: (1) the $y$-displacement of
the center of the top of the bottom lid, denoted by $u_y$, and (2) the radial deformation
of the outside of the tube, midway along the height, denoted by $u_r$.  
Figure~\ref{fig:pcap/mesh} depicts
the boundary
conditions and nodes of interest.

\begin{figure}
\centering
\includegraphics[scale=.3]{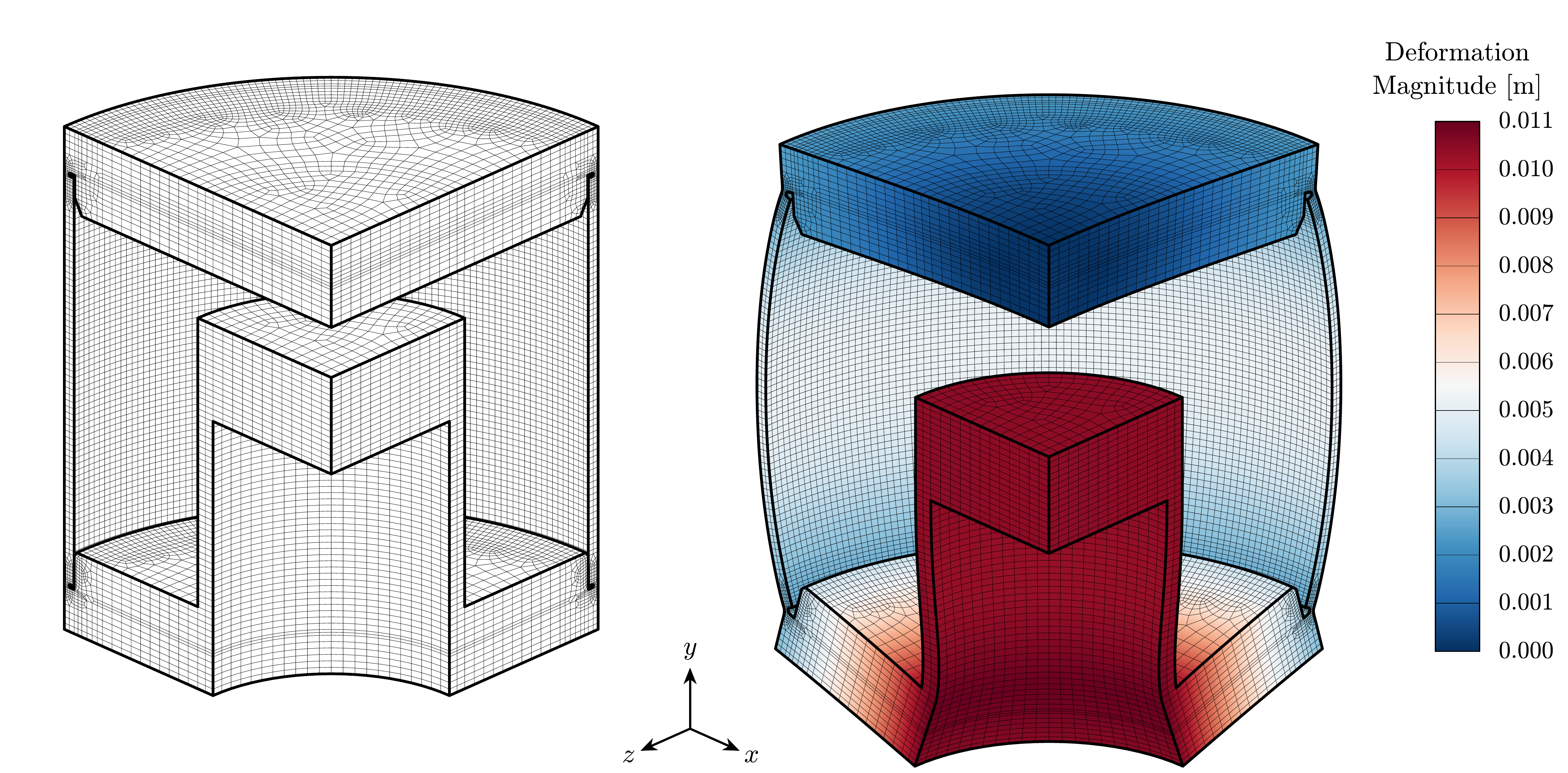}
\vspace{-0.4em}
\caption{PCAP: Largest simulated deformation (right) compared to undeformed state (left).}
\label{fig:pcap/large_def}
\end{figure}

The lid material has an elastic modulus of 110.3 GPa, a Poisson ratio of 0.32,
a yield strength of 88.02 MPa, and a hardening modulus of 232.8 MPa.  The tube
has a yield strength of 82.94 MPa, and a hardening modulus of 1.088 GPa.  
This experiment considers $\nparams=3$ parameters. 
The tube elastic modulus $\paramArg{1}=E$ is varied between 50.0 and 100.0 GPa,
the tube Poisson ratio $\paramArg{2}=\nu$ is varied between 0.200 and
0.350, and the applied pressure $\paramArg{3}=p$ is varied between 6.00 and 10.00 GPa.
The resulting parameter domain is thus $\paramDomain=[
	50.0\,\text{GPa}, 100.0\, \text{GPa}
]\times[0.200,
0.350]\times [
	6.00\,\text{GPa}, 10.00\,\text{GPa}
]$.
Across the parameter domain $\paramDomain$, $u_y$ varies between limits of $-0.0876$ cm and $-1.04$ cm,
which corresponds to 1.3\%--15.8\% of the undeformed height, and $u_r$ varies
between limits of 0.185
cm and 0.418 cm, which corresponds to 4.2\%--9.4\% of the undeformed radius.
Figure~\ref{fig:pcap/large_def} depicts the test case when undergoing the
greatest deformation.

\begin{figure}
\vspace{-4em}
\centering
\begin{subfigure}[b]{.49\textwidth}
\includegraphics[scale=.3]{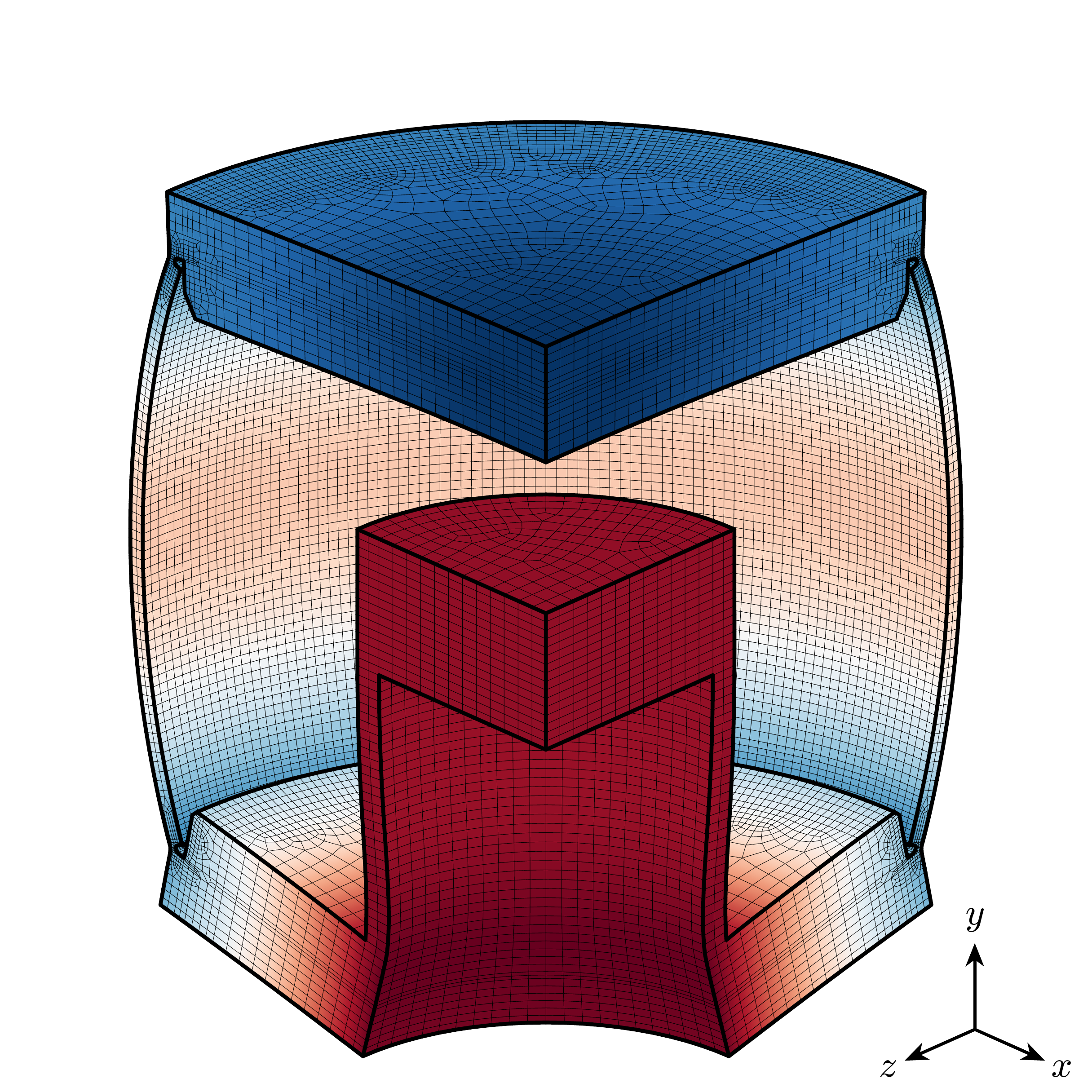}
\caption{First \reviewerTwo{POD basis vector}}
\label{fig:pcap/mesh_plots/basis1}
\vspace{-0.4em}
\end{subfigure}
\begin{subfigure}[b]{.49\textwidth}
\includegraphics[scale=.3]{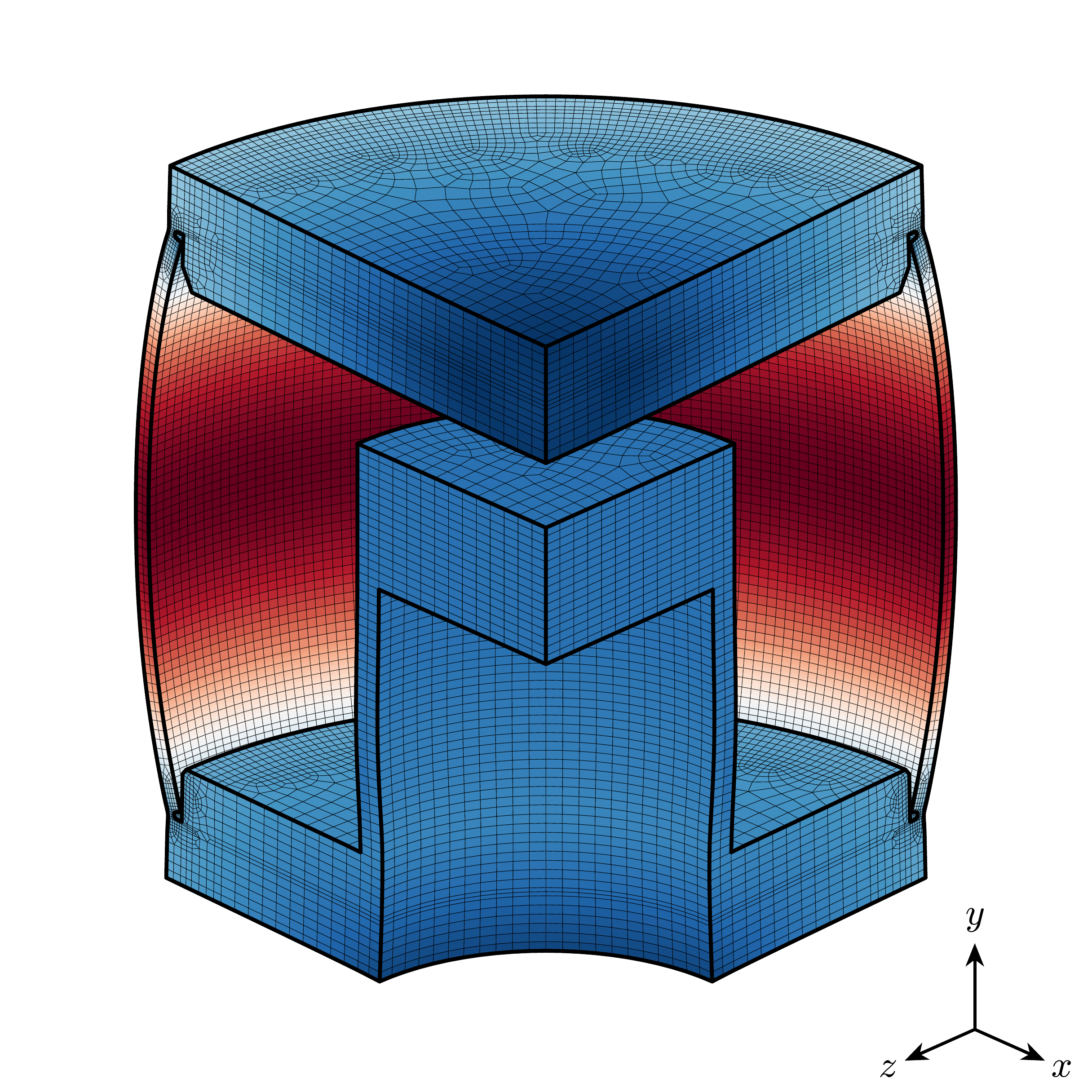}
\caption{Second \reviewerTwo{POD basis vector}}
\label{fig:pcap/mesh_plots/basis2}
\vspace{-0.4em}
\end{subfigure}
\\
\begin{subfigure}[b]{.49\textwidth}
\includegraphics[scale=.3]{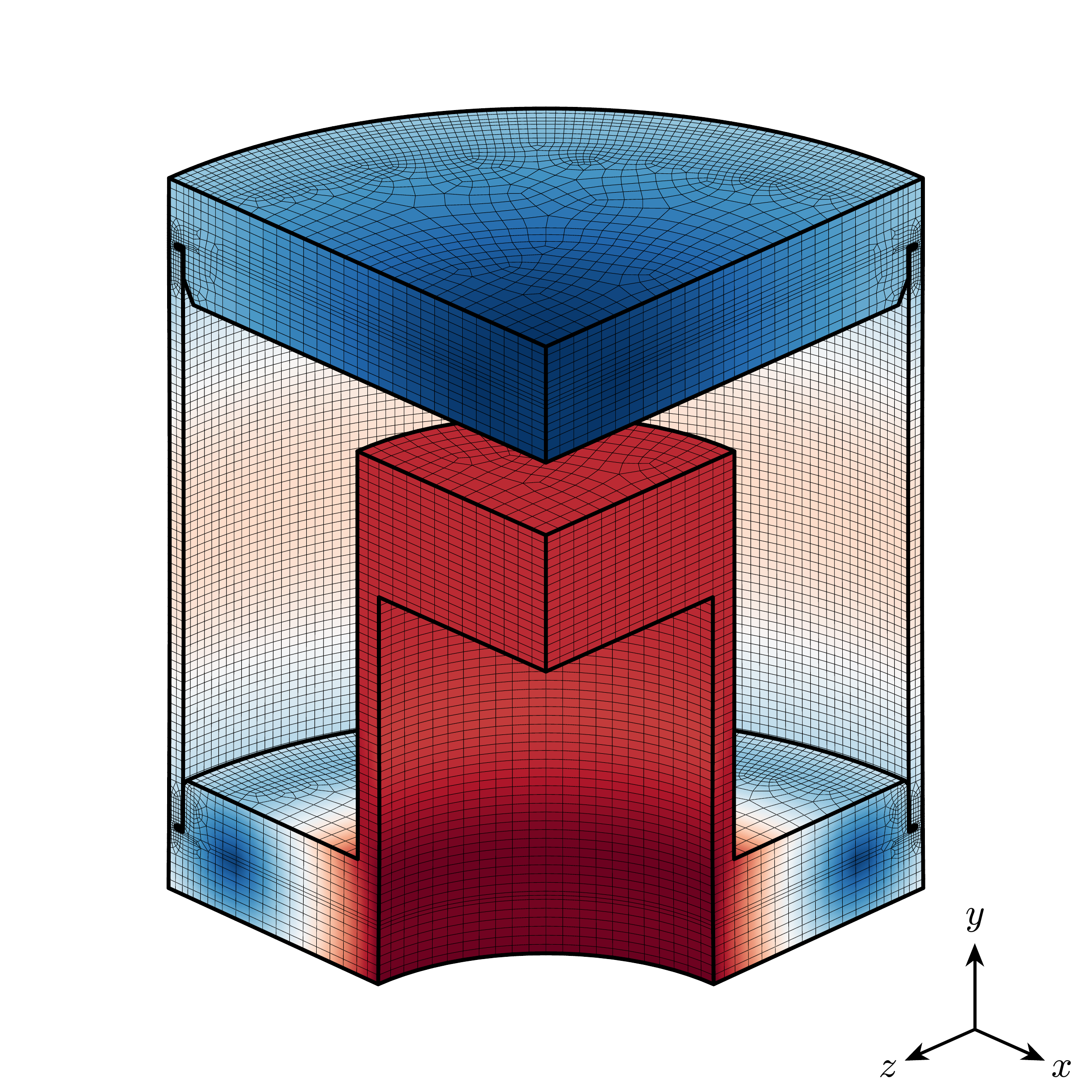}
\caption{Third \reviewerTwo{POD basis vector}}
\label{fig:pcap/mesh_plots/basis3}
\vspace{-0.4em}
\end{subfigure}
\begin{subfigure}[b]{.49\textwidth}
\includegraphics[scale=.3]{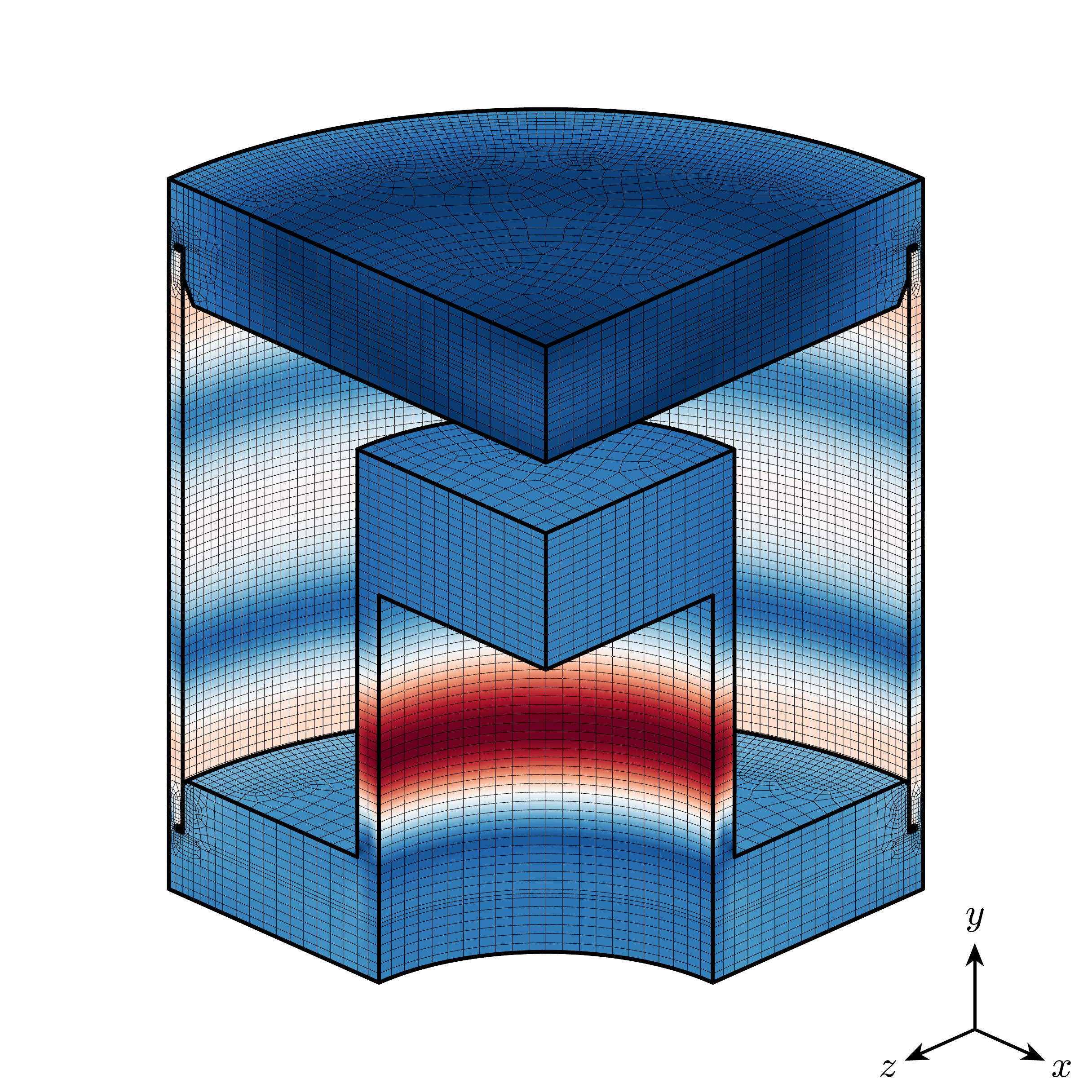}
\caption{Fourth \reviewerTwo{POD basis vector}}
\label{fig:pcap/mesh_plots/basis4}
\vspace{-0.4em}
\end{subfigure}
\\
\begin{subfigure}[b]{.49\textwidth}
\includegraphics[scale=.3]{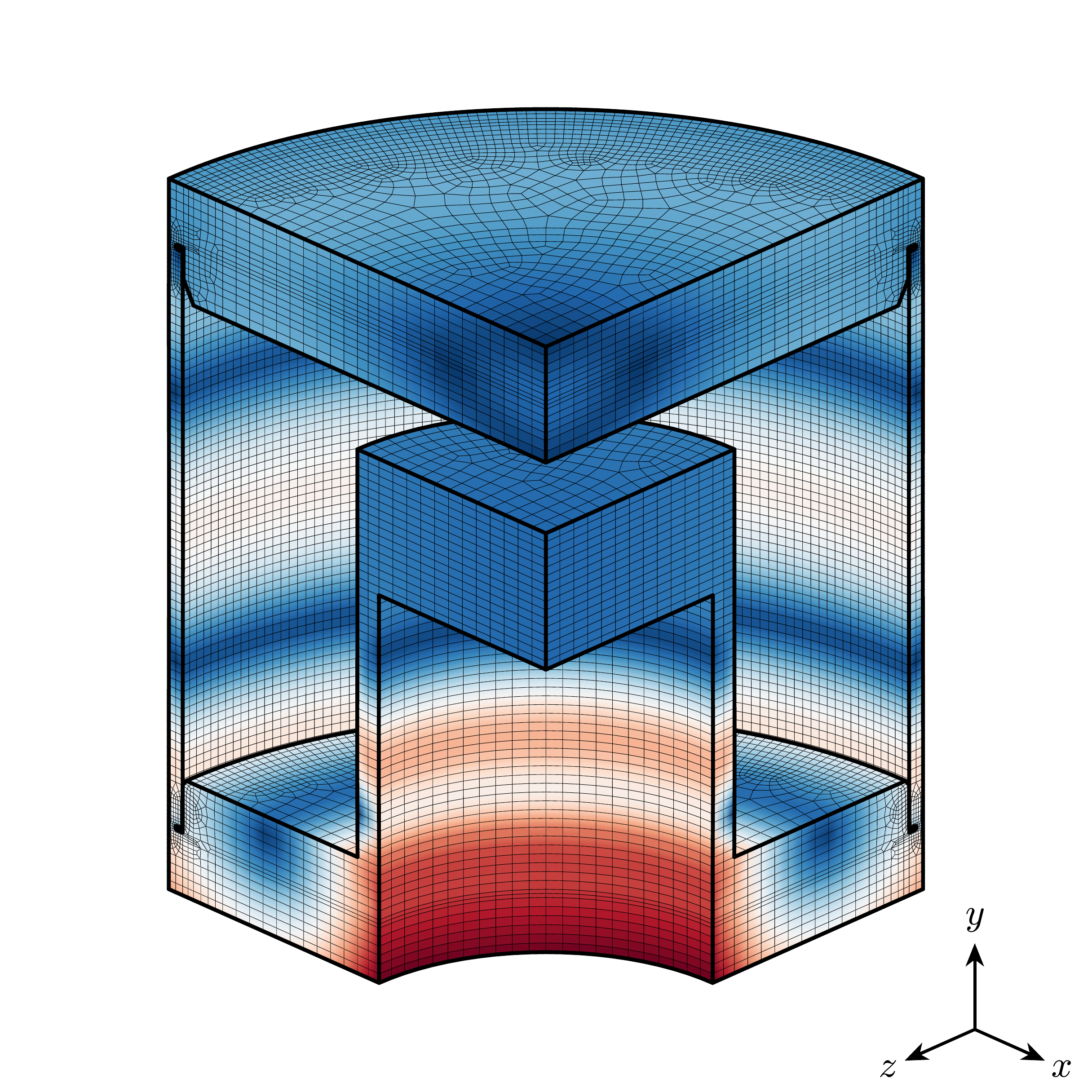}
\caption{Fifth \reviewerTwo{POD basis vector}}
\label{fig:pcap/mesh_plots/basis5}
\vspace{-0.4em}
\end{subfigure}
\begin{subfigure}[b]{.49\textwidth}
\caption*{}
\vspace{-0.4em}
\end{subfigure}
\caption{PCAP: \reviewerTwo{POD basis vectors}.  Red indicates largest deformation; blue indicates smallest deformation.}
\label{fig:pcap/basis_functions}
\end{figure}

To construct the snapshots, we first compute the high-fidelity-model
solution $\state(\params)$ for $\params\in\paramDomainPOD\subset \paramDomain$,
where $\paramDomainPOD$ comprises eight Latin hypercube samples 
such that $\nsnapshots=8$. We subsequently compute 
the trial-basis matrix $\trialbasis$  using proper orthogonal
decomposition (POD) by applying Algorithm~\ref{alg:pod} of \ref{app:pod} with inputs
$\{\state(\params)\}_{\params\in\paramDomainPOD}$.
The relative statistical energy captured by different basis
dimensions (see output of Algorithm \ref{alg:pod}) corresponds to $\energyThresholdArg{1}= 0.8503$,
$\energyThresholdArg{2}=0.9569$, 
$\energyThresholdArg{3}=0.9935$, 
$\energyThresholdArg{4}=0.9977$, 
$\energyThresholdArg{5}=0.9990$, 
$\energyThresholdArg{6}=0.9997$, 
$\energyThresholdArg{7}=0.9999$, and
$\energyThresholdArg{8}=1.0000$.
	Figure~\ref{fig:pcap/basis_functions} illustrates the first five POD \reviewerTwo{basis} vectors.

One set of 30 parameter instances is randomly sampled from the parameter domain
$\paramDomain$ to serve as the set of \pti{} $\paramDomainTrain$.  To assess
method performance for smaller amounts of training data,
we randomly create nested training sets from this set. Another set of 30
parameter instances is
randomly sampled to serve as the set of
parameter testing instances $\paramDomainTest$.
We consider approximate solutions generated by a Galerkin reduced-order model
(i.e., $\testbasis = \trialbasis$)
with $\stateRef=\zero$ using 
$\napprox=5$ different
basis dimensions: $\nbasisArg{i}=i$, $i=1,\ldots,5$.
We consider both \ref{data:pooled} (pooled) and \ref{data:unique} (unique) as described in Section \ref{subsection:trainingdata}.
The former case yields $\ntestData=150$ testing data points and up to
$\ntrainingData=150$ training points, while the latter yields
$\ntestData=30$ testing points and up to
$\ntrainingData=30$ testing points.
We executed the simulations on a workstation with 
\ourRereading{two 2.60-GHz Intel\textsuperscript{\tiny\textregistered}
Xeon\textsuperscript{\tiny\textregistered} E5-2660 v3 processors, which collectively contain 20 logical cores,} and 62.8
GB of RAM.
Table~\ref{tab:pcap/run_times} lists the average and median run times across
all training and testing FOM and ROM simulations.  \reviewerOne{For several of the parameter instances, for $\nbasisArg{2}$ and $\nbasisArg{3}$, the nonlinear solver required damping to reduce the step size during the Newton solve in order to converge.  For this reason, the average and median times deviate from the monotonically increasing trend among $\nbasisArg{1}$, $\nbasisArg{4}$, and $\nbasisArg{5}$}.  We note that the ROM yields
computational-cost savings in this case---despite the lack of
hyper-reduction---due to the low dimensionality of the POD basis.

\begin{table}
\centering
\begin{tabular}{r c r r r r r}
\toprule
&  & \multicolumn{5}{c}{Galerkin ROM} \\
\cmidrule(l){3-7}
& FOM & \multicolumn{1}{c}{$\nbasisArg{1}=1$} &
\multicolumn{1}{c}{$\nbasisArg{2}=2$} &
\multicolumn{1}{c}{$\nbasisArg{3}=3$} & \multicolumn{1}{c}{$\nbasisArg{4}=4$} &
\multicolumn{1}{c}{$\nbasisArg{5}=5$} \\
\midrule
Average time [seconds] & 2208 & 36.5 & 182.0 & 151.6 & 81.8 & 93.7 \\
 Median time [seconds] & 1840 & 35.5 &  82.4 & 108.0 & 82.4 & 87.3 \\
\bottomrule
\end{tabular}
\vspace{-.5em}
\caption{PCAP: FOM and ROM average and median run times.}
\label{tab:pcap/run_times}
\end{table}
The large dimension of this problem renders two feature-engineering methods
impractical: Feature-Engineering Method \ref{feat:DWR}, as it requires computing the dual-weighted
residual; and Feature-Engineering Method \ref{feat:paramres}, as it employs the entire
$\nstate$-dimensional residual vector as features. Thus, we do not consider
these methods in this set of experiments.  When using
Feature-Engineering Methods~\ref{feat:paramgappyPODres} and \ref{feat:paramsampleres}, three
sample sizes are used: $\nsamples=\{10,\, 100,\,1000\}$.

\subsubsection{Results} 
\label{subsubsec:pcap_results}
\showlocaltrue
\hasgausstrue
\ispcaptrue
\def \casetxt {pcap}
\def \casetex {PCAP}

\def \vartxta {er}
\def \vartexa {$\errorQoiArg{u_r}$}
\def \vartxtb {ey}
\def \vartexb {$\errorQoiArg{u_y}$}
\def \vartxtc {epsilon}
\def \vartexc {$\errorState$}
\def \nvar    {3}
\def\ntrain/{30}
\def\kpercentbetter/{13.49\%}
\def\qpercentbetter/{86.51\%}
\def\extratext/{}
\def\localtex/{\ref{data:unique}}
\def\globaltex/{\ref{data:pooled}}

Similarly to Section \ref{cube_results}, 
we first assess the difference in performance between the $k$- and
$q$-sampling approaches for computing the sampling matrix $\sampleMat$
employed by Feature-Engineering Methods \ref{feat:paramgappyPODres} and \ref{feat:paramsampleres}.
Figure~\ref{fig:\casetxt/kq_histogram} provides a comparison of the test MSEs
that arise from using $k$- and $q$-sampling with
$\card{\paramDomainTrain}=\ntrain/$ parameter training instances.
This
	figure is generated from 252 total data points, which aggregate test MSE
	values over three errors: $\errorQoiArg{u_r}$, $\errorQoiArg{u_y}$, and
	$\errorState$; Feature-Engineering Methods \ref{feat:paramgappyPODres} and
	\ref{feat:paramsampleres}; the seven regression techniques discussed
	in~\ref{subsection:regressionfunctionapprox}; three numbers of sample points
	$\nsamples\in\{10, 100, 1000\}$; and the two data-set approaches discussed
in Section \ref{subsection:trainingdata}.
In this case, $q$-sampling significantly outperforms $k$-sampling;  therefore, the remaining results for this
experiment only consider $q$-sampling.

\reviewerThree{
We trained all of the models on a workstation with
a 3.60-GHz Intel\textsuperscript{\tiny\textregistered} Core\textsuperscript{\tiny\texttrademark} i7-4790 processor, which contains eight logical cores, and 16.0 GB of RAM.
Figure~\ref{fig:\casetxt/time_matrix_t\ntrain/} shows $t_\text{train}$, the amount of wall time in seconds it took to train the regression model offline, for each
combination of feature-engineering method, regression technique, and data-set
approach for $\card{\paramDomainTrain}=\ntrain/$ \pti{}.  This time includes Steps \ref{cv_loop}--\ref{compute_regression} of Algorithm~\ref{alg:regressionNoiseModel} of \ref{app:methodalgorithms}, which account for cross-validation.
As expected, the simplest machine-learning technique, OLS: Linear, required the least amount of time, sometimes less than 16 milliseconds, whereas the most complex, ANN, required the most, as high as 1.64 hours.  Additionally, $\FeatureParResSampled$ with $\nsamples=1000$ required the most amount of time with ANN, due to having the greatest number of features.
After training, all combinations took less than 0.2 seconds to predict during the online stage.
}

For each of the ML regression techniques, Figure~\ref{fig:\casetxt/best_ind} reports how the 
test FVU varies with respect to the number of 
\pti{}
$\card{\paramDomainTrain}$
when using the best-performing feature-engineering method for each regression technique.  
Generally, as the number of \pti{} 
$\card{\paramDomainTrain}$
increases, the test FVU decreases,
and the best regression techniques are those that enable higher capacity:
ANN and SVR: RBF.  In
Figures~\ref{fig:\casetxt/best_ind_er_global},
\ref{fig:\casetxt/best_ind_er_local}, and~\ref{fig:\casetxt/best_ind_ey_local},
RF and $k$-NN perform the least competitively, most likely due to an insufficient
amount of training data. 
\ref{data:pooled} and \ref{data:unique} again perform similarly in this case.
Compared to the cube experiments discussed in Section \ref{cube_results}, the
test FVU values are noticeably larger in the present experiments, likely due
to the greater complexity (i.e., higher solution dimensionality, more complex
geometry) of the experiment.  Once more, predicting $\errorState$ is more
challenging than predicting $\errorQoiArg{u_r}$ and $\errorQoiArg{u_y}$; in
this case, ANN significantly outperforms other regression techniques, as shown
in Figures~\ref{fig:\casetxt/best_ind_epsilon_global}
and~\ref{fig:\casetxt/best_ind_epsilon_local}.


Whereas Figure~\ref{fig:\casetxt/best_ind} compares each regression technique when
using the best feature-engineering method for that technique,
Figure~\ref{fig:\casetxt/best_ml} compares each feature-engineering method when using the
best regression technique for that method.
Once again, it is immediately clear that 
using the feature $\FeatureResNorm$ alone yields
the highest test FVU, which does not improve as the amount of training data increases.
Additionally, features $\FeatureParResNorm$ and $\FeaturePar$ perform
quite poorly, with test FVU values orders of magnitude higher than
features
$\FeatureParResSampled$ and $\FeatureParResGappyPCA$. This is a notable
result, as employing features $\FeaturePar$ is the approach most commonly
adopted in the literature (i.e., via the `multifidelity correction' or `model
discrepancy' methods).


Figure~\ref{fig:\casetxt/matrix_t\ntrain/} reports the test FVU for each
combination of feature-engineering method, regression technique, and data-set
approach for $\card{\paramDomainTrain}=\ntrain/$ \pti{}.
First, we note that SVR: RBF and ANN consistently yield the best performance,
whereas OLS: Linear and OLS: Quadratic yield inconsistent performance.
In contrast to the cube experiments,
Figures~\ref{fig:\casetxt/matrix_er_global_t\ntrain/} and
\ref{fig:\casetxt/matrix_er_local_t\ntrain/} demonstrate that there is a
benefit to training using \ref{data:unique} rather than \ref{data:pooled}.
Generally, features $\FeatureParResGappyPCA$ and $\FeatureParResSampled$ yield
similar performance, with 
$\FeatureParResGappyPCA$ slightly outperforming $\FeatureParResSampled$ in
some cases. 
Because features $\FeatureResNorm$, $\FeatureParResNorm$, and $\FeaturePar$
correspond to a small number of features, their performance is nearly
insensitive to regression techniques, as all regression techniques yield
similarly low capacities for a small number of features.  None of these features
yield good performance in this case.
We remark that, in this experiment, only $\nsamples=10\ll
\nstate=278,301$ elements of the residual must be computed to realize
coefficients of determination above $r^2=0.99$ in all cases with ANN
regression; this again implies that excellent performance can be obtained
using only a very small number of cheaply computable features with the
proposed methodology.  


For $\card{\paramDomainTrain}=\ntrain/$ \pti{},
Figure~\ref{fig:\casetxt/exact_v_pred} compares the predicted values of
$\errorQoiArg{u_r}$, $\errorQoiArg{u_y}$, and $\errorState$ with the exact
values using the conventional feature choices
(i.e., $\FeatureResNorm$, $\FeaturePar$), as well as $\FeatureParResGappyPCA$, where performance of the best
regression technique for each of these feature-engineering methods is reported.
%
In each case, 
$\FeatureParResGappyPCA$ performs better than all
conventional approaches, with $r^2>0.998$ in every case, and with the test MSE
reduced by one to two orders of magnitude.

\reviewerThree{
Figure~\ref{fig:\casetxt/hist_gauss} assesses the accuracy of the noise model
$\noiseModel\sim\normal{0}{\varianceEstimate}$ computed in
\ref{step:noiseApproximation} in Section \ref{subsec:overview} for $\FeatureParResGappyPCA$
with $\nsamples=10$ using the best-performing regression technique.
As before, 
the figure reports the standard normal distribution compared to a histogram of
the prediction errors, which have been standardized according to the
hypothesized distribution.
These prediction errors are computed 
on a second set of independent test data
$\testDataTwo\defeq\{(\errorTestTwo{i},\featuresTestTwo{i})\}_{i=1}^{\card{\testData}}$,
which is constructed in an identical way to the the first set of test data
$\testData$.
Table~\ref{tab:\casetxt/validation_frequencies} lists the corresponding
validation frequencies $\frequency_\error(\frequency)$~\eqref{eq:validation_frequency} for multiple
$\frequency$-prediction intervals~\eqref{eq:prediction_interval}.

\begin{table}[ht!]
\centering
\reviewerThree{
\begin{tabular}{c c c c}
\toprule
& \multicolumn{3}{c}{$\frequency_\error(\frequency)$}\\
\cmidrule(l){2-4}
$\frequency$ & $\errorQoiArg{u_r}$ & $\errorQoiArg{u_y}$ & $\errorState$ \\
\midrule
0.80 & 0.8800 & 0.8133 & 0.8200 \\
0.90 & 0.9333 & 0.9000 & 0.9133 \\
0.95 & 0.9400 & 0.9200 & 0.9400 \\
0.99 & 0.9533 & 0.9667 & 0.9800 \\
\bottomrule
\end{tabular}
}
\vspace{-.5em}
\caption{\casetex: Validation frequencies for models using $\FeatureParResGappyPCA$
($\nsamples=10$) with SVR: RBF, ANN, and ANN, respectively for $\errorQoiArg{u_r}$, $\errorQoiArg{u_y}$, and $\errorState$.}
\label{tab:\casetxt/validation_frequencies}
\end{table}

Again, the data are not quite Gaussian, although they appear to match the
hypothesized distribution closer than in the cube experiments of Section
\ref{subsec:cube}. Once more, the tails of the distribution are heavier than what
would be predicted by a Gaussian distribution. The table indicates that, despite the data
being non-Gaussian, some of the prediction intervals are indeed accurate.
Namely, the 0.95-prediction
interval; the 0.90-prediction interval; and the 0.90-, 0.95-, and
0.99-prediction intervals are
reasonably accurate for the models associated with
responses $\errorQoiArg{u_r}$, $\errorQoiArg{u_y}$, and $\errorState$,
respectively.
}


\clearpage
\subsection{Burgers' equation: \ourRereading{inexact solutions} and coarse solution prolongation}
\label{subsec:burgers}
\FloatBarrier 
The previous two experiments demonstrate the proposed methodology in the
context of approximate solutions provided by a Galerkin reduced-order model.
Instead, this experiment demonstrates the methodology for two other types of
approximate solutions: inexact solutions as described in Section
\ref{subsec:inexact}, and lower-fidelity models arising from a coarse mesh as
described in Section \ref{subsec:lowerfidelity}.  We assess these approximate
solutions on a forced steady viscous Burgers' equation.

\subsubsection{Overview} 
We consider the governing equation
\begin{align}
uu_x - \frac{1}{R}u_{xx} = \alpha\sin 2\pi x,
\label{eq:burgers}
\end{align}
where $x\in[0,\,1]$, $u(0)=u_a$ and $u(1)=-u_a$. We adopt the
finite-difference discretization described in Ref.~\cite{tian_2007} 
with $2001$ uniformly spaced nodes, such that $\nstate=1999$.

We consider one quantity of interest: the slope
$m$ at $x=\frac{1}{2}$ approximated by finite differences as
\begin{align*}
m \defeq \frac{-u(\frac{1}{2}+2\Delta x)+8u(\frac{1}{2}+\Delta
	x)-8u(\frac{1}{2}-\Delta x)+u(\frac{1}{2}-2\Delta x)}{12\Delta x},
\end{align*}
where $\Delta x \equiv 1/(N_\mathbf{u}+1)$.

\begin{figure}
   \centering
      \centering
      \includegraphics[scale=\scalefactor/,clip=true,trim=0in 0in 0in 0in]{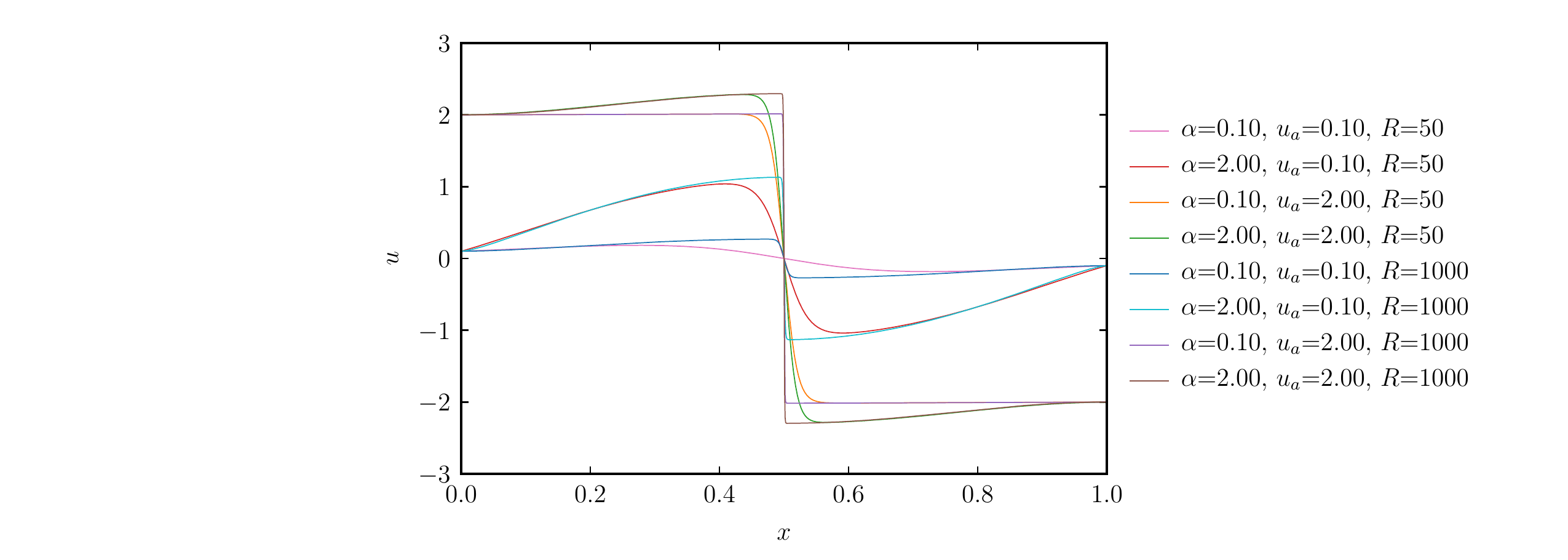}
      \vspace{-0.4em}
   \caption{Solutions to Burgers' equation~\eqref{eq:burgers} corresponding to the extrema of $\boldsymbol{\mu}$.}
   \label{fig:burgers}
\end{figure}

This experiment considers $\nparams=3$ parameters. The source parameter
$\paramArg{1} = \alpha$ varies between 0.10 and 2.00; the boundary-condition
magnitude $\paramArg{2} = u_a$ varies between 0.10 and 2.10; and the Reynolds
number $\paramArg{3} = R$ parameter varies between 50 and 1000.
Thus, the parameter domain is $\paramDomain=[0.10,2.00]\times[0.10,2.10]\times
[50,1000]$.
Figure~\ref{fig:burgers} depicts the solutions corresponding to the vertices of
the hypercube $\paramDomain$.

Because this experiment does not employ POD-based reduced-order models, we do
not require a snapshot set $\paramDomainPOD$.
One set of 100 parameter instances is randomly sampled from the parameter domain
$\paramDomain$ to serve as the set of \pti{} $\paramDomainTrain$.  To assess
method performance for smaller amounts of training data,
we randomly create nested training sets from this set. Another set of 100
parameter instances is
randomly sampled to serve as the set of
parameter testing instances $\paramDomainTest$.

For the inexact-solution case, we consider approximate solutions generated by
a $\napprox=2$ different inexact solutions corresponding to $\itConverge=1$
and $\itConverge=2$ total Newton iterations applied to solve the governing
equations from a linear initial guess 
$\stateIt{0}$ corresponding to
$u(x) = u_a(1 - 2x)$, $x\in[0,1]$. In this case, we consider
	only \ref{data:pooled} (pooled) as described in Section
	\ref{subsection:trainingdata}, 
yielding	
$\ntestData=100$ testing points and up to
$\ntrainingData=100$ testing points.

	For the lower-fidelity-model case, we consider approximate solutions
	generated from $\napprox=2$ different coarse meshes, corresponding to 
$\nstateLF=499$ and $\nstateLF=999$; each of these cases employs
the same discretization as the high-fidelity model and an equispaced grid.  In this case, we consider both \ref{data:pooled} (pooled) and
\ref{data:unique} (unique) as described in Section
\ref{subsection:trainingdata}. The former case yields 
$\ntestData=200$ testing data points and up to
$\ntrainingData=200$ training points, while the latter yields
$\ntestData=100$ testing points and up to
$\ntrainingData=100$ testing points.
When using Feature-Engineering Methods~\ref{feat:paramgappyPODres} and \ref{feat:paramsampleres},
three sample sizes are used: $\nsamples=\{10,\, 100,\,1000\}$.

\def \vartxta {em}
\def \vartexa {$\errorQoiArg{m}$}
\def \nvar    {1}
\def\ntrain/{100}
\FloatBarrier 
\showlocalfalse
\ispcapfalse
\hasgausstrue
\def \casetxt {burgers_unconverged}
\def \casetex {Inexact solutions}
\def\kpercentbetter/{40.48\%}
\def\qpercentbetter/{57.14\%}
\def\extratext/{Instances for which $\text{MSE}_q=\text{MSE}_k$ are not
plotted and comprise 2.38\% of cases.}
\def\globaltex/{\ref{data:pooled}}
\subsubsection{Results: inexact solutions} 
\label{BurgerUnconverged}
Again, we first assess the difference in performance between the $k$- and
$q$-sampling approaches for computing the sampling matrix $\sampleMat$
employed by Feature-Engineering Methods \ref{feat:paramgappyPODres} and \ref{feat:paramsampleres}.
Figure~\ref{fig:\casetxt/kq_histogram} provides a comparison of the test MSEs
that arise from using $k$- and $q$-sampling with $\card{\paramDomainTrain}=\ntrain/$ training runs.
This figure is generated from 42 total data points, which aggregate
test MSE
	values over
Feature-Engineering Methods \ref{feat:paramgappyPODres} and
	\ref{feat:paramsampleres}, the seven regression techniques discussed
	in~\ref{subsection:regressionfunctionapprox}, and three numbers of sample points
	$\nsamples\in\{10, 100, 1000\}$.
Yet again, $q$-sampling significantly outperforms $k$-sampling;  therefore,
the remaining results for this experiment only consider $q$-sampling.


For each of the ML regression techniques, Figure~\ref{fig:\casetxt/best_ind} reports how the 
test FVU varies with respect to the number of 
parameter training instances
$\card{\paramDomainTrain}$.
Consistent with observations in previous experiments, as the number of parameter training instances 
$\card{\paramDomainTrain}$
increases, the test FVU decreases,
and the best regression techniques are those that enable higher capacity:
ANN and SVR: RBF.
In this case, these two regression techniques yield roughly two orders of
magnitude smaller test FVU values than other regression techniques for
$\card{\paramDomainTrain}=\ntrain/$ training runs. The two linear methods OLS:
Linear and SVR: Linear are unable to reduce the test FVU below 0.1 (with any
feature choice), even with a large amount of training data, due to their
relatively low capacity.

Figure~\ref{fig:\casetxt/best_ml} compares each feature-engineering method when using the
best regression technique for that method.
Once again, it is clear that using the feature $\FeatureResNorm$ alone yields the 
highest test FVU, which does not improve as the amount of training data increases.
Generally, using features $\FeatureParResGappyPCA$ yields better performance
than using features $\FeatureParResSampled$.

Figure~\ref{fig:\casetxt/matrix_t\ntrain/} reports the test FVU for each
combination of feature-engineering method and regression technique
for $\card{\paramDomainTrain}=\ntrain/$ \pti{}.
First, we note that SVR: RBF and ANN again consistently yield the best
performance, while OLS: Linear and OLS: Quadratic yield inconsistent performance.
Furthermore, it is apparent that features $\FeatureParResGappyPCA$ yield
significantly better performance than features $\FeatureParResSampled$,
despite the fact they employ the same underlying samples of the residual
vector.
For this experiment, a remarkably small number of residual samples
$\nsamples=10$ is required to yield a coefficient of determination $r^2$ exceeding
0.999; this is achieved for features 
$[\params;\,\resApproxRedGappy(\params)]$
and ANN regression.

For $\card{\paramDomainTrain}=\ntrain/$ \pti{},
Figure~\ref{fig:\casetxt/exact_v_pred} compares the predicted values of $e_m$
with the exact values using the conventional feature choice $\FeatureResNorm$,
as well as features $\FeatureParResGappyPCA$,
each with their best-performing regression technique. 
%
Features $\FeatureParResGappyPCA$ yield drastically better performance than the
conventional choice in this case, as it yields a coefficient of determination
$r^2$ exceeding $0.9999$ and a test MSE that is several orders of magnitude
smaller than that obtained with the conventional feature choice. We emphasize
that the conventional approach in this case yields a poor coefficient of
determination $r^2 = 0.26913$.

\reviewerThree{
Figure~\ref{fig:\casetxt/hist_gauss} assesses the accuracy of the noise model
$\noiseModel\sim\normal{0}{\varianceEstimate}$ computed in
\ref{step:noiseApproximation} in Section \ref{subsec:overview} for $\FeatureParResGappyPCA$
with $\nsamples=10$ using the best-performing regression technique.
As before, 
the figure reports the standard normal distribution compared to a histogram of
the prediction errors, which have been standardized according to the
hypothesized distribution.
These prediction errors are computed 
on a second set of independent test data
$\testDataTwo\defeq\{(\errorTestTwo{i},\featuresTestTwo{i})\}_{i=1}^{\card{\testData}}$,
which is constructed in an identical way to the the first set of test data
$\testData$.
Table~\ref{tab:\casetxt/validation_frequencies} lists the corresponding
validation frequencies $\frequency_{\errorQoiArg{m}}(\frequency)$~\eqref{eq:validation_frequency} for multiple
$\frequency$-prediction intervals~\eqref{eq:prediction_interval}.

\begin{table}[ht!]
\centering
\reviewerThree{
\begin{tabular}{c c c c}
\toprule
$\frequency$ & $\frequency_{\errorQoiArg{m}}(\frequency)$ \\
\midrule
0.80 & 0.830 \\
0.90 & 0.900 \\
0.95 & 0.940 \\
0.99 & 0.970 \\
\bottomrule
\end{tabular}
}
\vspace{-.5em}
\caption{\casetex: Validation frequencies for model using $\FeatureParResGappyPCA$ ($\nsamples=10$) with ANN.}
\label{tab:\casetxt/validation_frequencies}
\end{table}

In this case, the data appear to be relatively close to Gaussian. This is
further evidenced by the accuracy of the 0.90- and 0.95-prediction intervals.
}


\FloatBarrier 
\showlocaltrue
\def \casetxt {burgers_coarse}
\def \casetex {Coarse solution prolongation}
\def\localtex/{\ref{data:unique}}
\def\globaltex/{\ref{data:pooled}}
\def\kpercentbetter/{35.71\%}
\def\qpercentbetter/{58.33\%}
\def\extratext/{Instances for which $\text{MSE}_q=\text{MSE}_k$ are not
plotted and comprise 5.95\% of cases.}
\subsubsection{Results: coarse solution prolongation} 

\label{sec:BurgerCoarse}
We again assess the difference in performance between the $k$- and
$q$-sampling approaches.
Figure~\ref{fig:\casetxt/kq_histogram} provides a comparison of the test MSEs
that arise from using $k$- and $q$-sampling with $\card{\paramDomainTrain}=\ntrain/$
parameter training instances.  
The 84 data points aggregate test MSE values over Feature-Engineering Methods 
\ref{feat:paramgappyPODres} and \ref{feat:paramsampleres},
the seven
regression techniques discussed
in~\ref{subsection:regressionfunctionapprox}, three numbers of sample points
$\nsamples\in\{10, 100, 1000\}$, and the two data-set approaches discussed in Section \ref{subsection:trainingdata}.
Once more, $q$-sampling outperforms $k$-sampling; therefore, the
remaining results for this experiment only consider $q$-sampling.


For each of the ML regression techniques, Figure~\ref{fig:\casetxt/best_ind}
reports how the test FVU varies with respect to the number of parameter
training instances $\card{\paramDomainTrain}$.
Again, the test FVU decreases as the number of parameter training instances 
$\card{\paramDomainTrain}$ increases, and the best regression techniques are
those that enable higher capacity: ANN and SVR: RBF. In addition, OLS:
Quadratic performs relatively well, although its performance is highly
sensitive to the feature choice as shown in Figure
\ref{fig:\casetxt/matrix_t\ntrain/}. 
RF and $k$-NN yield the worst performance, most likely due to an insufficient
amount of training data.


Figure~\ref{fig:\casetxt/best_ml} compares each feature-engineering method when using the
best regression technique for that method.
Yet again, using the feature $\FeatureResNorm$ alone yields the highest test
FVU, which does not improve as the amount of training data increases.
Generally, using features $\FeatureParResGappyPCA$ or $\FeatureParResPCA$
typically yields better performance that using commonly used features
$\FeaturePar$.  Additionally, in this case, employing \ref{data:unique} yields
significantly better performance than employing \ref{data:pooled}. This
implies that the feature--error relationship is quite different between the
two coarse meshes considered. As a result, constructing a regression model
from the union of their respective data results in relatively poor
performance---even with access to twice as much training data---compared with
regression models trained on each coarse-mesh model independently (see
discussion of \ref{data:pooled} in Section \ref{subsection:trainingdata}).


Figure~\ref{fig:\casetxt/matrix_t\ntrain/} reports the test FVU for each
combination of feature-engineering method, regression technique, and data-set
approach for $\card{\paramDomainTrain}=\ntrain/$ \pti{}.
Again, we observe the best performance from regression techniques SVR: RBF and ANN.
This figure again illustrates the benefit to employing \ref{data:unique}
rather than \ref{data:pooled}.
Again, we observe superior performance from 
features $\FeatureParResGappyPCA$ compared with features
$\FeatureParResSampled$, even though they use the same sampled elements of the residual.
As previously observed, feature $\FeatureResNorm$ is nearly insensitive to
regression technique and yields relatively poor performance, as it
corresponds to a single low-quality feature.
We observe that only $\nsamples=10$ samples are needed to achieve a
coefficient of determination $r^2$ exceeding $0.9999$, which is achieved for
$[\params;\,\resApproxRedGappy(\params)]$ and SVR: RBF regression using
\ref{data:unique}.


For $\card{\paramDomainTrain}=\ntrain/$ \pti{},
Figure~\ref{fig:\casetxt/exact_v_pred} compares the predicted values of $e_m$
with the exact values using the conventional feature choices:
$\FeatureResNorm$ and $\FeaturePar$, as well as $\FeatureParResGappyPCA$, where performance of the best
regression technique for each of these feature-engineering methods is
reported.
%
Features $\FeatureParResGappyPCA$ perform better than conventional features, yielding
test MSE values whose values are one to two orders of magnitude lower.

\reviewerThree{

Figure~\ref{fig:\casetxt/hist_gauss} assesses the accuracy of the noise model
$\noiseModel\sim\normal{0}{\varianceEstimate}$ computed in
\ref{step:noiseApproximation} in Section \ref{subsec:overview} for $\FeatureParResGappyPCA$
with $\nsamples=10$ using the best-performing regression technique.
As before, 
the figure reports the standard normal distribution compared to a histogram of
the prediction errors, which have been standardized according to the
hypothesized distribution.
These prediction errors are computed 
on a second set of independent test data
$\testDataTwo\defeq\{(\errorTestTwo{i},\featuresTestTwo{i})\}_{i=1}^{\card{\testData}}$,
which is constructed in an identical way to the the first set of test data
$\testData$.
Table~\ref{tab:\casetxt/validation_frequencies} lists the corresponding
validation frequencies $\frequency_{\errorQoiArg{m}}(\frequency)$~\eqref{eq:validation_frequency} for multiple
$\frequency$-prediction intervals~\eqref{eq:prediction_interval}.

\begin{table}[ht!]
\centering
\reviewerThree{
\begin{tabular}{c c c c}
\toprule
$\frequency$ & $\frequency_{\errorQoiArg{m}}(\frequency)$ \\
\midrule
0.80 & 0.935 \\
0.90 & 0.950 \\
0.95 & 0.965 \\
0.99 & 0.970 \\
\bottomrule
\end{tabular}
}
\vspace{-.5em}
\caption{\casetex: Validation frequencies for model using $\FeatureParResGappyPCA$ ($\nsamples=10$) with SVR: RBF.}
\label{tab:\casetxt/validation_frequencies}
\end{table}

Here, the data appear to be somewhat close to Gaussian, as evidenced by the
accuracy of the 0.95- and 0.99-prediction intervals. Nonetheless, future work
entails construction of more sophisticated noise models to more accurately
capture non-Gaussian or heteroscedastic behavior.
}



\FloatBarrier 
\subsection{Summary} 
We now summarize the dominant trends observed in the numerical experiments
reported in Sections \ref{subsec:cube}--\ref{subsec:burgers}.
\subsubsection{Feature-engineering methods} 

First, Feature-Engineering Method \ref{feat:params}, which employs features
$\FeaturePar$, is the most commonly adopted technique in the literature. The
numerical experiments observe that this technique yields good performance only
if applied with regression methods SVR: RBF or ANN; it performs poorly with
the standard technique OLS: Linear. Furthermore, this feature-engineering method
is always outperformed by one of the new feature-engineering methods proposed
in this work.
Second, we observe that Feature-Engineering Method \ref{feat:DWR}, which
employs dual-weighted-residual features $\dualweightedresidual(\params)$ and
was employed by the ROMES method \cite{drohmann_2015}, yields fairly good
performance for a small training set, as it is a single high-quality feature.
However, it is costly to compute, as it requires a dual solve, and is not
practical for application to large-scale problems.  Third, we observe that
Feature-Engineering Method \ref{feat:resnorm}, which employs features
$\resNorm{\params}$, yields poor performance, as it comprises a single
low-quality feature.  Additionally, its computation incurs an $\nstate$-dependent
operation count.  Ultimately, the best-performing feature-engineering
technique is Feature-Engineering Method \ref{feat:paramgappyPODres}, which is
newly proposed in this work and employs features
$[\params;\,\resApproxRedGappy(\params)]$; this approach yields a coefficient
of determination $r^2$ in excess of  $0.996$ in each experiment.  Moreover, these features can be computed inexpensively, as they require computing only a
small number of residual elements. In many cases, sampling only 10 elements of
the residual is sufficient to produce very accurate regression models with low-noise-variance results, and there is rarely much benefit to using more than
100 residual samples.

We additionally observe that for Feature-Engineering Methods
\ref{feat:paramgappyPODres} and \ref{feat:paramsampleres}, the best approach
for computing the required sampling matrix $\sampleMat$ is the \textit{q-sampling} approach proposed
in the model-reduction community \cite{drmacGugercinReview}; this technique
consistently outperforms \textit{k-sampling}, which is a linear univariate feature
selection approach developed in the machine-learning community. This result is
of interest, as it is one of the first examples of a technique developed in
the model-reduction community improving the performance of standard
machine-learning methods.

\subsubsection{Regression techniques} 
Overall, we observe that low-capacity linear regression methods OLS: Linear
and SVR: Linear tend to yield relatively
poor performance that does not significantly improve as the amount of training
data increases. We also observe that RF and $k$-NN tend to yield the worst
performance, as these high-capacity regression methods likely require more
data to generalize well.  Due to the strong structure imposed by OLS: Linear and
OLS: Quadratic techniques, the performance of these methods is highly dependent
on the choice of features; in particular, Feature-Engineering Method \ref{feat:paramgappyPODres} produces significantly better performance
for these regression techniques than does Feature-Engineering Method
\ref{feat:paramsampleres}, despite their reliance on the same sampled elements
of the residual.
Consistently, regression methods SVR: RBF and ANN yield the best performance.
Overall, the best results are typically obtained with these two regression
techniques and Feature-Engineering Method \ref{feat:paramgappyPODres} deployed
with a relatively small number of residual samples.

\subsubsection{Data-set construction methods} 
Overall, \ref{data:pooled} and \ref{data:unique} lead to similar performance, with \ref{data:unique} performing slightly better.  \ref{data:pooled} provides the advantage of more training data per high-fidelity solution; however, this advantage appears to be undermined by the disparity in error magnitudes across the different approximation solutions (e.g., $\nbasis$, $\itConverge$, or $\nstateLF$).

\section{Conclusions} 
\label{sec:conclusions}

This work has proposed a novel approach for quantifying both
quantity-of-interest errors and normed solution errors in approximate
solutions to parameterized systems of nonlinear equations. The
technique applies machine-learning regression methods (e.g., support vector
regression, artificial neural networks) to map features (i.e.,
error indicators such as residual gappy principal components) to a prediction
of the error; the noise model quantifies the epistemic uncertainty introduced
by the approximate solution and is modeled as a mean-zero, constant-variance
Gaussian random variable whose variance corresponds to the sample variance of
the approximate-solution error on a test set.

Experiments conducted on a range of computational-mechanics problems and types
of approximate solutions demonstrated that the best performance overall was
obtained by Feature-Engineering Method \ref{feat:paramgappyPODres}, which
employs features $[\params;\,\resApproxRedGappy(\params)]$, and regression
techniques corresponding to support vector regression with a Gaussian
radial-basis-function kernel (SVR: RBF) and a feed-forward artificial neural
network (ANN).  It is important to note that, even for problems with a quarter
of a million degrees of freedom, the proposed method is able to generate
inexpensive-to-evaluate machine-learning error models exhibiting coefficients
of determination exceeding $r^2=0.996$; computing the features required for
the error model necessitated evaluating only ten elements of the residual
vector.

Future work involves extending the methodology to dynamical systems,
developing higher-capacity heteroscedastic noise models, and deploying the
methodology on approximate solutions to coupled systems.

\section*{Acknowledgments} 
\label{sec:acknowledgments}
The authors thank Jeffrey Fike for his valuable assistance with Albany. 
This work was sponsored by Sandia's Advanced Simulation and Computing (ASC) Verification and Validation (V\&V) Project/Task \#103723/05.30.02.
This paper describes objective technical results and analysis. Any subjective views or opinions that might be expressed in the paper do not necessarily represent the views of the U.S. Department of Energy or the United States Government
Sandia National Laboratories is a multimission laboratory managed and operated by National Technology and Engineering Solutions of Sandia, LLC., a wholly owned subsidiary of Honeywell International, Inc., for the U.S. Department of Energy's National Nuclear Security Administration under contract DE-NA-0003525.

\appendix
\section{POD algorithm}
\label{app:pod}
Algorithm \ref{alg:pod} reports the POD algorithm considered in this work.
\begin{algorithm}[h!]
\caption{Proper orthogonal decomposition}
\begin{algorithmic}[1]\label{alg:pod}
\REQUIRE 
Snapshots $\{\state^i\}_{i=1}^{\nsnapshots}\subset\RR{\nstate}$
\ENSURE (Full) POD basis $\trialbasis\in\RRstar{\nstate\times\nsnapshots}$; statistical energy captured by different basis
dimensions $\energyThresholdArg{i}$, $i=1,\ldots,\nsnapshots$
\STATE Compute singular value decomposition $[\state^1\ \cdots\ \state^{\nsnapshots}] = \leftSing\Sing\rightSing^T$
with $\leftSing\equiv[\leftSingArg{1}\ \cdots\ \leftSingArg{\nsnapshots}]$, $\rightSing\equiv[\rightSingArg{1}\ \cdots\ \rightSingArg{\nsnapshots}]$, and $\Sing\equiv\diag(\SingArg{i})$ and $\SingArg{1}\geq \cdots\geq\SingArg{\nsnapshots}\geq 0$.
\STATE Determine statistical energy captured by different basis dimensions:
$\energyThresholdArg{i} = 
\sum_{k=1}^i\SingArg{k}^2/\sum_{k=1}^{\nsnapshots}\SingArg{k}^2$, $i=1,\ldots,\nsnapshots$.
\end{algorithmic}
\end{algorithm}

\section{Method algorithms}\label{app:methodalgorithms}
\reviewerThree{
This section provides the algorithms employed by the proposed method. The
offline stage executes Algorithm \ref{alg:dataGen}, followed by Algorithm
\ref{alg:regressionNoiseModel}. The online stage executes Algorithm
\ref{alg:online}.
\begin{algorithm}[h!]
\caption{Offline stage, step 1: data generation}
\begin{algorithmic}[1]
\label{alg:dataGen}
	\reviewerThree{
	\REQUIRE Parameter training instances $\paramDomainTrain\subset\paramDomain$, parameter test instances $\paramDomainTest\subset\paramDomain$, 
number of different types of approximate solutions $\napprox\geq 1$,
error of interest $\error\in\{\errorQoi,\errorState\}$, feature-engineering method (see Section \ref{subsec:indicators}), data-set method (see Section \ref{subsection:trainingdata})
\ENSURE Training data $\trainingData$, test data $\testData$, training data for residual PCA $\trainingDataRes$ (if Feature-Engineering Method
\ref{feat:paramresPCA} or \ref{feat:paramgappyPODres} is selected)
\FOR{$\params\in\paramDomainTrain\cup\paramDomainTest$}
\STATE Compute the state $\state(\params)$ by solving
	Eq.~\eqref{eq:fom_eqn}.
\FOR{$i=1,\ldots,\napprox$}
\STATE Compute the $i$th approximate solution $\stateApproxArg{i}(\params)$.
\STATE Compute the error
$$
\errorArg{i}(\params)=
\left\{
\begin{array}{@{} c l @{}}
\qoiFunc(\state(\params)) - \qoiFunc(\stateApproxArg{i}(\params)), & \text{if}\ \error=\errorQoi
\\[0.5em]
\|\state(\params)
	-\stateApproxArg{i}(\params)\|, & \text{if}\ \error = \errorState
\end{array}\right..
$$
\ENDFOR
\ENDFOR
\IF{Feature-Engineering Method \ref{feat:paramresPCA} or \ref{feat:paramgappyPODres} is selected}
\FOR{$\params\in\paramDomainTrain$}
\FOR{$i=1,\ldots,\napprox$}
\STATE Compute the residual
$\res (\stateApproxArg{i}(\params);\params
)$.
\ENDFOR
\ENDFOR
\STATE Define the training data for residual PCA $\trainingDataRes$
according to Eq.~\eqref{eq:trainingResDataset1} if \ref{data:pooled} is
	selected, or
according to Eq.~\eqref{eq:trainingResDataset2}
 if \ref{data:unique} is selected.
	\STATE Compute the residual-principal-component matrix $\basisresApprox$
	from the training data $\trainingDataRes$.
\ENDIF
\FOR{$\params\in\paramDomainTrain\cup\paramDomainTest$}
\FOR{$i=1,\ldots,\napprox$}
\STATE Compute the features $\featuresArg{i}(\params)$ resulting from the $i$th approximate solution $\stateApproxArg{i}(\params)$.
\ENDFOR
\ENDFOR
\STATE Define the training data $\trainingData$ and test data $\testData$
according to Eq.~\eqref{eq:trainingTestDataset1} if \ref{data:pooled} is
	selected, or
according to Eq.~\eqref{eq:trainingTestDataset2}
 if \ref{data:unique} is selected.
	}
\end{algorithmic}
\end{algorithm}

\begin{algorithm}[h!]
\caption{Offline stage, step 2: construction of regression and noise models}
\begin{algorithmic}[1]
\label{alg:regressionNoiseModel}
	\reviewerThree{
\REQUIRE Training data $\trainingData$, test data $\testData\defeq\{(\errorTest{i},\featuresTest{i})\}_{i=1}^{\ntestData}$, number of
	cross-validation folds $k$, regression model (see Section
	\ref{subsection:regressionfunctionapprox}), regression-model hyperparameters
	$\hyperparams$, cross-validation grid $\hyperparamsSet$
	\ENSURE Regression model $\regressionModel$, noise-model variance $\varianceEstimate$
\STATE Randomly divide the training data into $k$ non-overlapping sets
	$\trainingDataArg{i}\defeq\{(\errorArgs{i}{\ell},\featuresArgs{i}{\ell})\}_\ell\subset\trainingData$, $i=1,\ldots,k$ such
	that $\trainingDataArg{i}\cap\trainingDataArg{j} = \emptyset$ for $i\neq j$
	and $\cup_{i=1}^k\trainingDataArg{i} = \trainingData$.
\FOR{$i=1,\ldots,k$}\label{cv_loop}
\FOR{$\hyperparams\in\hyperparamsSet$}
	\STATE Compute the regression model $\regressionModelArgs{i}{\hyperparams}$
	using training set $\cup_{j\in\{1,\ldots,k\}\setminus i}\trainingDataArg{j}$ and hyperparameters
	$\hyperparams$.
	\STATE Compute the mean-squared-error cross-validation loss $\validationLossArgs{i}{\hyperparams} =
	\frac{1}{\card{\trainingDataArg{i}}}\sum_{\ell=1}^{\card{\trainingDataArg{i}}}(\errorArgs{i}{\ell}-\regressionModelArgs{i}{\hyperparams}(\featuresArgs{i}{\ell}))^2$.
\ENDFOR
\ENDFOR
	\STATE Select hyperparameters $\hyperparamsSelect =
	\arg\min_{\hyperparams\in\hyperparamsSet}\sum_{i=1}^k\validationLossArgs{i}{\hyperparams}$.
	\STATE Compute regression model $\regressionModel$ using the full training
	set $\trainingData$ and hyperparameters $\hyperparamsSelect$.\label{compute_regression}
	\STATE Compute the noise-model variance as the mean-squared-error test loss
	$\varianceEstimate =
\frac{1}{\ntestData}\sum_{i=1}^{\ntestData}
(\errorTest{i}-\regressionModel(\featuresTest{i}))^2$.
}
\end{algorithmic}
\end{algorithm}

\clearpage

\begin{algorithm}
\caption{Online stage}
\begin{algorithmic}[1]
\label{alg:online}
	\reviewerThree{
\REQUIRE Online parameter instance $\params\in\paramDomain$, regression model $\regressionModel$,
	noise-model variance $\varianceEstimate$
	\ENSURE Approximate solution $\stateApprox(\params)$, quantity-of-interest
	approximation $\qoiApprox(\params)$, error-model prediction $\errorModel(\params)$
\STATE Compute the approximate solution $\stateApprox(\params)$,
	quantity-of-interest approximation $\qoiApprox(\params) \defeq\qoiFunc(\stateApprox(\params))$,
	and resulting features $\features(\params)$.
	\STATE Compute the error-model prediction
$\errorModel(\params) = \regressionModel(\features(\params)) +
	\noiseModel$ with
$\noiseModel\sim\normal{0}{\varianceEstimate}$.
}
\end{algorithmic}
\end{algorithm}

}

\addcontentsline{toc}{section}{\refname}
\bibliographystyle{elsarticle-num}
\bibliography{Brian_references}

\end{document}